\documentclass [12pt] {report}
\usepackage[dvips]{graphicx}
\usepackage{epsfig, latexsym,graphicx}
\usepackage{amssymb}
\newcommand {\eqdef} {\stackrel{\rm def}{=}}
\newcommand {\D}[2] {\displaystyle\frac{\partial{#1}}{\partial{#2}}}

\newcommand {\Dd}[3] {\displaystyle\frac{\partial^2{#1}}{\partial{#2}\partial{#3}}}
\newcommand {\al} {\alpha}

\newcommand {\ga} {\gamma}
\newcommand {\la} {\lambda}

\newcommand {\de} {\delta}

\newcommand {\fr} {\displaystyle\frac}
\newcommand {\wt} {\widetilde}
\newcommand {\be} {\begin{equation}}
\newcommand {\ee} {\end{equation}}
\newcommand {\ba} {\begin{array}}
\newcommand {\ea} {\end{array}}
\newcommand {\bp} {\begin{picture}}
\newcommand {\ep} {\end{picture}}
\newcommand {\bc} {\begin{center}}
\newcommand {\ec} {\end{center}}
\newcommand {\bt} {\begin{tabular}}
\newcommand {\et} {\end{tabular}}
\newcommand {\lf} {\left}
\newcommand {\rg} {\right}

\newcommand {\cC} {{\cal C}}
\newcommand {\cF} {{\cal F}}

\newcommand {\cR} {{\cal R}}
\newcommand {\cS} {{\cal S}}
\newcommand {\cT} {{\cal T}}

\newcommand {\ses} {\medskip}

\newcommand {\e} {\mathop{\rm e}\nolimits}

\newcommand {\cE} {{\cal E}}

\newcommand {\bibit} {\bibitem}
\newcommand {\nin} {\noindent}

\newcommand {\De} {\Delta}

\newcommand {\Rho} {\mbox{\large$\rho$}}

\newcommand {\wh} {\widehat}

\newcommand {\const} {\mathop{\rm const}\nolimits}

\newcommand {\bfC}{{\bf C}}

 \usepackage{amsmath}

\newcommand {\cD} {{\cal D}}

\def\2#1#2#3{{#1}_{#2}\hspace{0pt}^{#3}}
\def\3#1#2#3#4{{#1}_{#2}\hspace{0pt}^{#3}\hspace{0pt}_{#4}}
\newcounter{sctn}
\def\sec#1.#2\par{\setcounter{sctn}{#1}\setcounter{equation}{0}
                  \noindent{\bf\boldmath#1.#2}\bigskip\par}

\begin {document}

\begin {titlepage}

\vspace{0.1in}

\begin{center}

\ses

\ses

{\large \bf Finsler connection preserving the two-vector  angle  under

\ses

\ses

the  indicatrix-inhomogeneous treatment}

\ses

\ses

\end{center}

\vspace{0.3in}

\begin{center}

\vspace{.15in} {\large G.S. Asanov\\} \vspace{.25in}
{\it Division of Theoretical Physics, Moscow State University\\
119992 Moscow, Russia\\
{\rm (}e-mail: asanov@newmail.ru{\rm )}} \vspace{.05in}

\end{center}

\begin{abstract}

 The
Finsler spaces in which the tangent Riemannian spaces are conformally flat
prove to be  characterized by
 the condition
that   the indicatrix  is a space of constant curvature.
In such spaces
the  Finslerian normalized
two-vector angle can  be explicated
from the  respective two-vector angle of the associated  Riemannian space.
Therefore the way is opening  to propose explicitly
the  connection preserving the  angle
even
at  the indicatrix-inhomogeneous level, that is,
when
 the indicatrix curvature value
 $  {\mathcal  C}_{\text{Ind.}}   $
is permitted
 to be  an arbitrary smooth function of the
indicatrix position point $x$.
The  connection obtained is metrical with the deflection part
 which is proportional to the gradient of the function
$H(x)$
entering the equality
 $
{\mathcal  C}_{\text{Ind.}} \equiv H^2.$
Also the connection
is
covariant-constant.
When   the transitivity of covariant derivative is used,
from the commutators of covariant derivatives the associated
curvature tensor is found. Various useful representations  have been  developed.
The Finsleroid space has been explicitly outlined.

\end{abstract}

\end{titlepage}

\vskip 1cm

\ses

\ses

\bc   {\bf Motivation and   Introduction}  \ec

\ses

\ses

In the Finsler space
the tangent bundle $TM$ over the base manifold $M$
is geometrized by means of the
 Finsler metric function $F(x,y)$,
such that
at each point $x\in M$
the tangent vectors $y\in T_xM$
are used, where   $T_xM$ is
the tangent space   supported by the $x$.

The embedded position  of the indicatrix
${\cal I}_x\subset T_xM$ in the tangent Riemannian space
$\cR_{\{x\}}=\{T_xM, g_{\{x\}}(y)\}$
(where $g_{\{x\}}(y)$ denotes the Finslerian metric tensor with $x$ considered fixed
and $y$ used as being the variable) induces the Riemannian metric on the indicatrix
through  the well-known method  (see, e.g., Section 5.8 in [1])
and in this sense makes the indicatrix
 a Riemannian space.
Therefore, the geodesics can be introduced on the indicatrix by applying
the conventional Riemannian methods.

In any (sufficiently smooth) Finsler space
 the two-vector angle
 $ \al_{\{x\}}(y_1,y_2)$
 can locally be determined  with the help of the indicatrix geodesic arc,
 which provokes the important question
 whether the Finsler geometry can
 be profoundly  settled down
 by developing and applying the connection
 which preserves the angle.

\ses

In general,
the angle
 $ \al_{\{x\}}(y_1,y_2)$
  is complicated and cannot  be determined in an explicit form,
except for rare Finsler metric functions.
The
lucky example
is given
by
the Finsler space
$\cF^N$
which is characterized by the condition that the indicatrix
is a space of constant curvature.
In the space
$\cF^N$
the angle
 $ \al_{\{x\}}(y_1,y_2)$
 can be found in the explicit and simple form.
 Namely,
it is possible to prove that (under attractive conditions)
in any dimension
$N\ge3$,
 the  tangent Riemannian  space
$\cR_{\{x\}}$
is conformally flat if and only if the  indicatrix
is a space of constant curvature.
The respective transformation
$  y={\bf C}(x,\bar y)$
is positively homogeneous of the degree which we shall denote by $H$.
 The remarkable equality
 ${\mathcal  C}_{\text{Ind.}} \equiv H^2$
arises,
where
  ${\mathcal  C}_{\text{Ind.}}={\mathcal  C}_{\text{Ind.}}(x)$
denotes  the value of the curvature of the indicatrix
${\cal I}_x\subset T_xM$.

{

Under this transformation
$  y={\bf C}(x,\bar y)$
each tangent Riemannian space
$\cR_{\{x\}}$
is conformally changed  to
become
a
 Euclidean space
$\cE_{\{x\}}$.
The distribution of the last spaces
$\cE_{\{x\}}$
over the base manifold $M$
composes the
{\it associated Riemannian space},
which we denote by
$\cR^N=(M,S)$,
where
$ S=\sqrt{a_{mn}(x) y^m y^n} $
is the Riemannian metric constructed from
the  metric tensor
$a_{mn}(x)$
of the space
$\cE_{\{x\}}$.
We are entitled to
 induce the  angle
 $\al^{\text{Riem}}_{\{x\}}$
 conventionally defined in
the Riemannian
space
$\cR^N$
into
the Finsler
space
$\cF^N$,
obtaining
 simply
$
  \al_{\{x\}}(y_1,y_2)=\bigl(1/H(x)\bigr)  \al^{\text{Riem}}_{\{x\}}(\bar y_1,\bar y_2).
$

\ses

To explicate  the coefficients $N^m{}_n$ of nonlinear connection
from the  Finsler angle
$\al=\al_{\{x\}}(y_1,y_2)$,
we should successfully  propose the preservation equation.
The nearest possibility
is to formulate  the  equation
 $d_i\al=0$
 in accordance with the formulas (I.1.12) and (I.1.15),
applying the separable operator $d_i$ indicated in (I.1.11).

This possibility has been realized in the preceding work [10,11].
Namely, in that work
the separable preservation equation
 $d_i\al=0$
has been  solved in the Finsler space
$\cF^N$
under the assumption that
 ${\mathcal  C}_{\text{Ind.}} =\const$,
 and whence
 $H  =\const$.
  The explicit
 coefficients $N^m{}_n$
 have been obtained.

In general the indicatrix curvature value
${\mathcal  C}_{\text{Ind.}} $
may depend on the points $x\in M$ which support the indicatrix.
We call the space   ${\mathcal F}^N$
 indicatrix-homogeneous,
  if the value is a constant,
  whence
   $H  =\const$.
If the dependence
${\mathcal  C}_{\text{Ind.}} ={\mathcal  C}_{\text{Ind.}} (x)$
does hold,
we say that
the space   ${\mathcal F}^N$
is  indicatrix-inhomogeneous,
in which case $H_i\ne0$, where
$H_i=\partial H/\partial x^i$.
 The representations obtained in the previous work [10,11]
 are the $(H_i\to0)$-limits
of their generalized counterparts developed in the present study.

 {

It appears  that in general the angle preservation equation  formulated in the separable way
does not permit any solution for the coefficients
$N^m{}_n$.

This conclusion can be drawn from the implications
which are derivable
by the help
of
 the  coincidence-limit  method
(see Section 3.2 in [12])
which extracts
 the information from behavior  of Riemannian geodesics.
To this end we should use the distance function $E=E(x,y_1,y_2)$
with
$
E=(1/2)\al^2.
$
Namely,
evaluating various partial derivatives of the
function
with respect to
$y_1$ and $y_2$
and finding   the coincidence limits
when
${y_2}\to{y_1}$,
we can obtain a valuable information on
the derivatives of the Finsler metric tensor of
the Finsler space.
Performing the required evaluations
on the level of
the second-order partial  derivatives
$\partial^2/\partial y_1^m\partial y_2^n$,
and, then,
 applying   the operation
${y_2}\to{y_1}$
to the resultant expressions,
it is possible to arrive at the following general conclusion:
{\it In any  Finsler space
the vanishing assumption
 $d_i\al=0$
of the separable type
 entails
the equality
}
$$
\cD_ih_{mn}
=
\fr 2{F}
h_{mn}     d_iF.
$$

If we additionally postulate
 $d_iF=0,$
 we obtain
 $\cD_ih_{mn}=0$ and, therefore,
the  metricity $\cD_ig_{mn}=0$
which is formulated with the covariant derivative
$\cD$ arisen from the deflectionless connection.
We can apply
the derivative
$\partial^2/\partial y^l\partial y^h$
to the equality
$\cD_ig_{mn}=0$,
which leads after simple evaluations
to the vanishing
$
\cD_i
S_n{}^k{}_{jm}
=0.
$
Here,
the
$
S_n{}^k{}_{jm}
$
is the tensor
which describes the curvature of indicatrix
(see Section 5.8 in [1]).

\ses

Clearly,
the vanishing
$
\cD_i
S_n{}^k{}_{jm}
=0
$
can be realized in but  rare particular cases of the Finsler space.
The vanishing
is realized
in
the indicatrix-homogeneous case of
the space   ${\mathcal F}^N$,
 and cannot be realized in
the indicatrix-inhomogeneous  space   ${\mathcal F}^N$.

Therefore,
the   account for the dependence
$H=H(x)$
in
the
Finsler space
${\mathcal F}^N$
 is neither  straightforward
 nor
 trivial   task.

{

These important
(and rather unexpected?)
implications
enforce us to look for more capable ideas to formulate the preservation of angle.
The  attractive idea is to substitute
the normalized angle
$  {\al}^{\{H(x)\}}_{\{x\}}(y_1,y_2)
=
H(x)
{ \al}_{\{x\}}(y_1,y_2)$
(see (I.1.26))
 with the initial angle
${\al}_{\{x\}}(y_1,y_2)$
in the separable preservation law,
according to (I.1.27).
The law obtained is of  the recurrent-type
(I.1.28),
namely
$
d_i \al+ (1/H)H_i\al=0.
$
It appears that
this  preservation
is reconciled
with
the indicatrix-inhomogeneous
Finsler space
${\mathcal F}^N$
 at any
 scalar $H=H(x)$.
The reason thereto
is the following assertion obtainable by the help of the
 coincidence-limit  method:
{\it
In any Finsler space
the
vanishing assumption
$
d_i \al+ (1/H)H_i\al=0
$
 entails the equality}
$$
\cD_ih_{mn}
=
\fr 2{F}
   h_{mn}
   d_iF
-\fr2HH_i h_{mn}.
$$
\ses\\
In these patterns the vanishing  $d_iF=0$ yields
the equality
$\cD_i g_{mn}=-(2/H)H_i h_{mn}$,
which
entails the extension
of
the previous vanishing
$
\cD_i
S_n{}^k{}_{jm}
=0
$
such that the right-hand part of the extension
is just the expression which is obtained when
the characteristic representation of the tensor
$S_n{}^k{}_{jm}$
of the space
$\cF^N$
under study is inserted under the action of the covariant derivative
$\cD_i.$

Thus,
the recurrent-type equation
$
d_i \al+ (1/H)H_i\al=0
$
of preservation
 of the angle
is reconciled
with
the indicatrix-inhomogeneous
Finsler space
${\mathcal F}^N$
 at any
 scalar $H=H(x)$
 (see  Proposition I.1.2 in Section I.1),
and therefore is accepted in the present work to apply.
We solve the equation with respect to the coefficients
 $N^m{}_n(x,y)$.

\ses

{

The $N^m{}_n(x,y)$ thus appeared to read (I.2.15)
can naturally be interpreted as the
{\it coefficients of the   non-linear connection produced by the  angle
in the space}
${\mathcal F}^N$
{\it
studied  on the
general
indicatrix-inhomogeneous level.}

Because of the conformal flatness
of
the  tangent Riemannian spaces
$\cR_{\{x\}}$,
the Finsler space
${\mathcal F}^N$
involves
 the
 associated Riemannian space
$\cR^N$
and, therefore,
the Riemannian connection
coefficients
$     L^m{}_{ij}=a^m{}_{ij}+S^m{}_{ij}$
 (shown  in (I.1.14))
in which
the entered  Christoffel symbols
$a^m{}_{ij}$
are to be constructed from the Riemannian metric tensor $a_{mn}(x)$
of the   space $\cR^N$;
the notation
$S^m{}_{ij}$
is the torsion tensor.

\ses

With the knowledge of the coefficients
$N^m{}_n(x,y)$,
 we can
straightforwardly evaluate the derivative coefficients
$N^k{}_{im}$
and
 express the Finslerian connection coefficients
$T^k{}_{im}$
through
the Riemannian connection coefficients
$  L^m{}_{ij}=  L^m{}_{ij}(x)$
and the function $H=H(x)$
(by the help of
the formulas (I.1.33) and (I.2.18)).

The  coefficients
$
T^k{}_{im}
$
 involve
 the  deflection tensor
$
\De^k{}_{im} =    -  N^k{}_{im} - T^k{}_{im}$
which is non-vanishing
as far as $H_i\ne0$, namely
$
\De^k{}_{im}
=
(1/H)
H_ih^k_m.
$
There arises the   covariant derivative
$\cT$,
which properties are  listed in (I.1.37)-(I.1.40).
In distinction from the connection
developed in the indicatrix-homogeneous case,
the $\cT$-connection obtained  is no  more deflectionless.
Nevertheless,
the $\cT$-connection
is metrical
and the equality
$
N^m{}_j=- T^m{}_{ji}y^i
  $
holds.

{

In this way,
the
 metrical  non-linear Finsler   connection
${\mathcal  FN}=\{N^m{}_i,T^m{}_{ij}\}$
is  induced  in the
 space    ${ \mathcal   F}^N$ from
the  metrical  linear   connection
$
{\mathcal  RL}=\{L^m{}_j,L^m{}_{ij}\}
$
evidenced in the
Riemannian
space ${\mathcal  R}^N$,
where
$
  L^m{}_j=- L^m{}_{ji}y^i.
  $
The involved function $H=H(x)$ may be an arbitrary smooth function of $x$.

\ses

The Finsler
    connection
$
{\mathcal  FN}=\{N^m{}_i,T^m{}_{ij}\}
$
can be understood
  to be a result of an appropriate nonlinear
deformation
of
the
connection
$
{\mathcal  RL}=\{L^m{}_j,L^m{}_{ij}\}.
$
It is the
 transformation $  y={\bf C}(x,\bar y)$
that represents  the deformation said.

\ses

 In other words,
in    the  Finsler space
${\mathcal  F}^N$
we evidence the phenomenon
that  the metrical non-linear angle-preserving connection
 is  the $\bfC$-deformation of
 the
  metrical  linear   connection
applicable in
the  Riemannian
space
 ${\mathcal  R}^N$:
$$
{\mathcal  FN}=\bfC\cdot {\mathcal  RL}.
$$

We shall show that
the  $\bfC$-deformation  is
$\cT$-{covariant constant:
$
{\cT}\cdot{\bfC}=0.
$
Also, the covariant derivative $\cT$
is  the
manifestation of the
{\it transitivity}
of the connection  under   this transformation,
in short,
$
\cT=\cC\cdot\nabla,
$
where
$\nabla$
is
 the
  covariant derivative
 applicable
in
the   Riemannian space
 ${\mathcal  R}^N$.

\ses

In the Riemannian geometry we have merely $H=1$.
In the Finsler space
$\cF^N$,
the scalar $H(x)$
plays the role of the parameter
which changes the indicatrix curvature value.
Varying the scalar $H(x)$ evokes  the changes in the Finsler space
$\cF^N$.

{

In the theory of Finsler spaces
 the notion of connection
was  studied
on the basis of various convenient sets of axioms
(see [1-5] and references therein).
Regarding the significance of the angle notion,
the important  step was made in [6]
were in processes of studying  implications of the
two-vector angle defined by area,
 the theorem was proved
 which states that a diffeomorphism
between two Finsler spaces is an isometry iff it keeps the  angle.
 This  Tam{\'a}ssy's  theorem
clearly
substantiates
the  idea  to develop the Finsler connection from the
Finsler two-vector angle, possibly on the analogy of  the Riemannian geometry.

To meet new methods of applications,
the interesting chain of linear connections
was introduced and studied in [3].
It was emphasized that in the Riemannian geometry we have naturally
the metrical and linear connection
applicable on the tangent bundle
of the variables $x,y$.
Like to the
constructions developed in the preceding work [10,11] dealt with
the indicatrix-homogeneous case,
in the present indicatrix-inhomogeneous study
of the space $\cF^N$
the export of  this connection
generates  the   required Finsler connection.

{

By performing the comparison  between the commutators
of the obtained Finsler  covariant derivative $\cT$ and
 the commutators of  the underlined
Riemannian
covariant derivative $\nabla$,
 {\it not} assuming $H=\const$ so  that
 $H(x)$ is permitted to be an arbitrary smooth function of $x$,
the associated curvature tensor
 $\Rho_k{}^n{}_{ij}$
 can  straightforwardly  be derived.

\ses

The Finsleroid case of the  space ${\cF}^N$
provides us with the example
when the
key transformation
$  y={\bf C}(x,\bar y)$
is known explicitly.
Therefore, we can straightforwardly
 apply
 the developed
  indicatrix-inhomogeneous
  theory
taking the metric function of the Finsleroid type.
The explicit representations
 for the respective Finsleroid  coefficients
$ N^m{}_n$, as well as for
 the entailed derivative coefficients
$ N^k{}_{im}$  and  $ N^k{}_{imn}$,
are found.
Thus we have got prepared the connection
${\mathcal  FN}$
in the Finsleroid space at our disposal
with
 an arbitrary input scalar $H(x)$.

\ses

Below we are interested  in spaces of the dimension $N\ge3$.
The two-dimensional case has been studied in
[8,9].

 {

\ses

\ses

\ses

\ses

{\large \bf Chapter I. The Outline}

\ses

\ses

\ses

\ses

\setcounter{equation}{0}

 {\bf I.1.  Basic representations}

\bigskip

For a given  function Finsler metric function $F=F(x,y)$
we  can construct
the covariant tangent vector $\hat y=\{y_i\}$
and the  Finslerian metric tensor $\{g_{ij}\}$
in the conventional  way:
$
y_i :=(1/2)\partial {F^2}/\partial{y^i}
$
and
$
g_{ij} :
=\partial  y_i/\partial  y^j.
$
The contravariant tensor $\{g^{ij}\}$ is
defined by the reciprocity conditions
$g_{ij}g^{jk}=\de^k_i$,
where $\de$ stands for the Kronecker symbol.
The indices $i,j,\dots$ refer to local admissible coordinates $\{x^i\}$ on the base
manifold $M$.
We shall  also use
 the tensor
$
C_{ijk}=(1/2)\partial g_{ij}/\partial y^k$.
By $l$ we shall denote the unit vectors, namely, $l=y/F(x,y)$,
such that $F(x,l)=1$.

\ses

Let $ U_x $ be a simply connected and geodesically complete region on the indicatrix
${\cal I}_x$ supported by a point $x\in M$.
Any point pair $u_1,u_2\in U_x$ can be joined by the respective
arc ${\cal A}_{\{x\}}(l_1,l_2) \subset {\cal I}_x$
of the  Riemannian  geodesic  line
drawn on $U_x$.
By identifying the  length of the arc  with the angle notion we arrive at
the {\bf geodesic-arc angle}
$ \al_{\{x\}}(y_1,y_2)$,
where $y_1,y_2\in T_xM$  are  two vectors issuing from the origin $0\in T_xM$ and
possessing the property that their direction rays $0y_1$ and $0y_2$
intersect
 the indicatrix at the  point pair $u_1,u_2\in U_x$.
We obtain
\be
\al_{\{x\}}(y_1,y_2)=||{\cal A}_{\{x\}}(l_1,l_2) ||.
\ee

\ses

The   coefficients
$N^k{}_i =N^k{}_i (x,y)$
are
required to construct the operator
\be
d_i~:=\D{}{x^i}
+
N^k{}_{i}
\D{}{y^k}.
\ee
These coefficients are assumed naturally to be positively homogeneous of degree 1
with respect to the vector argument $y$.

{

The  derivative coefficients
\be
N^k{}_{nm}=
\D{N^k{}_{n} }{y^m},  \qquad
N^k{}_{nmj}=
\D{N^k{}_{nm} }{y^j}
\ee
possess the identities
\be
N^k{}_{nm}y^m=N^k{}_{n}, \qquad N^k{}_{nmj}y^m= N^k{}_{nmj}y^j=0,
\qquad N^k{}_{nmj}= N^k{}_{njm}.
\ee

The coefficients are used to construct
  the covariant derivatives
\be
\cD_kF~:=d_kF, \qquad
\cD_k l^m~:=d_kl^m-N^m{}_{kn}l^n, \qquad
\cD_k l_m~:=d_kl_m+N^h{}_{km}l_h,
\ee
and
\ses\\
\be
\cD_k g_{mn}~:= d_kg_{mn} +N^h{}_{km}g_{hn} + N^h{}_{kn}g_{mh}.
\ee
\ses\\
The identities
\be
\D{\cD_k F}{y^m}
=\cD_k l_m,
\qquad
\D{\cD_k l_m}{y^n}
=\cD_k g_{mn}+l_hN^h{}_{kmn}
\ee
are obviously valid,
together with
\be
\D{\cD_i g_{mn}}{y^j}
=
2\cD_iC_{mnj}
+
N^t{}_{imj}g_{tn}
+
N^t{}_{inj}g_{mt},
\ee
\ses\\
where
\be
{\cD}_iC_{mnj}~:=
d_iC_{mnj}
+
N^t{}_{ij}C_{mnt}
+
N^t{}_{im}C_{tnj}
+
N^t{}_{in}C_{mtj}.
\ee

In addition to the Finsler metric tensor
$g_{mn}$, we shall use also the tensor
\be
h_{mn} =g_{mn} -l_ml_n
\ee
\ses\\
which possesses the property
$
h_{mn}y^n=0.
$
The covariant derivative of this tensor will be constructed in the manner similar to (1.6),
namely
$$
\cD_k h_{mn}~:= d_kh_{mn} +N^h{}_{km}h_{hn} + N^h{}_{kn}h_{mh}.
$$

{

To    deal with  the two-vector angle
$\al= \al_{\{x\}}(y_1,y_2)$,
we
  merely
  extend the operator
 $d_i$ in the
 {\it separable way}, namely
\ses
\be
d_i=\D{}{x^i}
+N^k{}_{i}(x,y_1)
\D{}{y^k_1}
+
N^k{}_{i}(x,y_2)
\D{}{y^k_2},         \qquad y_1,y_2\in T_xM,
\ee
and introduce the covariant derivative $\cD_i \al$
according to
\be
\cD_i \al=d_i \al.
\ee

\ses

In the Riemannian geometry we have the  separable  operator
\be
d^{\text{Riem}}_i=\D{}{x^i}
+L^k{}_i{(x,y_1)} \fr{\partial}{\partial y_1^k}
+L^k{}_i{(x,y_2)} \fr{\partial}{\partial y_2^k}, \qquad y_1,y_2\in T_xM,
\ee
with
the linear coefficients
$
  L^k{}_i{(x,y_1)} =-L^k{}_{ij}(x)y_1^j
$
and
$    L^k{}_i{(x,y_2)} =-L^k{}_{ij}(x)y_2^j
$
obtained from the Riemannian connection coefficients
\be
   L^m{}_{ij}=a^m{}_{ij}+S^m{}_{ij},
\ee
where  $a^m{}_{ij}=a^m{}_{ij}(x)$ stands for the Christoffel symbols
constructed from the Riemannian metric tensor $a_{mn}(x)$
and
 $S^m{}_{ij}=S^m{}_{ij}(x)$ is an arbitrary {\it torsion tensor}:
$S_{mij}=-S_{jim}$
with
$S_{mij}=a_{mh}S^h{}_{ij}$.
When applied to the Riemannian
 two-vector angle
$\al^{\text{Riem}}_{\{x\}}(y_1,y_2)=a_{mn}(x)y^m_1y^n_2/S_1S_2,$
where $S_1=\sqrt{a_{mn}(x)y^m_1y^n_1}$ and $S_2=\sqrt{a_{mn}(x)y^m_2y^n_2}$,
the operator reveals the  fundamental vanishing property
$$
d^{\text{Riem}}_i   \al^{\text{Riem}}_{\{x\}}(y_1,y_2)=0,
\qquad y_1,y_2\in T_xM.
$$

{

By analogy,
 one may
 assume that
the  Finsler coefficients
$N^k{}_{i}$
fulfill
the {\it separable angle-preservation equation}
\ses\\
\be
\cD_i \al=0
\ee
\ses\\
 to try developing the theory in which
the  properties
\be
\cD_kF=0, \qquad
\cD_k l_m=0,
\qquad \cD_k l^m=0,
\ee
\ses\\
together with the {metricity}
\be
\cD_k g_{mn}=0
\ee
hold fine.
This  metricity, taken in conjunction with
the identities
indicated in (1.7), just entails
the vanishing
\be
l_h
N^h{}_{kmn}
=
0.
\ee

{

The following  valuable implication
can be deduced from angle
by applying the  coincidence-limit  method
exposed in  Section 3.2 in [12]:
{\it
In any  Finsler space
the vanishing assumption
 $d_i\al=0$
of the separable type
 entails
the equality
\ses\\
\be
\cD_ih_{mn}
=
\fr 2{F}
h_{mn}     d_iF
\ee
}
\ses\\
(take below  the formula (1.29),
keeping $H=\const$).
If we additionally postulate
 $d_iF=0,$
 we obtain
 $\cD_ih_{mn}=0$ and, therefore,
 $\cD_ig_{mn}=0$ .

\ses

\ses

Thus, starting with
the  separable angle-preservation equation
leads to   the following implication:

\ses

\ses

\be
 PRESERVATION  ~ ~ OF ~ ~ ANGLE ~~ AND ~ ~ LENGTH  ~ ~ ~ \Longrightarrow ~ ~ ~  METRICITY,
 \ee

\ses

\ses

\nin
that is, the two conditions $d_i\al=0$ and $d_iF=0$ entail $\cD_i g_{mn}=0$.

{

\ses

\ses

\ses

When
$
{\cD}_i g_{mn}
=0,
$
  the identity (1.8) communicates the
validity of the vanishing
\be
2
\cD_iC_{mnj}
+
N^t{}_{imj}g_{tn}
+
N^t{}_{inj}g_{mt}
=0,
\ee
which in turn entails that,
because the tensor
$C_{mnj}$
is totally symmetric,
 the tensor
$$
N_{nimj}
~:=N^t{}_{imj}g_{tn}
$$
must be  totally symmetric
 with respect to
the subscripts
$n,m,j$:
\be
N_{nimj}= N_{minj} = N_{jimn}= N_{nijm},
\ee
\ses\\
and whence
\be
N^k{}_{imn}
=
-
\cD_i
C^k{}_{mn},
\ee
\ses\\
where
$$
{\cD}_iC^k{}_{mn}~:=
d_iC^k{}_{mn}
-
N^k{}_{it}C^t{}_{mn}
+
N^t{}_{im}C^k{}_{tn}
+
N^t{}_{in}C^k{}_{mt}.
$$
With the representation (1.23),
the vanishing (1.18) can be regarded as
 a direct implication of the identity
$y^kC_{knj}=0$ shown by the    tensor  $C_{knj}$.

\ses

\ses

Thus,
{\it in any  Finsler space  the two  conditions     }
 $d_i\al=0$
 {\it and }
 $d_iF=0$
{\it entail the representation }
(1.23)
{\it for the coefficients   }
$N^k{}_{imn}$.

\ses

\ses

By differentiating these coefficients with respect to  $y^j$
and making the interchange of the indices
$m,j$,
and also noting that
$
\partial N^k{}_{imn}/\partial y^j-\partial N^k{}_{ijn}/\partial y^m=0
$
and
$$
\D{C^k{}_{mn}}{y^j}  -  \D{C^k{}_{jn}}{y^m}
=
-
2
\lf( C^h{}_{nm}    C^k{}_{hj} - C^h{}_{nj}    C^k{}_{hm} \rg),
$$
from
(1.23)
we can arrive at the following vanishing after a short evaluation:
\be
\cD_i
S_n{}^k{}_{jm}
=0,   \quad \text{where} ~~
S_n{}^k{}_{ij}
=
\lf(
C^h{}_{nj}    C^k{}_{hi}
-
C^h{}_{ni}    C^k{}_{hj}
\rg)
F^2.
\ee

\ses

{

However, there are no  reasons to trust
that the separable form (1.15) for the angle preservation
 is applicable in general to any Finsler space.
For it might happen
that the equation (1.15)
doesn't  permit any solution
with respect to
 the coefficients
$N^k{}_{i}=N^k{}_{i}(x,y)$.
Indeed, the formula (1.24) tells us that
the following proposition is true.

\ses

\ses

{\bf Proposition I.1.1.}
{\it
One is entitled to hope to determine the coefficients
$N^m{}_n$ of a Finsler space
 from the separable equation
 $d_i\al=0$ supplemented by the condition $d_iF=0$
if only the Finsler space possesses the property
}
$
\cD_i
S_n{}^k{}_{jm}
=0.
$

\ses

\ses

Clearly,
the vanishing
$
\cD_i
S_n{}^k{}_{jm}
=0
$
can be realized in but  rare particular cases of the Finsler space.

\ses

{

In this connection it can be of help to introduce
a {\it characteristic  indicatrix scale} $R(x)$
in each tangent space to normalize the angle.
If the volume
$V_{{\cal I}_x}$
of
the Finslerian indicatrix
$
{\cal I}_x\subset T_xM
$
is finite, it is attractive to obtain the scale
by the help of the equality
\be
V_{{\cal I}_x}
=
C_1 (R(x))^{N-1}, \qquad C_1=\const.
\ee
In this case the $R(x)$ has the geometrical meaning of the
{\it radius of the indicatrix} supported by p. $x$.

\ses

In this respect, there is the deep
qualitative distinction of the Finsler geometry from
the Riemannian geometry.
Namely, in the latter geometry  we
have simply $V_{{\cal I}_x}=\const$, and whence
$R=\const$.
The new reality that the value of
$V_{{\cal I}_x}$
may vary from point to point of the background manifold $M$
arises in the Finsler geometry,
in which case the $R$ may be a function of $x$.

The $R(x)$ thus appeared proposes naturally  the  scale factor
in
the  tangent Riemannian space
$\cR_{\{x\}}$
supported by  the point $x$.

\ses

This motivation suggests the idea to replace
the above angle
$
{ \al}_{\{x\}}(y_1,y_2)
$
by
 the {\it  normalized angle}
\be
 {\al}^{\{H(x)\}}_{\{x\}}(y_1,y_2)
~:=
H(x)
{ \al}_{\{x\}}(y_1,y_2), \qquad y_1,y_2\in T_xM,
\ee
where we have introduced the scalar
$
H(x)=1/R(x),
$
to use
the  preservation equation
 \be
d_i {\al}^{\{H(x)\}}_{\{x\}}(y_1,y_2)=0
\ee
instead of
$d_i{ \al}_{\{x\}}(y_1,y_2)=0$
formulated in  (1.15).
The  preservation law
(1.27)     can be written
in the {\it recurrent} form
\be
d_i \al+ \fr1HH_i\al=0.
\ee
The  $d_i$ is the operator (1.11)
and
$
H_i=\partial H/\partial x^i.
$

\ses

Since the angle $
{ \al}_{\{x\}}(y_1,y_2)
$
is measured by the indicatrix arc length,
it seems quite natural to normalize the angle by means of the characteristic scale factor,
according to (1.26).

{

To elucidate patterns,
it proves being of great help to apply the  coincidence-limit  method
(see Section 3.2 in [12]).
Namely,
with the function
$
E=(1/2)\al^2
$
the recurrent preservation
$
d_i \al+ (1/H)H_i\al=0
$
proposed by (1.28) entails
the following $E$-equation
$$
\D{E}{x^i}  +N^k{}_{1i}   \D{E}{y^k_1} +   N^k{}_{2i}   \D{E}{y^k_2}
=
-
\fr 2HH_i
 E,
$$
\ses\\
where
$
N^k{}_{1i}  =N^k{}_{i}(x,y_1)
$
and
$
N^k{}_{2i}  =N^k{}_{i}(x,y_2).
$
Evaluating various partial derivatives of this $E$-equation
with respect to
$y_1$ and $y_2$
and finding   the coincidence limits
when
${y_2}\to{y_1}$,
we can obtain a valuable information of
the tensors of
the Finsler space.
Performing the required evaluations
on the level of
the second-order partial  derivatives
$\partial^2/\partial y_1^m\partial y_2^n$,
and, then,
 applying   the operation
${y_2}\to{y_1}$
to the resultant expressions,
it is possible to arrive at the general conclusion
that
{\it
in any Finsler space
the
vanishing assumption
$
d_i \al+ (1/H)H_i\al=0
$
 entails the equality}
\be
\cD_ih_{mn}
=
\fr 2{F}
   h_{mn}
   d_iF
-\fr2HH_i h_{mn}.
\ee

The formula (1.29)  has been  derived in Appendix E  in all detail
by  performing required long substitutions
(see  (E.37)  in Appendix E).

{

By differentiating the equality (1.29) with respect to $y^j$,
it is possible to
obtain the coefficients
$N^k{}_{imn}$.
In this way,
when the vanishing $d_iF=0$
is also keeping valid,
simple direct evaluations yield  the representation
\be
N^k{}_{imn}
=
\fr2HH_i
\fr1F
l^kh_{mn}
-
\cD_i
C^k{}_{mn},
\ee
which extends the previous (1.23).
The symmetry
(1.22) is now replaced by
$$
N^t{}_{imj}g_{tn}
-
\fr2HH_i
\fr1F
h_{mj}l_n
=
N^t{}_{imn}g_{tj}
-
\fr2HH_i
\fr1F
h_{mn}l_j.
$$
Instead of the vanishing (1.18)
 we obtain
\be
FN^k{}_{inm}l_k=
\fr2HH_i h_{mn}.
\ee
The vanishing
$\cD_i
S_n{}^k{}_{jm}
=0$
indicated in (1.24) is now extended,
namely
the above representation (1.30) straightforwardly
 entails  the equality
 $$
\cD_i
S_n{}^k{}_{jm}
=
-
\fr2HH_i
 \lf(h^k_jh_{mn}-h^k_mh_{jn}\rg).
$$

\ses

{

\ses

From (1.29) we can conclude that  when $d_iF=0$
we have
\be
{\cD}_i g_{mn}
=
-
\fr2HH_i
 h_{mn}
\ee
at an arbitrary smooth function $H=H(x)$.

\ses

The equality
(1.32)
suggests us to introduce
the {\it total connection coefficients}
\be
T^k{}_{im}
=
-N^k{}_{im}
-
\fr1HH_ih^k_m,
\ee
so that the {\it deflection tensor}
\be
\De^k{}_{im}\eqdef  -  N^k{}_{im} - T^k{}_{im}
\ee
\ses\\
is non-vanishing
as far as $H_i\ne0$, namely
\be
\De^k{}_{im}
=
\fr1HH_ih^k_m.
\ee
It follows that
\be
T^k{}_{im}y^m=
-N^k{}_{im}y^m
\equiv
-
N^k{}_i, \qquad
l_kT^k{}_{im}=-l_kN^k{}_{im}.
\ee

There arises the {\it total covariant derivative}
$\cT_i$,
showing
the properties
 \be
\cT_i F=0, \qquad \cT_i l_m=0, \qquad \cT_i l^m=0,
\ee
\ses\\
and
the {\it metricity}
\be
\cT_i g_{nm}
=
0,
\ee
\ses\\
where
\be
\cT_i F\eqdef d_i F, \qquad
\cT_i l_{m}\eqdef d_il_m-T^h{}_{im}l_h, \qquad
\cT_i l^{m}\eqdef d_il^m+T^m{}_{ih}l^h,
\ee
and
\ses\\
\be
\cT_i g_{nm}\eqdef d_ig_{nm}-T^h{}_{im}g_{hn}-T^h{}_{in}g_{hm}.
\ee

\ses

{

In all the previous formulas started with (1.26),
the $H(x)$ was an arbitrary smooth scalar not related anyhow to the indicatrix curvature,
the constancy of the indicatrix curvature was not implied,
and the Finsler space was arbitrary.

\ses

If the indicatrix of a Finsler space is a space of constant curvature
at any point $x\in M$,
we say that
the  Finsler space is
the ${\cF}^N$-{\it space},
where $N\ge3$ is the dimension of the space.

\ses

The interest to the Finsler space
$\cF^N$
is motivated by the following important observations.
Given an arbitrary
 Finsler space
of  any dimension
$N\ge3$.
The  tangent Riemannian  space
$\cR_{\{x\}}\subset T_xM$
is conformally flat if and only if the  indicatrix
${\cal I}_x\subset T_xM$
is a space of constant curvature,
 assuming naturally
  that the involved conformal multiplier
 is  homogeneous with respect to the argument
$y$.
The dependence of the conformal multiplier  on the variable $y$
is
 presented by the  power of the Finsler metric function.
The remarkable equality
 ${\mathcal  C}_{\text{Ind.}} \equiv H^2$
ensues.
These observations form the content of Proposition II.2.1
(formulated and proved in Section II.2 of Chapter II),
which
extends  Proposition 2.1 of the preceding work [10,11]
in the following essential aspect.

\ses

In [10,11],
the assumption was made that
 the respective conformal multiplier
is  of the  power dependence on the Finsler metric function,
in accordance with the  representations indicated in the formula  (II.2.3) of Section II.2.
In proving  Proposition II.2.1 in Section II.2,
we  outline    the reasoning line which explains that the representations
are actually the direct consequences of
the property
that the
 indicatrices are  spaces of constant
curvature.

{

\ses

We say that
the Finsler space
${\mathcal F}^N$
is {\it indicatrix-homogeneous}
if
$
{\mathcal  C}_{\text{Ind.}}=\const$.
In this case,
the  deflectionless
connection has been derived from
 the separable angle-preservation equation
in the preceding work [10,11].

Alternatively,
the Finsler space
${\mathcal F}^N$
is said to be
{\it indicatrix-inhomogeneous}
if ${\mathcal  C}_{\text{Ind.}}={\mathcal  C}_{\text{Ind.}}(x)$.
On  this level,
    because of the equality
     ${\mathcal  C}_{\text{Ind.}} \equiv H^2,$
we have
 $H=H(x)$ and  $H_i \ne 0$.

\ses

{

On the  indicatrix-inhomogeneous level
of study of the Finsler space
${\mathcal F}^N$
with
$d_iF=0$  the separable preservation law for the angle is
  {\it impossible to introduce}.
Indeed,
the  law entails the
  metricity
$\cD_i g_{mn}=0$
of the deflectionless type (see (1.17) and the definition (1.6)),
together with the representation (1.23) for the coefficients
$N^k{}_{imn}$
and  the  vanishing
$
\cD_i
S_n{}^k{}_{jm}
=0,
$
where
$
S_n{}^k{}_{ij}
=
\lf(
C^h{}_{nj}    C^k{}_{hi}
-
C^h{}_{ni}    C^k{}_{hj}
\rg)
F^2
$
(see  (1.24)).
 It is known that the indicatrix is a space of constant curvature if and only if the last tensor
 fulfills the equality
$ S_n{}^k{}_{ij} =C(h_{nj}h^k_i-h_{ni}h^k_j)$ with the  factor $C$ which is independent of $y$,
in which case
$
{\mathcal  C}_{\text{Ind.}}=1-C
$
(see Section 5.8 in [1]).
In  the Finsler space
${\mathcal F}^N$,
we have
$
{\mathcal C}_{\text{Ind.}}=H^2.
$
The two vanishings
$\cD_i g_{mn}=0$ and
$\cD_i F=0$ entail
$\cD_i h_{mn}=0$.
Whence from $\cD_i
S_n{}^k{}_{jm}=0$ it follows that
$H_i=0$.

\ses

If, however, we start with recurrent preservation law
supplemented by the vanishing condition
$\cD_i F=0$,
then from
(1.29) we have
${\cD}_i h_{mn}
=
-
(2/H)
H_i
 h_{mn}.
 $
Applying the covariant derivative
$\cD_i$
to the tensor
$ S_n{}^k{}_{ij} =C(h_{nj}h^k_i-h_{ni}h^k_j)$
and taking into account that
$C=1-H^2$,
after short evaluations we now  arrive at
 the equality
$$
\cD_i
S_n{}^k{}_{jm}
=
-
\fr2H
H_i
 \lf(h^k_jh_{mn}-h^k_mh_{jn}\rg)
$$
which is equivalent to the implication
written below   (1.31).
 Thus,
 the following proposition is valid.

\ses

\ses

{\bf Proposition I.1.2.}
{\it  The recurrent-type preservation}
 (1.28)
{\it of the angle,
 that is,
}
$
d_i \al+ (1/H)H_i\al=0,
$
{\it
is reconciled
with
the indicatrix-inhomogeneous
Finsler space
}
${\mathcal F}^N$
{\it at any scalar}
 $H=H(x)$
{\it obtainable from
the identification}
$
{\mathcal C}_{\text{Ind.}}=H^2.
$

\ses

\ses

The observations motivate us  to go
to the preservation law (1.27)
 which is not separable
 from the standpoint
 of the indicatrix-arc angle
  $ \al_{\{x\}}(y_1,y_2)$,
  whenever
 $H\ne const.$

{

In so doing,
the coefficients
$N^{m}{}_n$
of the Finsler space
${\mathcal F}^N$
are obtained to read
   (I.2.16) in Section I.2.
   They
don't involve explicitly the gradients $H_n$.
If, however, we expand the partial derivatives $\partial/\partial x^n$
which enter the right-hand part of (I.2.16),
the coefficients will break down  into two parts:
\be
 N^m{}_n  = N^{{{\rm I}}m}{}_n+\breve N^m{}_n,   \qquad \breve N^m{}_n=\breve N^mH_n.
 \ee
Here,
the first part
$N^{{{\rm I}}m}{}_n$
are the coefficients of the indicatrix-homogeneous case
(given by the formula (2.30) in [10],
and by the formula (2.36) in [11])
in which the constant $H$
has been
merely
replaced by arbitrary $H(x)$,
and the vector field
$\breve N^m$
does not involve any gradient of $H(x)$.
We may say that the coefficients
$ N^m{}_n $
are of the {\it  linear dependence} on the gradient $H_n$.

\ses

The entailed
  coefficients
$N^k{}_{mn}$
are given by the representation (II.3.32) of Chapter II
which is
 applicable to any
 indicatrix-inhomogeneous Finsler space
${\mathcal F}^N$.
It is also possible to evaluate explicitly
 the derivative coefficients
$N^k{}_{mni}=\partial{N^k{}_{mn}}/\partial{y^i}$.
The required evaluations
lead straightforwardly
to
the validity of the representation
(1.30)
in the
 ${\mathcal  F}^N$-space
with
  an arbitrary smooth function $H(x)$,
provided the vanishing $d_nF=0$
is assumed
(see Proposition II.3.5 in Chapter II).

Having evaluated the coefficients
$N^k{}_{mn}$,
we  obtain from (1.33) the
 total connection coefficients
$
T^k{}_{im}
$
thereby solving the problem of finding the connection in
 the
 ${\mathcal  F}^N$-space at
the   indicatrix-inhomogeneous level.
The  coefficients
$
T^k{}_{im}
$
 involve
 the  deflection tensor
$
\De^k{}_{im}
$
indicated in (1.34) and (1.35).
There arises the   covariant derivative
$\cT$,
which properties are  listed in (1.36)-(1.40).

{

Section I.2
gives a brief summary
of
 Chapter II.

 \ses

The formula (I.2.16)
 indicates the representation of
the coefficients
$N^{m}{}_n$
which is
valid for an arbitrary  Finsler space of the type
 ${\cF}^N$.
The representation involves the vector field
$U^i$ which realizes the
key transformation
$  y={\bf C}(x,\bar y)$
indicated in (I.2.1).
Given a particular
 Finsler space of the type
 ${\cF}^N$,
 the formula (I.2.16)
yields
the coefficients
$N^{m}{}_n$
in a completely explicit way
when the respective
field
$U^i$
is known.

The Finsleroid case to which  Section I.3 is devoted
provides us with such an example,
for
the required
field
$U^i$
is explicitly given,
namely  by means of the representation
(I.3.20) (which was earlier found in Section 6 of [7]).
Therefore, we can straightforwardly
 apply
 the developed theory
 of the
 ${\mathcal  F}^N$-space
 to
the metric function of the Finsleroid type.
The expansion (1.41) for the respective Finsleroid  coefficients
$ N^m{}_n$
has been evaluated.
The explicit representation of the entailed derivative coefficients
$ N^k{}_{im}$
is  indicated.
The respective validity of the representations  (1.29) and (1.30)
of the tensors $\cD_ih_{mn}$
and
$
N^k{}_{imn}
$
on the  indicatrix-inhomogeneous level
of study of the Finsleroid space
 has been   verified
by direct evaluations presented in detail.

\ses

Several Appendices
are added
in which numerous fragments
of the underlined evaluations
have been displayed.

{

\ses

\ses

 \setcounter{equation}{0}

{\bf I.2.  Indicatrix of constant curvature }

\ses

\ses

\ses

Let $M$ be the base manifold, such that $\cF^N=(M,F)$,
where $F=F(x,y)$ is the Finsler metric function
and $N\ge3$ is the dimension of the space.
If the indicatrix of a Finsler space is a space of constant curvature,
we say that
the  Finsler space is
the ${\cF}^N$-{\it space}.
Denote by
${\mathcal  C}_{\text{Ind.}}$
the value of curvature of the indicatrix supported by the point $x\in M$.
If
${\mathcal  C}_{\text{Ind.}}$ is a constant over the manifold $M$, we say that the
space ${\cF}^N$ is of the {\it indicatrix-homogeneous} case.

In general, the value
${\mathcal  C}_{\text{Ind.}}$ may vary from point to point of $M$,
in which case we say that
the  space ${\cF}^N$ is
of the {\it indicatrix-inhomogeneous} type.
The possibility is characterized by a function
${\mathcal  C}_{\text{Ind.}}={\mathcal  C}_{\text{Ind.}}(x)$
such that the derivative
$\partial{\mathcal  C}_{\text{Ind.}}/\partial x^i$
does not vanish identically.

\ses

\ses

{

In such spaces,
 the transformation
\be
  y={\bf C}(x,\bar y), \quad  y,\bar y\in T_xM,
\ee
can be proposed
which
 maps
 the tangent vectors $y\in T_xM$ into the tangent vectors
of the same tangent space $T_xM$,
subject to the following conditions.
The  transformation is non-linear with respect to $\bar y$. Non-singularity and
sufficient smoothness are  implied.
Also, the  transformation
is  positively   homogeneous
of a degree $H(x)$
regarding dependence
on
 tangent vectors $y$.
Each tangent Riemannian space
$\cR_{\{x\}}=\{T_xM, g_{\{x\}}(y)\}$
is conformally transformed to Euclidean space,
to be denoted by
$\cE_{\{x\}}$.
The distribution of the last spaces
$\cE_{\{x\}}$
over the base manifold $M$
composes the
{\it associated Riemannian space},
which we denote by
$\cR^N=(M,S)$,
where
$ S=\sqrt{a_{mn}(x) y^m y^n} $
is the Riemannian metric constructed from
the  metric tensor
$a_{mn}(x)$
of the space
$\cE_{\{x\}}$.

Under these conditions,
the scalar $H(x)$ can be taken
from the identification
\be
{\mathcal  C}_{\text{Ind.}} \equiv H^2.
 \ee
The equality
\be
S(x,\bar y)=\lf(F(x,y)\rg)^{H(x)}
\ee
arises
(see (II.2.10)),
which validates the indicatrix correspondence to the Euclidean sphere;
$ S(x,\bar y)=\sqrt{a_{mn}(x)\bar y^m\bar y^n}. $
The relevant  conformal multiplier $p^2$
is
 constructed from the Finsler metric function $F$,
 according to
\be
p=
\fr1H
F^{1-H}.
\ee
We take $1>H>0$ for definiteness,
the extension of the approach to other  values
of $H$
being a straightforward task.

{

\ses

If  $f(x,y)$ is the involved conformal multiplier
in the
tangent Riemannian space
$\cR_{\{x\}}$,
 then
 the equality
$$
g_{\{x\}}(y)=f(x,y) u_{\{x\}}(y)
$$
should introduce  the tensor $u_{\{x\}}(y)$  which associated
 Riemannian curvature tensor
vanishes identically.
The function $f(x,y)$ is assumed naturally to be homogeneous with respect to the argument
$y$.
Denoting
the
homogeneity degree
of $f(x,y)$
by means of $2a(x)$,
we just conclude that
the difference $1-a$ is exactly the homogeneity degree of the transformation (2.1) considered,
that is,
$$
H=1-a.
$$

The following assertions are  valid.
A Finsler space  is  the
${\mathcal F}^N$-space if and only if  the indicatrix of the Finsler space
 is a space of constant curvature.
The dependence of the multiplier $f$ on the variable $y$
is presented by the  power of the Finsler metric function
$F$
(see  Proposition II.2.1 in Section 2 of Chapter II).

\ses

{

The respective  two-vector angle
 $ \al_{\{x\}}(y_1,y_2)$
proves to be obtainable from
 the angle
$ \al^{\text{Riem}}_{\{x\}}(y_1,y_2) $
operative  in the Riemannian space,
namely the simple equality
\be
  \al_{\{x\}}(y_1,y_2)=\fr1{H(x)}  \al^{\text{Riem}}_{\{x\}}(\bar y_1,\bar y_2)
\ee
(see (II.2.51)-(II.2.52))
is valid.

\ses

{

We locally represent the transformation (2.1) by means of the functions
\be
y^i=y^i(x,t), \qquad t^n \equiv \bar y^n.
\ee
The homogeneity entails
$
y^i(x,k t)=k^{1/H} y^i(x,t)
$
with
$ k>0$ and $ \forall t,
$
together with
$
y^i_nt^n=(1/H)y^i,
$
where
$ y^i_n=\partial y^i/\partial t^n$.

The definition
\be
U^i\eqdef (1/S)\bar y^i \equiv (1/F^H)\bar y^i
\ee
introduces  the normalized vector,
which is obviously unit:
$
U_iU^i=1$ and $    U_i=a_{ij}U^j.
$
The zero-degree homogeneity
$
U^i(x,k y)= U^i(x,y)$
with
$ k>0$ and $ \forall t$
holds, entailing the identity
$
U^i_ny^n=0
$
with
\be
 U^i_n~:=\D{ U^i}{ y^n}
 =
\fr1{F^H} t^i_n-\fr1FH U^il_n,
\ee
where
$t^i_n=\partial t^i/\partial y^n$.
It follows that
 \be
F^H U^h_sy^k_{h}
 =h^k_s, \qquad
F^H U^i_ky^k_t
 =\de^i_t -U^iU_t, \qquad
 U_i U^i_n=0.
\ee

The vanishing
\be
U_i\lf( \D{U^i}{x^n}+L^i{}_{kn} U^k\rg)=0
\ee
 holds obviously,
where
$L^i{}_{nk}$
are the Riemannian connection coefficients (I.1.14).

\ses

{

The representation  (2.5) of the angle
takes on the simple form
\be
   \al_{\{x\}}(y_1,y_2)=\fr1{H(x)}\arccos\la, \quad  \text{with} ~~
\la =
a_{mn}(x)U_1^m  U_2^n,
  \ee
\ses\\
where $U_1^m = U^m (x,y_1)$    and
$U_2^m = U^m (x,y_2).$

\ses

When the recurrent preservation
$
d_i \al+ (1/H)H_i\al=0
$
proposed by (1.28)
is applied to the angle given in (2.11),
we obtain simply
\be
d_i\la=0,
\ee
where $d_i$
is the separable operator (1.11).
That is,
the recurrent preservation law formulated for
the Finsler
  ${\mathcal  F}^N${\it-space angle}
$\al_{\{x\}}$ given by
(2.11) is tantamount to the separable
preservation law
for the Euclidean angle
$
  \al^{\text{Riem}}_{\{x\}}=  \arccos\la$,
  whence to
the separable
preservation law  (2.12).

The form of the right-hand part in the formula
$
\la = a_{mn}(x)U_1^m  U_2^n
$
 is such that
the law (2.11) is obviously equivalent to the
vanishing
\be
{\mathcal D}_nU^i=0
\ee
for the field $U^i=U^i(x,y)$,
where
we introduced the covariant derivative
\be
{\mathcal D}_nU^i~:=d_nU^i+L^i{}_{nk} U^k.
\ee

Since
$$
d_nU^i=\D{U^i}{x^n}+N^k{}_nU^i_k,
$$
we arrive at  the conclusion that
in  the ${\mathcal F}^N$-space,
the   coefficients
$N^{m}{}_n $
can unambiguously be found
from the equation
$
d_n\lf(H(x)
{ \al}_{\{x\}}(y_1,y_2)\rg)=0
$
 to be given explicitly by
the  representation

{

\be
N^{m}{}_n
=-y^m_iF^H
 \lf(\fr HFU^i\D{F}{x^n}+
 \D{U^i}{x^n}+
 \lf(a^i{}_{nk}+S^i{}_{nk}\rg) U^k\rg)
 +
l^m
d_nF
\ee
(see (II.3.12) in Chapter II).
Here,
$a^i{}_{nk}$
are the Riemannian  Christoffel symbols;
 $S^i{}_{nk}=S^i{}_{nk}(x)$ is an arbitrary  torsion tensor,
 that is, the tensor possessing the skew-symmetry property
$ S_{ink}=-S_{kni},$
where
$ S_{ink}=a_{ij}S^j{}_{nk}$.

\ses

Whenever $d_nF=0$,
the  representation (2.15)
takes on the form
\be
N^{m}{}_n
=-l^m\D{F}{x^n}
-y^m_iF^H
 \lf(
 \D{U^i}{x^n}+\lf(a^i{}_{nk}+S^i{}_{nk}\rg) U^k\rg)
\ee
(see (II.1.19) in Chapter II).
These coefficients
$N^{m}{}_n$
present the {\it general solution}
to the couple equations
$
d_n\lf(H(x)
{ \al}_{\{x\}}(y_1,y_2)\rg)=0
$
and
$d_nF=0$,
so that no problem of
 uniqueness of connection coefficients may  be questioned.
The entrance of the torsion tensor
$S^i{}_{kn}$ is the only freedom,
in complete analogy to the connection coefficients of the Riemannian space.

\ses

The evaluations performed in Section II.3 of Chapter II
have  arrived also  at the representation
\ses\\
\be
N^m{}_{n}
=
d^{\text{Riem}}_n  y^m(x,t)+  \fr1H H_ny^m\ln F
\ee
(see (II.3.29) in Chapter II)
which is alternative to (2.16); here, $y^m=y^m(x,t)$ are the functions
(2.6).

\ses

The representations
(2.15)-(2.17) involve the gradient $H_n$
and are applicable to any
 indicatrix-inhomogeneous Finsler space
${\mathcal F}^N$.

{

The  coefficients
$N^k{}_{mn}$
can be evaluated from (2.16) to read
\ses\\
$$
N^k{}_{mn}
=
-
\fr1F
h^k_n\D{F}{x^m}
-
l^k\D{l_n}{x^m}
-
 C^k{}_{ns}
  N^s{}_m
+
\fr1{F}
\lf(
l_nh^k_s
-
(1-H)
l^kh_{ns}
\rg)
 N^s{}_m
$$

\ses

\ses

\be
-
y^k_hF^H
 \lf( \D{ U^h_n}{x^m}+ L^h{}_{ms}U^s_n\rg)
\ee
(see  Proposition II.3.4 in Chapter II).
With these coefficients,
the validity of the representation
(1.30) for the entailed coefficients $N^k{}_{imn}$
can straightforwardly be verified
(see  Proposition II.3.5 in Chapter II).

{

 The  space ${\mathcal  F}^N$ is obtainable
  from
  the Riemannian
 space
 ${\mathcal  R}^N$
  by means of the  deformation
$  y={\bf C}(x,\bar y)$
(see (II.2.1) in Chapter II)
which
can be presented by the {\it  deformation tensor}
\be
{ C}^i_m~:=p\bar y^i_m,
\ee
so that
\be
 g_{mn}= { C}^i_m  { C}^j_n a_{ij}
\ee
\ses
and the zero-degree homogeneity
\be
{ C}^i_m(x,ky)=  { C}^i_m(x,y), \qquad k>0,\forall y,
\ee
holds,
together with the identity
\be
{ C}^i_m(x,y)y^m=\lf(F(x,y)\rg)^{1-H}\bar y^i
\ee
\ses
(see (II.2.24)-(II.2.27)).
In Section II.4 we show that
the  $\bfC$-deformation  is
$\cT$-{covariant constant:
\be
{\cT}\cdot{\bfC}=0,
\ee
\ses\\
where $\cT$ designates the covariant derivative introduced by the help of the formulas
(I.1.33)-(I.1.40)
(see Proposition II.4.1).

\ses

{

Also, the covariant derivative $\cT$
is  the
manifestation of the
{\it transitivity}
of the connection  under   the $\cC$-transformation,
in short,
\ses\\
\be
\cT=\cC\cdot\nabla,
\ee
where
$\nabla$
is
 the
  covariant derivative
 applicable
in
the background  Riemannian space
 ${\mathcal  R}^N$
 (see Proposition II.4.2).
 In other words,
in    the  Finsler space
${\mathcal  F}^N$
 the metrical non-linear angle-preserving connection
 is  the $\bfC$-export of
 the
  metrical  linear   connection
(II.1.2)
applicable in
the  space
 ${\mathcal  R}^N$.

{

\ses

\ses

In  Section II.5 we perform the attentive
comparison between the commutators
of the involved Finsler  covariant derivative $\cT$ and
 the commutators of  the underlined
Riemannian
covariant derivative $\nabla$,
 {\it not} assuming $H=\const$, such that
 $H(x)$ can be an arbitrary smooth function of $x$.
In this way,  we derive
the associated curvature tensor
 $\Rho_k{}^n{}_{ij}$.
Important properties of the tensor are elucidated.

\ses

{

\ses

\ses

 \setcounter{equation}{0}

{\bf I.3. Reduction to the Finsleroid  space }

\ses

\ses

\ses

\ses

In the Finsleroid  case,
we make the notation change $H(x)\to h(x)$.

\ses

The scalar $g(x)$ obtained through
\be
h(x)=\sqrt{1- \fr{g^2(x)}4}, \qquad \text{with} \quad -2<g(x)<2,
\ee
plays the role of the characteristic parameter.

It follows that
\ses\\
\be
g_i=-\fr{4h}gh_i,
\ee
\ses\\
where $g_i=\partial g/\partial x^i$
and
$h_i=\partial h/\partial x^i$.

\ses

We  assume that in addition to a
Riemannian  metric
 $\sqrt{a_{ij}(x)y^iy^j}$
the manifold $M$ admits a non-vanishing 1-form
$ b=b_i(x)y^i$ of the unit length:
\be
a_{ij}(x)b^i(x)b^j(x)=1,
\ee
where
$b^i(x)=a^{ij}(x)b_j(x).$
The tensor $a^{ij}(x)$ is reciprocal to $a_{ij}(x)$, so that
$a_{ij}a^{jn}=\de^n_i$, where $\de^n_i$ stands for the Kronecker symbol.
We need also the quadratic form
\be
B=b^2+gbq+q^2\equiv \lf(b+\fr12gq\rg)^2+h^2q^2,
\ee
where
\be
q=\sqrt{r_{mn}y^my^n} \quad \text{with} \quad     r_{mn}=a_{mn}-b_mb_n,
\ee
so that
\be
 a_{ij}(x)y^iy^j=b^2+q^2.
 \ee

{

\ses

We shall also use the scalar
\be
\chi=\fr1h \Bigl(
-\arctan   \fr G2   +\arctan\fr{L}{hb}\Bigr),    ~  {\rm if}  ~ b\ge 0;
\quad
\chi=\fr1h \Bigl(
 \pi-\arctan
\fr G2
+\arctan\fr{L}{hb}\Bigr),
~   {\rm if}
~ b\le 0,
\ee
 with the function
$    L =q+(g/2) b $
fulfilling   the identity
\be
 L^2+h^2b^2=B.
 \ee

The definition range
$$
0\le\chi\le\fr1h\pi
$$
is of value to describe all the tangent space.
The normalization in (3.7)
is such that
\be
\chi\bigl|_{y=b}\bigr. =0.
\ee
The quantity (3.7) can conveniently be written as
\be
\chi  =  \fr1h  f
\ee
with
the function
\be
f=\arccos \fr{ A(x,y)}   {\sqrt{B(x,y)}}
\ee
ranging as follows:
\be
0\le f\le \pi.
\ee
The Finsleroid-axis vector $b^i$ relates to the value $f=0$, and
the opposed vector $-b^i$ relates to the value $f=\pi$:
\be
f=0 ~~  \sim ~~ y=b;  \qquad f=\pi ~~ \sim ~~ y=-b.
\ee

{

With these ingredients,
we construct the Finsler metric function
\be
K=\sqrt B\, J, \qquad \text{with} ~~ J=\e^{-\frac12 g \chi}.
\ee
The normalization is such that
\be
K(x,b(x))=1
\ee
(notice that
 $q=0$ at $y^i=b^i$).
The positive  (not absolute) homogeneity  holds:
 $K(x, \ga y)=\ga K(x,y)$ for any $\ga >0$   and all admissible $(x,y)$.

 Under these conditions, we call $K(x,y)$
the {\it  ${\mathbf\cF\cF^{PD}_{g}}$-Finsleroid   metric function},
obtaining  the ${\mathbf\cF\cF^{PD}_{g}}$-{\it Finsler space}
\be
{\mathbf\cF\cF^{PD}_{g}} :=\{M;\,a_{ij}(x);\,b_i(x);\,g(x);\,K(x,y)\}.
\ee

\ses

\ses

 {\large  Definition}.  Within  any tangent space $T_xM$, the  metric function $K(x,y)$
  produces the {\it    ${\mathbf\cF\cF^{PD}_{g} } $-Finsleroid}
 \be
 \cF\cF^{PD}_{g;\,\{x\}}:=\{y\in   \cF\cF^{PD}_{g; \, \{x\}}: y\in T_xM , K(x,y)\le 1\}.
  \ee

\ses

 \ses

 {\large  Definition}. The {\it    ${\mathbf\cF\cF^{PD}_{g} } $-Indicatrix}
 $ {\cal I}\cF^{PD}_{g; \, \{x\}} \subset T_xM$ is the boundary of the
    ${\mathbf\cF\cF^{PD}_{g} } $-Finsleroid, that is,
 \be
{\cal I}\cF^{PD}_{g\, \{x\}} :=\{y\in {\cal I}\cF^{PD}_{g\, \{x\}} : y\in T_xM, K(x,y)=1\}.
  \ee

\ses

 \ses

 {\large  Definition}. The scalar $g(x)$ is called
the {\it Finsleroid charge}.
The 1-form $b=b_i(x)y^i$ is called the  {\it Finsleroid--axis}  1-{\it form}.

\ses

\ses

The entailed components $y_i :=(1/2)\partial {K^2}/ \partial{y^i}$)
of the  covariant tangent vector $\hat y=\{y_i\}$
can be found in the simple form
\be
y_i=(u_i+gqb_i) J^2,
\ee
where $u_i=a_{ij}y^j$.

{

Let us elucidate the structure of the coefficients
$N^k{}_{m} $
in the Finsleroid  case proper.
From  (6.26) of [7]
it follows that the quantity $U^i=(1/K^h)\bar y^i$
can  explicitly  be given by
\be
U^i=\lf[hv^i  +   \lf( b+\fr12gq\rg)b^i    \rg] \fr 1{\sqrt B},
\ee
where
$
v^i=y^i-bb^i.
$
So we have
\ses\\
$$
\D {U^i}g=
-\fr g{4h}v^i \fr 1{\sqrt B}
 +  \fr12qb^i     \fr 1{\sqrt B}
-
\fr1{2B}U^iqb,
$$
\ses\\
or
$$
\D {U^i}g=
-\fr g{4h^2}
U^i
+\fr g{4h^2}
( b+\fr12gq)b^i     \fr 1{\sqrt B}
 +  \fr12qb^i     \fr 1{\sqrt B}
-
\fr1{2B}U^iqb.
$$

{

Since
$$
K^hy^m_iU^i=\fr1hy^m
$$
(a consequence of the homogeneity involved)
and
\ses\\
$$
K^hy^m_i b^i=
\Biggl[
b^m
+
\fr1{B}\Biggl(
\fr{1}{h}
\lf(b+\fr12gq\rg)
-b -gq
\Biggr)
y^m
\Biggr]
\sqrt B
$$
(see  (D.12) in  [7]),
we can straightforwardly evaluate the contraction
$$
K^hy^m_i\D {U^i}g=
-\fr g{4h^2}
\fr1hy^m
-
\fr1{2B}\fr1hqb
y^m
+
\fr g{4h^2}
( b+\fr12gq)
b^m
+
\fr g{4h^2}
\fr1{B}
\fr{1}{h}
\lf(b+\fr12gq\rg)^2
y^m
$$

\ses

\ses

$$
-
\fr g{4h^2}
\fr1B( b+\fr12gq)
(b +gq)
y^m
 +
   \fr12q
\Biggl[
b^m
+
\fr1{B}\Biggl(
\fr{1}{h}
\lf(b+\fr12gq\rg)
-b -gq
\Biggr)
y^m
\Biggr].
$$

{

\nin
Using the equality
$$
\lf(b+\fr12gq\rg)^2=B-h^2q^2
$$
(see (3.4))
leads to
the representation
\ses\\
$$
K^hy^m_i\D {U^i}g
=
\fr g{4h^2}
\lf( b+\fr12gq\rg)
b^m
-
\fr g{4}
\fr1{B}
\fr{1}{h}
q^2y^m
-
\fr g{4h^2}
\fr1B
\lf( b+\fr12gq\rg)
(b +gq)
y^m
$$

\ses

\ses

$$
 +  \fr12q
\Biggl[
b^m
+
\fr1{B}\Biggl(
\fr{1}{h}
\fr12gq
-b -gq
\Biggr)
y^m
\Biggr],
$$
\ses\\
which can be simplified as follows:
$$
K^hy^m_i\D {U^i}g
=
\fr g{4h^2}
\lf( b+\fr12gq\rg)
b^m
-
\fr g{4h^2}
\fr1B\lf( b+\fr12gq\rg)
(b +gq)
y^m
$$

\ses

\ses

$$
 +  \fr1{2h^2}q
\Biggl[
b^m
\lf(1-\fr{g^2}4\rg)-
\fr1{B}
(b+gq)
y^m
\lf(1-\fr{g^2}4\rg)
\Biggr]
$$

\ses

\ses

\ses

$$
=
\fr 1{2h^2}
\lf( q+\fr12gb\rg)
b^m
-
\fr 1{2h^2}
\fr1B
\lf( q+\fr12gb\rg)
(b +gq)
y^m,
$$
\ses
so that
\ses\\
$$
K^hy^m_i\D {U^i}g
=
\fr 1{2h^2}
\fr1B
\lf( q+\fr12gb\rg)
\lf[
Bb^m
-
(b +gq)
y^m
\rg].
$$
\ses
By comparing this result with the representation
$$
A^m=\fr N2g\fr 1{qK}
\Bigl[q^2b^m-(b+gq)v^m\Bigr]  \equiv KC^{mn}{}_n
$$
(see  (A.27)  in  [7]),
we come to
\be
K^hy^m_i\D {U^i}g
=
\fr 1{h^2}
\fr qB
\lf( q+\fr12gb\rg)
\fr{K}{Ng}
A^m.
\ee

{

Therefore, in the Finsleroid  case
the coefficients
$N^k{}_{i} $
proposed by (I.2.16)
are  the sum
\be
 N^k{}_i  = N^{{{\rm I}}k}{}_i+\breve N^k{}_{i}, \qquad \breve N^k{}_{i}=\breve N^k{}g_i,
 \ee
 \ses\\
 where
 \be
\breve N^k{}
=
-
\fr1{h^2}
\fr{q}B\lf(q+\fr12 gb\rg)
\fr K{Ng}
A^k
-
\fr12
\bar My^k
\ee
\ses
with $\bar M$ coming from
\ses
\be
 \D {K^2}g=\bar M  K^2.
\ee
The torsion tensor $S^k{}_{ij}=S^k{}_{ij}(x)$ has been  neglected.
The
 $
N^{{{\rm I}}k}{}_i
$ are the coefficients  (6.48)  of  [7] ( they can also be found in [10,11]),
namely,
\be
N^{{{\rm I}}k}{}_i=\Biggl[\lf(b\!-\!\fr1h\lf(b+\fr12gq\rg)\rg)
\eta^{kj}
\!+\!
\lf(
\fr1{q^2} v^k \lf(b \!  -\!  \fr1h(b+gq)\rg)
\!+
\!\lf(\fr1h\!-\!1\rg)b^k
\rg)
y^j
\Biggr]
 \nabla_ib_j
-     a^k{}_{ij}y^j.
\ee
\ses\\
They don't involve the gradient $g_i$.
The tensor
\be
\eta^{kn}=a^{kn}-b^kb^n-\fr1{q^2}v^kv^n
\ee
enters the representation.
This tensor obeys the nullification
\be
y_k\eta^{kn}=b_k\eta^{kn}=0.
\ee
\ses\\
The designation
$ \nabla_i$ stands for the Riemannian covariant derivative
constructed with the help of the Riemannian Christoffel symbols
$a^k{}_{ij}=a^k{}_{ij}(x)$.

\ses

The
 $
N^{{{\rm I}}k}{}_i
$
 are the coefficients
$N^k{}_{i} $
obtained when the condition $h=\const$ which specifies the
indicatrix-homogeneous case
is postulated.

\ses

{

For the coefficients
$$
\breve N^k{}_{im}=\D{\breve N^k{}_{i}}{y^m}
$$
\ses\\
the representation
$$
\breve N^k{}_{im}
=
\fr1{h^2}g_i
\fr{q^2}{2B}\lf(1+\fr12 g\fr bq-2h^2 \rg)
\fr2{Ng}
A_m
 l^k
+
\fr1{h^2}g_i
\fr{q^2}{2B}
\lf(1+\fr12 g\fr bq\rg)
\lf(\fr bq+g\rg)
h^k_m
$$

\ses

\ses

\be
+
\fr1{h^2}g_i
\fr{q^2}{2B}
\lf(
\fr bq+\fr12g
\rg)
\fr2{Ng}
\fr2{Ng}
A_m  A^k
-
\fr12
g_i
\bar Mh^k_m
+
\fr1K
l_m\breve N^k{}_{i}
\ee
\ses\\
is obtained (see Appendix A).

Using (3.28) we find  straightforwardly  that
\be
y_k\Dd{\breve N^k{}_i}{y^m}{y^n}=
\fr2h
h_i
h_{mn}.
\ee

For the coefficients
$$
\breve N^k{}_{imn}=\D{\breve N^k{}_{im}}{y^n}
$$
\ses\\
the  representation
\ses\\
\be
\breve N^k{}_{imn}
=
-
\fr g{2h^2}g_i
\fr1K
h_{mn}
 l^k
-
\fr1{gh^2}g_i
\fr1K
A^k{}_{mn}
\ee
\ses\\
can explicitly  be derived
(see Appendix A);
$A^k{}_{mn}=KC^k{}_{mn}$.

\ses

\ses

The full coefficients read
\be
N^k{}_{imn}
=
\fr2hh_i
\fr1K
l^kh_{mn}
-
\fr1K
\cD_i
A^k{}_{mn}
\ee
(see Appendix A).
Thus in the Finsleroid  case proper
we have straightforwardly verified the validity of the representation (1.30).

{

\ses

\ses

\ses

{\large \bf Chapter II. Phenomenon of indicatrix of constant  curvature

\ses

 with
$ {\mathcal  C}_{\text{Ind.}}={\mathcal  C}_{\text{Ind.}}(x)$}

\ses

\ses

\setcounter{equation}{0}

\bigskip

 {\bf II.1. Motivation}

\ses

\ses

In any dimension $N\ge3$
the Finsler metric function $F$
geometrizes the tangent bundle $TM$ over the base manifold $M$
such that
at each point $x\in M$
the tangent space $T_xM$ is  endowed with the curvature tensor
constructed from the respective Finslerian metric tensor  $g_{\{x\}}(y)$
by means of the conventional rule of the Riemannian geometry
considering $y$ to be the variable argument.
There arises the Riemannian space
$\cR_{\{x\}}=\{T_xM, g_{\{x\}}(y)\}$
 supported by the point $x\in M$
 such that $T_xM$ plays the role of the base manifold for the space.
We call
$\cR_{\{x\}}$
the {\it tangent Riemannian space}.

Given an $N$-dimensional  Riemannian space
${\mathcal R}^N=(M,S)$, where $S$ denotes the Riemannian metric function,
one may endeavor to  obtain a Finsler space
${\mathcal F}^N=(M, F)$ by applying an appropriate
transformation
 ${\bf C}$ to tangent spaces.
The  base manifold $M$ is keeping the same for both the spaces,
${\mathcal R}^N$ and ${\mathcal F}^N$.

We assume that the transformation ${\bf C}$ is {\it restrictive}, in the sense that no point $x\in M$ is shifted under the transformation, so that in each tangent space
$T_xM$ the deformation maps tangent vectors $y\in T_xM$ into the tangent vectors
of the same
$T_xM$:
\be
  y={\bf C}(x,\bar y), \quad  y,\bar y\in T_xM.
\ee
In general, this transformation is non-linear with respect to $\bar y$. Non-singularity and
sufficient smoothness are always  implied.

{

We may evidence in the
Riemannian
space ${\mathcal  R}^N$
the {\it metrical  linear Riemannian  connection}  ${\mathcal  RL}$, which
in terms of local coordinates $\{x^i\}$ introduced in $M$
is given by
 \be
{\mathcal  RL}=\{L^m{}_j,L^m{}_{ij}\}: \qquad    L^m{}_j=- L^m{}_{ji}y^i,
   \quad
   L^m{}_{ij}=a^m{}_{ij}+S^m{}_{ij},
\ee
where  $a^m{}_{ij}=a^m{}_{ij}(x)$ stands for the Christoffel symbols
constructed from the Riemannian metric tensor $a_{mn}(x)$
of the space ${\mathcal  R}^N$
and
 $S^m{}_{ij}=S^m{}_{ij}(x)$ is an arbitrary {\it torsion tensor}:
$S_{mij}=-S_{jim}$
with
$S_{mij}=a_{mh}S^h{}_{ij}$.
The respective covariant derivative $\nabla$ can be introduced
in the natural way.
Namely,
 considering the (1,1)-type
tensor
$ W^n_m(x,y)$
on the tangent bundle
associated to  the space ${\mathcal  R}^N$,
we can  take the definition
\be
\nabla_iW^n{}_m
=
 d^{\text{Riem}}_iW^n{}_m
+
L^n{}_{hi}W^h{}_m
-L^h{}_{mi}W^n{}_h,
\ee
\ses\\
which
involves the action of  the operator
\be
d^{\text{Riem}}_i=
\D{}{x^i}+L^k{}_i\D{}{y^k}.
\ee

{

In the  tangent Riemannian space
$\cR_{\{x\}}$
we can construct from the metric tensor
$g_{ij}=g_{ij}(x,y)$
the  curvature tensor
$
\wh R_{\{x\}}=\{\wh R_n{}^m{}_{ij}(x,y)\}
$
 by the help of  the ordinary Riemannian method,
regarding $\{y^i\}$ as variables.
Namely,
we obtain the representation
$$
\wh R_n{}^m{}_{ij}
=
\D{C^m{}_{ni}}{y^j}
-
\D{C^m{}_{nj}}{y^i}
 +
C^h{}_{ni}    C^m{}_{hj}
-
C^h{}_{nj}    C^m{}_{hi}.
$$
Since
$
\partial {C^m{}_{ni}}/\partial{y^j}
-
\partial{C^m{}_{nj}}/\partial{y^i}
\equiv
-2
\lf( C^h{}_{ni}    C^m{}_{hj} - C^h{}_{nj}    C^m{}_{hi} \rg),
$
we have simply
$$
\wh R_n{}^m{}_{ij}=
\fr1{F^2}
S_n{}^m{}_{ij},
$$
where
$$
S_n{}^m{}_{ij}
=
\lf(
C^h{}_{nj}    C^m{}_{hi}
-
C^h{}_{ni}    C^m{}_{hj}
\rg)
F^2.
$$
The tensor $S_n{}^m{}_{ij} $ describes the curvature of indicatrix (see Section 5.8 in [1]).

{

We need
 the {\it metrical  non-linear Finsler   connection}
 ${\mathcal  FN}$,
 such that
 \be
{\mathcal  FN}=\{N^m{}_i,T^m{}_{ij}\}: \qquad   N^m{}_i=N^m{}_i(x,y),
\qquad T^m{}_{ij}=T^m{}_{ij}(x,y),
\ee
where the objects
$N^m{}_i(x,y)$
and
$T^m{}_{ij}(x,y)
$
are to depend on the variable $y$ in an essentially non-linear way.
The adjective ``metrical'' means that the action of the entailed
covariant derivative  on the Finsler metric function,
 and also on the Finsler metric tensor, yields identically zero.
The coefficients $N^m{}_i$ and $T^m{}_{ij}$
are assumed to be positively homogeneous regarding the dependence on  vectors $y$,
respectively of degree 1 and degree 0.

{

In the Riemannian limit of the  Finsler space,
 the spaces
 ${\mathcal R}_{\{x\}}$
 are Euclidean spaces
and
the tensor
 $g_{\{x\}}(y)$ is independent of $y$.
 The conformally flat structure of the spaces
 ${\mathcal R}_{\{x\}}$
can naturally  be taken to treat as  the next  level of generality
of the Finsler space.
Can the metrical connection preserving the two-vector angle be introduced
on that level?

The  deformation of the Riemannian space to the Finsler space
proves to be
the  convenient method of consideration  to apply.
Namely,  when the Riemannian space
can be  deformated to  the Finsler space
characterized by
the  conformally flat structure of the spaces
 ${\mathcal R}_{\{x\}}$
the positive and clear answer to the above question can be arrived at.
The respective conformal multiplier is shown to be  a power of the Finsler
metric function.

We shall evidence the phenomenon that the used non-linear deformation
\be
{\mathcal  FN}=\bfC\cdot {\mathcal  RL}
\ee
of the Riemannian connection
 yields the Finsler connection
${\mathcal  FN}$ which preserves the Finslerian
two-vector angle
$ \al_{\{x\}}(y_1,y_2)$.

{

\ses

\ses

\setcounter{equation}{0}

\ses

{\bf II.2. Key observations}

\ses

\ses

Below, {\it  any dimension }  $N\ge 3$ is allowable.

\ses

Let $M$ be an $N$-dimensional
$C^{\infty}$
differentiable  manifold, $ T_xM$ denote the tangent space to $M$ at a point $x\in M$,
and $y\in T_xM\backslash 0$  mean tangent vectors.
Suppose we are given on the tangent bundle $TM$ a Riemannian metric
 $S$.
 Denote by
${\mathcal R}^N=(M,S)$
the obtained $N$-dimensional Riemannian space.
Let additionally a Finsler metric function  $F$ be introduced on this  $TM$,
yielding
a Finsler space
 ${\mathcal F}^N=(M,F)$.
We shall study the Finsler space
 ${\mathcal F}^N$
can be specified according to the following definition.

\ses

\ses

INPUT DEFINITION.
The  Finsler space ${\mathcal F}^N$ under consideration is the  {\it   deformated
 Riemannian space}
${{\mathcal R}}^N$:
\be
{\mathcal F}^N=\bfC\cdot {\mathcal R}^N,
\ee
specified by the condition  that in each tangent space $T_xM$
the metric tensor   $g_{\{x\}}(y)$
produced  by the Finsler metric
is  the $\bfC$-transformation of  the tensor which is {\it conformal}
to  the Euclidean metric tensor
 entailed by the Riemannian metric of the space
$ {\mathcal R}^N$.
It is assumed that  the applied $\bfC$-transformations (1.1)
 do not influence any point $x\in M$ of the base manifold  $M$
 and that they are  sufficiently   smooth
and invertible.
It is also natural to require that
the $\bfC$-transformations (1.1)
send  unit vectors to unit vectors:
\be
{\mathcal  IF}_{\{x\}}= \bfC \cdot    {\mathcal  S}_{\{x\}}.
\ee
Additionally, we subject the $\bfC$-transformation
 to the condition of positive homogeneity
with respect to tangent vectors $y$, denoting the degree of homogeneity
 by $H$.

\ses

{

If  $f(x,y)$ is the involved conformal multiplier
in the
tangent Riemannian space
$\cR_{\{x\}}$,
 then
 the equality
$
g_{\{x\}}(y)=f(x,y) u_{\{x\}}(y)
$
should introduce  the tensor $u_{\{x\}}(y)$  which associated
 Riemannian curvature tensor
vanishes identically.
The function $f(x,y)$ is assumed naturally to be homogeneous with respect to the argument
$y$.
Denoting
the
homogeneity degree
of $f(x,y)$
by means of $2a(x)$,
we just conclude that
the difference $1-a$ is exactly the homogeneity degree of the transformation (2.1) considered,
that is,
$
H=1-a.
$

The following proposition is valid.

\ses

\ses

{\bf Proposition II.2.1.}
{\it
A Finsler space  is  the
${\mathcal F}^N$-space if and only if  the indicatrix of the Finsler space
 is a space of constant curvature.
The dependence of the multiplier $f$ on the variable $y$
proves  to be presented by the  power of the Finsler metric function
$F$, such that
\be
g_{\{x\}}(y)=p^2u_{\{x\}}(y),  \quad p=c_1(x)\lf(F(x,y)\rg)^{a(x)},
  \quad c_1(x)>0.
\ee
The equality
}
$
{\mathcal  C}_{\text{Ind.}} = H^2
$
{\it ensues.

}

\ses

\ses

The proposition is of the local meaning in both the base manifold and the
tangent space.

\ses

\ses

{

{\bf Proof}. Given an arbitrary Finsler space of  any dimension  $N\ge3$.
The  tangent Riemannian  space
$\cR_{\{x\}}$
is conformally flat if and only if the  indicatrix
is a space of constant curvature.
Indeed,
in dimensions $N\ge4$ the conformal flatness holds if and only if the respective
Weyl tensor
$W_{ijmn}$ vanishes identically.
By evaluating the tensor and considering  the direct implications of
the contraction vanishing
$
W_{ijmn}l^nl^j=0,
$
we immediately obtain the representation
$ S_{nmij} =C(h_{nj}h_{mi}-h_{ni}h_{mj})$
which is characteristic of the constancy of the indicatrix curvature.
In the   dimension
$N\ge3$,
the conformal flatness of the space
$\cR_{\{x\}}$
is tantamount to the identical vanishing of
 the respective Cotton-York tensor.
  Considering the vanishing attentively leads again to the
representation
$ S_{nmij} =C(x)(h_{nj}h_{mi}-h_{ni}h_{mj})$.
These observations prove the first  part of   Proposition
II.2.1.
All the involved computations are explicitly represented in Appendix B.

To get the required conclusions concerning  the form of the respective conformal multiplier
we can  start with
the  tensor
$
u_{ij}
=
z(x,y)
(c_1(x))^{-2} F^{-2a(x)}
g_{ij},
$
where $z$ is a test smooth positive function
homogeneous of the  degree zero
with respect to the argument $y$.
We  evaluate
the respective curvature tensor
$
\wt  R_{\{x\}}
$
and assume
$\wt  R_{\{x\}}=0 $
to determine   the tensor
$
S_n{}^m{}_{ij}
=
\lf(
C^h{}_{nj}    C^m{}_{hi}
-
C^h{}_{ni}    C^m{}_{hj}
\rg)
F^2.
$
After that, we consider the implications of the vanishing
$
S_n{}^m{}_{ij}l^ml^j=0
$
and arrive at the representation
$$
S_{nmij}
=
a(2-a)
(h_{nj}h_{mi}-h_{ni}h_{mj})
+
F^2
 \fr1{2z^2}
  \bigl(z_h g^{hs} z_s\bigr)
(h_{nj}h_{mi}-h_{ni}h_{mj})
$$

\ses

$$
+  \fr {a-1}
{2z}
\Bigl(
z_n(
l_ih_{mj}-l_jh_{mi})
- z_m
(l_ih_{nj}-l_jh_{ni})
+
 l_n
(z_ih_{mj}-z_jh_{mi})
-
 l_m
(  z_ih_{jn} -z_jh_{in})
\Bigr)
F,
$$
where
 $z_k=\partial z/\partial y^k$.
The tensor
$
S_n{}^m{}_{ij}
$
must obviously possess the property
$
S_{nmij}l^i=0.
$
Therefore, we must fulfill the equation
$
(a-1)\bigl( z_{n} h_{mj} -
 z_{m} h_{nj}
\bigr)
=0.
$
Because of  $a\ne1$,
we can take only
 $z_n=0$, which means that the function
 $z$ is independent of $y$.
Without any loss  of generality  we can take $z=1$.
Thus we have proved  the second part in   Proposition II.2.1.
From the above representation of the tensor $S_{nmij}$
we just obtain
$
{\mathcal  C}_{\text{Ind.}} =1- a(2-a)\equiv (1-a)^2.
$
Since the difference $1-a$ is equal  to $H$,
the identification
$
{\mathcal  C}_{\text{Ind.}} = H^2
$
is valid.
 All the computations which are required to trace the validity of
the formulas   exposed
can be found in Appendix C.
~ Proposition II.2.1 is valid.
To have the equality
$S(x,\bar y)=\lf(F(x,y)\rg)^{H(x)}$, we make the choice
$c_1=1/H$.

{

\ses

Let the $\bfC$-transformation (I.2.1) proposed in Chapter I
be assigned locally by means of the differentiable functions
\be
\bar y^m=\bar y^m(x,y),
\ee
subject to the required   homogeneity
\be
\bar y^m(x,ky)=k^H\bar y^m(x,y), \qquad k>0, \forall y.
\ee
This entails the identity
\be
\bar y^m_ky^k=H\bar y^m,
\ee
where
$  \bar y^m_k=\partial \bar y^m/\partial y^k$.
Fulfilling  (2.1) means locally
\be
g_{mn}(x,y)=c_{ij}(x,\bar y) \bar y^i_m   \bar y^j_n, \qquad  c_{ij}(x,\bar y)=\lf(p(x,y)\rg)^2a_{ij}(x).
\ee
\ses\\
If we contract this tensor by $y^my^n$ and use the homogeneity identity
(2.6), we obtain the equality
\be
p(x,y)=\fr1{H(x)}\fr{F(x,y)}{S(x,\bar y)}.
\ee

{

On every punctured tangent space $T_xM\setminus0$,
the Finsler metric function $F$ is assumed to be positive,
 and also  positively homogeneous
of degree 1:
$$
F(x,ky)=kF(x,y), \quad  k>0, \forall y.
$$
The entailed Finsler metric tensor is
 positively homogeneous
of degree 0.
 Therefore, to comply the representation (2.7) with the stipulation (2.3),
we must put
\be
H=1-a.
\ee
With this observation, comparing (2.3) with (2.8) yields the equality
\be
c_1S=\fr1HF^H.
\ee
To comply with   the  indicatrix correspondence (2.2), we should put $c_1=1/H$,
which leads to the equality
 $S=F^H$      indicated in (I.2.3).

{

\ses

Denote
by
\be
y^i=y^i(x,t), \qquad t^n \equiv \bar y^n,
\ee
the inverse transformation, so that
$$
y^i(x,k t)=k^{1/H} y^i(x,t), \qquad k>0, \forall t,
$$
and
\be
y^i_nt^n=\fr1Hy^i,
\ee
where
$ y^i_n=\partial y^i/\partial t^n$.
The inverse to (2.7) reads:
\be
g_{kh} y^k_my^h_n  =  c_{mn}.
\ee

\ses

The following useful relations can readily be arrived at:
\ses\\
\be
y_my^m_n=\fr{F^2}{HS^2}t_n
\equiv
\fr1HF^{2(1-H)}t_n, \qquad t_n=a_{nh}t^h,
\ee
\ses\\
and
$$
y_my^m_{nl}t^l_j
+
g_{mj}y^m_n
=
2
\lf(\fr1H-1\rg)
F^{-2H}y_jt_n
+
\fr1H
  F^{2(1-H)}a_{nh}t^h_j,
$$
\ses
where
$t^l_j=\bar y^l_j$ and
$y^m_{nl}=\partial { y^m_n}/\partial { y^l}.
$
Alternatively,
\be
t_ht^h_n=\fr{HS^2}{F^2}y_n
\equiv
HF^{2(H-1)}y_n
\ee
and
\be
t_ht^h_{nu}y^u_i
+
a_{hi}t^h_n
=
2
(H-1)
F^{-2}t_iy_n
+
HF^{2(H-1)}g_{nu}y^u_i,
\ee
\ses\\
where $t^h_{nu}=\partial t^h_{n}/\partial y^u.$
We may  also write
\ses\\
\be
t_ht^h_{ni}
=
H(1-H)F^{2(H-1)}(g_{ni}-2l_nl_i).
\ee

\ses

From (2.13) it follows that
$$
g_{nm}y^m_i=p^2t^j_n a_{ij}, \qquad
y^k_{hp} t^p_n=- y^k_{p} t^p_{nv}y^v_h.
$$

{

Differentiating (2.7) with respect to $y^k$ yields the following
representation for the  tensor
$C_{mnk}=(1/2)\partial g_{mn}/\partial y^k$:
\be
2C_{mnk}=
(1-H)\fr2{F}l_k
g_{mn}
+p^2 (t^i_{mk}t^j_n +t^i_mt^j_{nk}) a_{ij}.
\ee
Contracting this tensor by $y^n$ results in the equality

\be
p^2 t^i_{mk}t^j  a_{ij}=
\lf(\fr1H-1\rg)(h_{km} -l_k  l_m ),
\ee
where the vanishing $C_{mnk}y^n=0$ and the homogeneity identity (2.6)
have
been taken into account.

\ses

Symmetry of the tensor $C_{mnk}$ demands
\be
(1-H)\fr2{F}(l_kg_{mn}-l_mg_{kn})
+p^2 (t^i_mt^j_{nk}-t^i_kt^j_{nm}) a_{ij}
=0,
\ee
so that we may alternatively write
\be
C_{mnk}=
(1-H)\fr1{F}(l_kg_{mn}+l_ng_{mk}-l_mg_{nk})
+p^2 t^i_mt^j_{nk} a_{ij}.
\ee

{

Contracting  the last tensor
by $g^{nk}$ yields
\ses\\
$$
2C_{m}=
(1-H)\fr2{F}l_m
+g^{nk}
p^2 (t^i_{nk}t^j_m +t^i_nt^j_{mk}) a_{ij}
\equiv 2C_{mnk}g^{nk},
$$
from which it ensues that
$$
2C_{m}=
(1-H)\fr2{F}l_m
+
2
g^{nk}
p^2 t^i_{nk}t^j_m  a_{ij}
+
g^{nk}
p^2 (t^i_nt^j_{mk}-t^i_mt^j_{nk})
 a_{ij},
$$
or
\ses\\
$$
2C_{m}=
(1-H)\fr2{F}l_m
+
2
g^{nk}
p^2 t^i_{nk}t^j_m  a_{ij}
-
(1-H)
g^{nk}
\fr2{F}(l_mg_{nk}-l_ng_{mk})
 a_{ij}.
$$
It is also convenient to use the representation
\ses\\
\be
FC_{m}=
-(N-2)(1-H)l_m
+
Fg^{nk}
p^2 t^i_{nk}t^j_m  a_{ij}.
 \ee

\ses

Since
$
y^m_i=p^2t^j_n a_{ij}g^{nm},
$
we can write
\ses\\
\be
FC^{m}=
-(N-2)(1-H)l^m
+
Fg^{nk}
 t^i_{nk}y^m_i.
\ee

{

 The  space ${\mathcal  F}^N$ is obtainable
  from
  the Riemannian
 space
 ${\mathcal  R}^N$
  by means of the  deformation
which, owing to  (2.7),
can be presented by the {\it conformal deformation tensor}
\be
{ C}^i_m~:=p\bar y^i_m,
\ee
so that
\be
 g_{mn}= { C}^i_m  { C}^j_n a_{ij}.
\ee
\ses
The zero-degree homogeneity
\be
{ C}^i_m(x,ky)=  { C}^i_m(x,y), \qquad k>0,\forall y,
\ee
holds,
together with
\be
{ C}^i_m(x,y)y^m=\lf(F(x,y)\rg)^{1-H}\bar y^i.
\ee

\ses

The indicatrix  correspondence  (2.2) is a direct implication of the
equality  $S=F^H$ .
We may apply the  transformation (1.1) to the unit vectors:
\be
l=\bfC\cdot L: ~ ~ ~  l^i=y^i(x,L);
~ \qquad
L=\bfC^{-1}\cdot l: ~ ~  ~  L^i=t^i (x,l),
\ee
where
$l^i=y^i/F(x,y)$ and $L^i=t^i/S(x,t)$ are components of the respective
Finslerian and Riemannian unit vectors,
which possess the properties $F(x,l)=1$ and $S(x,L)=1$.
We have $L^m=t^m(x,l)$.
On the other hand, from (2.7)  it just follows that
\be
g_{mn}(x,l)=\fr1{H^2}a_{ij}(x)t^i_m(x,l) t^j_n(x,l),
\ee
\ses
so that  under the transformation (2.28) we have
\be
g_{mn}(x,l)dl^mdl^n=\fr1{H^2}a_{ij}(x)dL^i dL^j.
\ee

{

\ses

\ses

{\bf Note.}
The deformation performed by
the formulas
 (2.24) and (2.25) is {\it unholonomic},
in the sense that
\be
\D{ C^i_m}{y^n}-\D{C^i_n}{y^m}\ne0.
\ee
The vanishing appears  if only the factor $p =F^{1-H}/H$ is independent of the vectors $y$,
that is, when $H=1$ (which is the Riemannian case proper).
Regarding the $y$-dependence,
the tensor  $C^i_m$ is homogeneous of degree zero, in accordance with (2.26).
If we divide the tensor by $p$, we obtain from (2.24) the tensor
$\bar y^i_m$ which is the derivative tensor, namely
$  \bar y^i_m=\partial \bar y^i/\partial y^m$. However, such a property cannot be addressed to
the tensor  $C^i_m$.
It is the reason why we start with the stipulation that the underlined
transformation
(which is downloaded locally
by the formulas (2.4)-(2.7))
be homogeneous of the degree $H$ with respect to the variable $y$.
By proceeding in this way, it proves possible to come to the
conformal representation (2.30) of $g_{mn}(x,l)dl^mdl^n$
which is of the key significance
to obtain the angle and the connection coefficients.

\ses

\ses

{

No support vector enters  the right-hand part of (2.30).
Therefore,
any two nonzero tangent vectors
  $y_1,y_2\in T_xM$ in  a fixed tangent space  $T_xM$
form the  ${\mathcal  F}^N${\it-space angle}
  \be
 \al_{\{x\}}(y_1,y_2) =
 \fr1{H(x)}
 \arccos\la,
  \ee
 where the scalar
\be
\la =
  \fr{a_{mn}(x)t_1^m  t_2^n}
 { S_1S_2 },
\qquad \text{with} \quad  t_1^m = t^m (x,y_1)  \quad   \text{and} \quad
t_2^m = t^m (x,y_2)    ,
   \ee
is of the entire Riemannian  meaning in the
space ${\mathcal  R}^N$;
the notation
$
S_1=\sqrt{a_{mn}(x)t^m_1t^n_1}
$
and
$
 S_2=\sqrt{a_{mn}(x)t^m_2t^n_2}
$
has been used.

From (2.33) it follows that
$$
\D{\la}{x^i}=
  \fr{a_{mn,i}t_{1}^m  t_2^n}
{S_1S_2 }+
\fr1{ S_1S_2}
  a_{mn}\lf(\D{t_{1}^m}{x^i}  t_2^n
  +
t_{1}^m  \D{t_2^n  }{x^i}
\rg)
$$

\ses

$$
-
\fr12\la
\lf[
\fr1{ S_1S_1 }
\lf(a_{mn,i}t_{1}^mt_1^n
+ 2 a_{mn} \D{t_{1}^m}{x^i} t_1^n
 \rg)
    +
\fr1{S_2S_2 }
\lf(a_{mn,i}t_{2}^mt_2^n
+ 2 a_{mn} \D{t_{2}^m}{x^i}  t_2^n
 \rg)
\rg],
$$
\ses
where   $a_{mn,i}=\partial a_{mn}/\partial x^i$,
and

$$
\D{\la}{y^k_1}=
\lf[
  \fr{a_{mn}  t_2^n}
{S_1S_2 }
 -      \fr{a_{mn} t_1^n}
{S_1S_1}\la
\rg]t_{1k}^m ,
\qquad
\D{\la}{y^k_2}=    \lf[
 \fr{a_{mn}  t_1^n}
{ S_2S_1 }
 -      \fr{a_{mn} t_2^n}
{S_2S_2}\la
\rg]t_{2k}^m.
$$

\ses

When the recurrent preservation
$$
d_i \al+ (1/H)H_i\al=0
$$
proposed by (I.1.28)
is applied to the angle given in (2.32),
we obtain simply
\be
d_i\la=0,
\ee
where $d_i$
is the separable operator (I.1.11).
That is,
the recurrent preservation law formulated for
the Finsler
  ${\mathcal  F}^N${\it-space angle}
(2.32) is tantamount to the separable
preservation law
(2.34) for the Euclidean angle
$ \arccos\la$.

{

We note also that
$$
A_1^k\D{\la}{y^k_1}=
F_1g_1^{nh}
 t_1{}^i_{nh}y_1{}^k_i
 \D{\la}{y^k_1}
 =
F_1g_1^{nh}
 t_1{}^i_{nh}
\lf[
  \fr{ t_{2i}}
{S_1S_2 }
 -      \fr{ t_{1i}}
{S_1S_1}\la
\rg].
$$

\ses

{

\setcounter{equation}{0}

\ses

\ses

 {\bf II.3.  Derivation and properties of  the  coefficients
$ N^{m}{}_n$
in the  ${\mathcal  F}^N$-space}

\ses

\ses

Let us start from (2.11) and introduce
  the vector
$U^i=U^i(x,y)$
according to
\ses
\be
U^i\eqdef \fr1St^i \equiv \fr1{F^H}t^i,
\ee
\ses\\
which is obviously unit:
\be
U_iU^i=1,   \qquad    U_i=a_{ij}U^j.
\ee
\ses\\
The zero-degree homogeneity
\be
U^i(x,k y)= U^i(x,y), \qquad k>0, \forall y,
\ee
holds, entailing the identity
\be
U^i_ny^n=0,
\ee
where
\be
 U^i_n~:=\D{ U^i}{ y^n}
 =
\fr1{F^H} t^i_n-\fr1FH U^il_n.
 \ee

From (2.14) it follows that
 \be
F^H U^h_sy^k_{h}
 =h^k_s, \qquad
F^H U^i_ky^k_t
 =\de^i_t -U^iU_t, \qquad
 U_i U^i_n=0.
\ee

The vanishing
\be
U_i\lf( \D{U^i}{x^n}+L^i{}_{kn} U^k\rg)=0
\ee
 holds obviously,
 where
 $L^i{}_{kn}=L^i{}_{kn}(x)$ are the Riemannian connection coefficients appeared in (1.2).

{

\ses

The representation  (2.33) takes on the simple form
\be
\la =
a_{mn}(x)U_1^m  U_2^n,
\ee
\ses\\
with
\be
U_1^m = U^m (x,y_1), \qquad
U_2^m = U^m (x,y_2).
\ee

\ses

{

The form of the rght-hand part in the formula (3.8) which
represents the scalar $\la$ is such that
the
preservation law
$d_i\la=0$
written in (2.34)
 is obviously equivalent to the
vanishing
\be
{\mathcal D}_nU^i=0
\ee
for the field $U^i=U^i(x,y)$,
with
 the covariant derivative
\be
{\mathcal D}_nU^i~:=d_nU^i+L^i{}_{nk} U^k.
\ee
Since
$$
d_nU^i=\D{U^i}{x^n}+N^k{}_nU^i_k,
$$
we obtain the representation
\be
N^{m}{}_n
=-y^m_iF^H
 \lf(\fr HFU^i\D{F}{x^n}+
 \D{U^i}{x^n}+L^i{}_{nk} U^k\rg)
 +
l^m
d_nF
\ee
\ses\\
which was indicated in (I.2.15).

{

We have arrived at the following proposition.

\ses

\ses

{\bf Proposition II.3.1.}
{\it Given  an arbitrary smooth function $H(x)$,
the angle preservation  equation
 $d_n\al+(1/H)H_n\al=0$
in the  ${\mathcal  F}^N$-space
entails  the representation}
 (3.12)
 {\it  for the coefficients
}
$ N^{m}{}_n.$

\ses

\ses

 {

By differentiating (3.10) with respect to $y^m$
 we may conclude that
 the covariant derivative
\be
{\mathcal D}_n U^i_m ~:= d_n U^i_m-D^h{}_{nm}U_h^i+L^i{}_{nl}U^l_m, \qquad D^h{}_{nm}=-N^h{}_{nm},
\ee
vanishes identically:
\be
{\mathcal D}_n U^i_m= 0.
\ee

\ses

{

Below, we shall assume that
$$d_nF=0.$$

\ses

Using
$t^i=F^H U^i$
together with
\be
{\mathcal D}_nt^i~:=d_nt^i+L^i{}_{kn}t^k,
\ee
from (3.10) we find
\be
{\mathcal D}_nt^i=t^iH_n\ln F.
\ee

Differentiating (3.16) with respect to $y^m$
leads to the conclusion  that
 the covariant derivative
\be
{\mathcal D}_n t^i_m ~:= d_n t^i_m-D^h{}_{nm}t_h^i+ L^i{}_{nl}t^l_m
\ee
possesses the property
\be
{\mathcal D}_n t^i_m=
\lf(
t^i_m
\ln F
+
t^il_m\fr1 {F}
\rg)
H_n.
\ee

{

\ses

With
$
p=
(1/H)
F^{1-H}
$
from (3.18)
we get
\be
{\mathcal D}_n (pt^i_m)=
p
t^il_m\fr1 {F}
H_n
-\fr1HH_npt^i_m,
\ee
so that,
\be
{\mathcal D}_n (pt^i_m)=
-\fr pHH_nh^k_mt^i_k.
\ee

Consider (2.7):
$$
 g_{mn}(x,y)=p^2t^i_m  t^j_na_{ij}(x).
$$
We  obtain
$$
{\mathcal D}_k g_{mn}=
p^2
\lf(
t^il_m\fr1 {F}
-\fr1Ht^i_m
\rg)H_k
 t^j_na_{ij}
+
p^2
t^i_m
\lf(
t^jl_n\fr1 {F}
-\fr1Ht^j_n
\rg)H_k
 a_{ij}.
$$
\ses
Using
$
t_ht^h_n
=
HF^{2(H-1)}y_n
$
(see (2.34))
leads to
\ses\\
$$
{\mathcal D}_k g_{mn}= -\fr2HH_kg_{mn}
+
p^2
l_m
H_k
HF^{2(H-1)}l_n
+
p^2
t^i_m
t^jl_n\fr1 {F}
H_k
 a_{ij}.
$$

In this way we arrive at the following result
after the direct evaluations performed.

\ses

\ses

{\bf Proposition II.3.2.}
{\it Given  an arbitrary smooth function $H(x)$
in the  ${\mathcal  F}^N$-space,
the angle preservation
 $d_n\al+(1/H)H_n\al=0$
taken in conjunction with the preservation $d_nF=0$ of the metric function
entails that
 the covariant derivative of the metric tensor
 reads}
\ses\\
\be
{\mathcal D}_k g_{mn}
=
-
\fr2HH_k
 h_{mn}.
\ee

\ses

\ses

{

Now, we contract (3.19)  by $y^m_k$, getting
$$
y^m_k{\mathcal D}_n t^i_m=
y^m_k d_n t^i_m
-
D^h{}_{nm}y^m_kt_h^i
+
L^i{}_{nk}
=
\lf(
\de^i_k
\ln F
+
t^it_k
\fr 1{HF^{2H}}
\rg)
H_n.
$$

\ses

Since
$ y^n_k t^k_j=\de^n_j$,
the previous identity can be transformed to
\ses\\
$$
t^i_md_ny^m_k
+
D^h{}_{nm}y^m_kt_h^i
-
L^i{}_{nk}
=
-
\lf(
\de^i_k
\ln F
+
t^it_k
\fr 1{HF^{2H}}
\rg)
H_n.
$$

Contract this equality by $y^j_i$, obtaining
the equality
\ses\\
\be
 d_n
y^j_k
+
D^j{}_{nm}y^m_k
-
L^i{}_{nk}
y^j_i
=
-
\lf(
y^j_k
\ln F
+
y^jt_k
\fr 1{H^2F^{2H}}
\rg)
H_n,
\ee
\ses\\
which can be written simply as
\be
 d_n
\lf(\fr1p  y^j_k\rg)
+
T^j{}_{nm}y^m_k  \fr1p
-
L^i{}_{nk}
y^j_i\fr1p
=
0,
\ee
\ses
where
$T^j{}_{nm}$ are the coefficients introduced in (I.1.33).
Taking into account the representations
(3.15)-(3.16) together with the identity
$$
t^m\D{\fr1p  y^j_k}{t^m}=0
$$
ensuing from the homogeneity,
the equality (3.23) becomes
\be
 d^{\text{Riem}}_n
\lf(\fr1p  y^j_k\rg)
+
T^j{}_{nm}y^m_k  \fr1p
-
L^i{}_{nk}
y^j_i\fr1p
=
0.
\ee
\ses\\
We  have
used the Riemannian operator
$
d^{\text{Riem}}_n
$
introduced in (II.1.4).

{

We know that
\be
 d^{\text{Riem}}_n
t^k
+
L^k{}_{nm}t^m
=
0.
\ee
Therefore, contracting (3.24) by $t^k$ yields
\be
 d^{\text{Riem}}_n
\lf(  \fr1{Hp} y^j\rg)
-
N^j{}_{n}  \fr1{Hp}
=
0.
\ee
\ses\\
Here we have
\be
Hp=F^{1-H}=S^{1/H}\fr1S,
\ee
 so that
\be
 d^{\text{Riem}}_n
 \fr1{Hp}
=
 \fr1{Hp}
\fr1{H^2}
 H_n
\ln S
\equiv
\fr1{Hp}
\fr1{H}
H_n
\ln F.
\ee

\ses

We arrive at the following proposition.

\ses

\ses

{\bf Proposition II.3.3.}
{\it With   an arbitrary smooth function} $H(x)$,
{\it
in the
} ${\mathcal  F}^N$-{\it  space
with $d_nF=0$
 the representation
}
\be
N^m{}_{n}
=
d^{\text{Riem}}_n  y^m(x,t)+  \fr1H H_ny^m\ln F
\ee
{\it written by the help of   the Riemannian operator
}
$d^{\text{Riem}}_n$
{\it is valid.
}

\ses

\ses

{

\ses

The derivative coefficients
$N^k{}_{mn}=\partial{N^k{}_{m}}/\partial{y^n}$
can straightforwardly be evaluated from the coefficients
$N^k{}_m$ written in (I.2.16).
We obtain
$$
N^k{}_{mn}
=
-\fr1F
h^k_n\D{F}{x^m}
-
l^k\D{l_n}{x^m}
-
\lf(y^k_{hp}t^p_n+\fr HFy^k_hl_n\rg)
F^H \lf( \D{U^h}{x^m}+ L^h{}_{ms}U^s\rg)
$$

\ses

\ses

$$
-
y^k_hF^H
 \lf( \D{ U^h_n}{x^m}+ L^h{}_{ms}U^s_n\rg).
 $$
\ses\\
Owing to (3.10)),
we have
$$
\D{U^h}{x^m}+ L^h{}_{ms}U^s=-N^s{}_mU^h_s,
$$
\ses\\
so that using (3.6) we observe that the coefficients
$N^k{}_{mn}$ are equal to
\ses\\
$$
-\fr1F
h^k_n\D{F}{x^m}
-
l^k\D{l_n}{x^m}
+
U^h_sy^k_{hp}
t^p_n
F^H N^s{}_m
+
\fr HFy^k_hl_n
F^H N^s{}_mU^h_s
-
y^k_hF^H
 \lf( \D{ U^h_n}{x^m}+ L^h{}_{ms}U^s_n\rg),
 $$
\ses\\
or
$$
N^k{}_{mn}
=
-
\fr1F
h^k_n\D{F}{x^m}
-
l^k\D{l_n}{x^m}
-
h^v_s y^k_{j} t^j_{nv} N^s{}_m
+
\fr HF h^k_sl_n
 N^s{}_m
-
y^k_hF^H
 \lf( \D{ U^h_n}{x^m}+ L^h{}_{ms}U^s_n\rg),
 $$
\ses\\
where the relation
 $$
F^H  U^h_sy^k_{hp}
t^p_n
=
-
F^H  U^h_s
 y^k_{p} t^p_{nv}y^v_h
=
-
h^v_s y^k_{p} t^p_{nv}
$$
has been used.

{

From (2.21) we have
$$
p^2 t^i_mt^j_{nk} a_{ij}
=
C_{mnk}
-(1-H)\fr1{F}(l_kg_{mn}+l_ng_{mk}-l_mg_{nk}),
$$
which is
\ses\\
$$
 g_{mv}y^v_jt^j_{nk}
=
C_{mnk}
-(1-H)\fr1{F}(l_kg_{mn}+l_ng_{mk}-l_mg_{nk}).
$$
We obtain
\ses\\
\be
y^k_{j} t^j_{nv}
=C^k{}_{nv}
-(1-H)\fr1{F}(l_v\de^k_n+l_n\de^k_v-l^kg_{nv})
\ee
\ses
and
\be
 h^v_sy^k_{j} t^j_{nv}
=C^k{}_{ns}
-(1-H)\fr1{F}(l_nh^k_s-l^kh_{ns}).
\ee

{

In this way we come to the representation
\ses\\
$$
N^k{}_{mn}
=
-
\fr1F
h^k_n\D{F}{x^m}
-
l^k\D{l_n}{x^m}
-
\lf( C^k{}_{ns}
-(1-H)\fr1{F}(l_nh^k_s-l^kh_{ns})
\rg)
 N^s{}_m
+
\fr HF h^k_sl_n
 N^s{}_m
$$

\ses

\ses

$$
-
y^k_hF^H
 \lf( \D{ U^h_n}{x^m}+ L^h{}_{ms}U^s_n\rg).
 $$

The eventual  result reads
\ses\\
$$
N^k{}_{mn}
=
-
\fr1F
h^k_n\D{F}{x^m}
-
l^k\D{l_n}{x^m}
-
 C^k{}_{ns}
  N^s{}_m
+
\fr1{F}
\lf(
l_nh^k_s
-
(1-H)
l^kh_{ns}
\rg)
 N^s{}_m
$$

\ses

\ses

\be
-
y^k_hF^H
 \lf( \D{ U^h_n}{x^m}+ L^h{}_{ms}U^s_n\rg).
\ee

\ses

Thus we can formulate the following assertion.

\ses

\ses

{\bf Proposition II.3.4.}
{\it With   an arbitrary smooth function} $H(x)$,
{\it
in the
} ${\mathcal  F}^N$-{\it  space
with $d_nF=0$
 the  coefficients
$N^k{}_{mn}$
can be given by means of the explicit representation written in}
(3.32).

\ses

\ses

We are also able to evaluate the entailed coefficients
$N^k{}_{mni}=\partial{N^k{}_{mn}}/\partial{y^i}$.
The required evaluations which have been presented in detail in Appendix D
lead straightforwardly
to
 the representation
\be
N^k{}_{mni}
=
\fr2HH_m
\fr1Fl^k
h_{ni}
-
\cD_m C^k{}_{ni}.
\ee

Thus we can formulate the following assertion.

\ses

\ses

{\bf Proposition II.3.5.}
{\it Given  an arbitrary smooth function} $H(x)$,
{\it
in the
} ${\mathcal  F}^N$-{\it  space
with $d_nF=0$
 the  coefficients
$N^k{}_{mni}$
admit the simple representation}
 (3.33)
 {\it in terms of the covariant derivative of the tensor
$C^k{}_{ni}$.}

\ses

\ses

{

\ses

\setcounter{equation}{0}

 {\bf II.4.  Properties of  covariant derivative   }

\ses

\ses

The equality
 (3.20)
can be written in the form
\be
\cT_i (pt^{m}_n)=0
\ee
with
\be
\cT_i (pt^{m}_n)= d_i (pt^{m}_n) - T^h{}_{in} pt^{m}_h
+
L^m{}_{il}pt^{l}_n,
\ee
where
$T^h{}_{in}$ are the coefficients (I.1.24).
If we contract the last vanishing
by $y^n$
and note
that
 $\cT_i y^n=0$
 (see (I.1.37)), we get
\be
\cT_i (Hpt^m)=0,
\ee
where
\be
\cT_i (Hpt^m)=d_i(Hpt^m)
+
L^m{}_{ik}Hpt^k.
\ee

\ses

We may write
\be
\cT_i (Hpt^m)= t^md_i(Hp)+Hp\cD_it^m,
\ee
\ses\\
where
$\cD_it^m=
d_i t^m+L^m{}_{ik} t^k$.

{

Owing to  the equality
$C^{m}_n=pt^{m}_n$
(see (2.24)),
from (4.1)
we are  entitled to formulate the following proposition.

 \ses

\ses

{\bf Proposition II.4.1.}
{\it
The  $\bfC$-deformation  is}
$\cT$-{covariant constant:}
\be
{\cT}\cdot{\bfC}=0.
\ee

\ses

\ses

In terms of local coordinates the previous vanishing reads
\be
\cT_n {C}^m_k=0,
\ee
 where
\be
\cT_n {{} C}^m_k= d_n {C}^m_k - T^h{}_{nk}{C}_h^m
+
L^m{}_{nl}
C^l_k.
\ee

\ses

The reciprocal coefficients
\be
\wt  C^n_m=\fr1p y^n_m
\ee
\ses\\
fulfills the similar vanishing
\be
\cT_n {\wt C}^m_k=0,
\ee
\ses\\
where
\be
\cT_n {\wt C}^m_k=
 d^{\text{Riem}}_n
{\wt C}^m_k
+
T^m{}_{nh}{\wt C}^h_k
-
L^i{}_{nk}
{\wt C}^m_i
\ee
\ses\\
(see (3.23)).

{

Let us realize the action of the $\cC$-transformation (2.1)-(2.2)
on tensors  by the help of the deformaton
\be
\{w(x,y)\}  =\cC\cdot \{W(x,t)\},
\ee
\ses\\
assuming that the tensors
$\{w(x,y)\}  $ are positively homogeneous of degree 0 with respect to the variable $y$,
and that  the tensors
$\{W(x,t)\}  $ are positively homogeneous of degree 0 with respect to the variable $t$.
Namely, in the scalar case
we use the identification
\be
w(x,y)  =W(x,t),
\ee
 obtaining  merely
\be
d_iw= d^{\text{Riem}}_iW
\ee
(because of
the vanishing
$\cD_i U^j=0$ indicated in (3.10)-(3.11)),
where
$d^{\text{Riem}}_i$
is  the operator defined by   (II.1.4).
Given a tensor $w_n(x,y)$ of the type (0,1)
we use  the transformation
\be
 w_n=  C^m_n W_m.
\ee
The  metrical  linear  connection
${\mathcal  RL}$ introduced by (1.2)
 may  be used to  define the
   covariant derivative  $\nabla$ in ${\cal R}^N$
according to the  conventional  rule:
\be
\nabla_iW_m=\D{W_m}{x^i}+L^k{}_i   \D{W_m}{t^k}
-
L^h{}_{im}W_h,  \qquad    L^k{}_j=- L^k{}_{ij}t^i,
   \ee
which can be written shortly
\be
\nabla_iW_m= d^{\text{Riem}}_iW_m
-
L^h{}_{im}W_h.
\ee
We have
\be
\nabla_iS=0, \qquad \nabla_i t^j=0, \qquad  \nabla_ia_{mn}=0.
\ee

By virtue of  the nullification $\cT_iC^m_n=0$
shown in (4.7),
 we obtain
 the {\it transitivity property}
\be
\cT_iw_n=C^m_n\nabla_iW_m
\ee
for the covariant derivatives.

{

The method can be repeated in  case
of the covariant vectors
$w^n(x,y)$  and $W^n(x,t)$,
 namely we write
\be
 w^n= \wt  C^n_m W^m,
\ee
obtaining
\be
\cT_iw^n=\wt C^n_m\nabla_iW^m,
\ee
where
the reciprocal coefficients
$
\wt  C^n_m=(1/p) y^n_m
$
defined by (4.9) have been arisen.

{

The method can also
be extended to more general tensors in a direct manner.
For example, considering the (1,1)-type
tensors
$\{w^n_m(x,y), W^n_m(x,t)\}$
of the zero-degree positive
homogeneity with respect to the variables $y$ and $t$,
we can use the  covariant
derivative
\be
\nabla_iW^n{}_m=\D{W^n{}_m}{x^i}+L^k{}_i   \D{W^n{}_m}{t^k}
+
L^n{}_{hi}W^h{}_m
-L^h{}_{mi}W^n{}_h
\equiv
 d^{\text{Riem}}_iW^n{}_m
+
L^n{}_{hi}W^h{}_m
-L^h{}_{mi}W^n{}_h
\ee
\ses\\
and the deformation
\be
w^n{}_m=\wt C^n_h C^k_mW^h{}_k
\ee
to obtain  the transitivity property
\be
\cT_iw^n{}_m=\wt C^n_h C^k_m\nabla_iW^h{}_k
\ee
for the covariant derivatives
$\cT_i$ and
$\nabla_i$.

{

\ses

Now we may formulate
 the following proposition.

 \ses

\ses

{\bf Proposition II.4.2.}
{\it
The covariant derivative $\cT$
is  the
manifestation of the
{\it transitivity}
of the connection  under   the $\cC$-transformation.
}

\ses

\ses

In short,
\ses\\
\be
\cT=\cC\cdot\nabla.
\ee

{

\setcounter{equation}{0}

\ses

\ses

\ses

   {\bf II.5. Entailed curvature tensor}

\ses

\ses

Henceforth, the torsion tensor $S^m{}_{ij}$  (entered the initial connection (1.2))
is not accounted for.

\ses

Given a tensor $w^n{}_k=w^n{}_k(x,y)$ of the tensorial type (1,1),
commuting the covariant derivative
\be
\cT_iw^n{}_k~:=   d_iw^n{}_k+T^n{}_{ih}w^h{}_k-T^h{}_{ik}w^n{}_h
\ee
 yields the equality
\be
\lf[\cT_i\cT_j-\cT_j\cT_i\rg] w^n{}_k=M^h{}_{ij}\D{w^n{}_k}{y^h}  -E_k{}^h{}_{ij}w^n{}_h
+E_h{}^n{}_{ij}w^h{}_k
\ee
with the  tensors
\be
 M^n{}_{ij}~:= d_iN^n{}_j-d_jN^n{}_i
\ee
and
\be
E_k{}^n{}_{ij}~: =
d_iT^n{}_{jk}-d_jT^n{}_{ik}+T^m{}_{jk}T^n{}_{im}-T^m{}_{ik}T^n{}_{jm}.
\ee

\ses

By applying the commutation rule (5.2)
to the particular choices $\{F, y^n, y_k,g_{nk}\}$
and noting the vanishing
$\{ {\cT}_i F =  {\cT}_i y^n=  {\cT}_i y_k = {\cT}_ig_{nk}=0\}$,
we obtain the identities
\be
y_n M^n{}_{ij}=0,
\qquad
y^kE_k{}^n{}_{ij}  = -   M^n{}_{ij},   \qquad
y_nE_k{}^n{}_{ij}  = M_{kij},
\ee
and
\be
E_{mnij}+E_{nmij}=2C_{mnh}
M^h{}_{ij} \quad \text{with} \quad C_{mnh}=\fr12\D{g_{mn}}{y^h}.
\ee

{

It proves pertinent to replace in the commutator (5.2)
the partial derivative
$\partial  w^n{}_k/\partial y^h$
by the definition
\be
\cS_h{w^n{}_k}~:=\D{w^n{}_k}{y^h}+C^n{}_{hs}w^s{}_k-C^m{}_{hk}w^n{}_m
\ee
which has the meaning of the covariant derivative in the tangent
  space
supported by the point $x\in M$.
In particular,
$$
\mathcal S_h g_{nk}~:=\D{g_{nk}}{y^h} -C^m{}_{hn}g_{mk} -C^m{}_{hk}g_{nm}=0.
$$
With
the {\it curvature tensor}
\be
\Rho_k{}^n{}_{ij}~:=E_k{}^n{}_{ij}
-
M^h{}_{ij}C^n{}_{hk},
\ee
the commutator (5.2) takes on the form
\be
\lf(\cT_i\cT_j-\cT_j\cT_i\rg) w^n{}_k=
M^h{}_{ij}\cS_h w^n{}_k  -\Rho_k{}^h{}_{ij}w^n{}_h
+\Rho_h{}^n{}_{ij}w^h{}_k.
\ee

\ses

We denote
$\Rho_{knij}=
g_{mn}\Rho_k{}^m{}_{ij}.$
The skew-symmetry
\be
 \Rho_{mnij}=-\Rho_{nmij}
\ee
 holds (cf.  (5.6)).

\ses

The equalities
\be
y^k \Rho_k{}^n{}_{ij}  = -   M^n{}_{ij},   \qquad
y_n \Rho_k{}^n{}_{ij}  = M_{kij},
\ee
obviously hold.

{

Let us evaluate the tensor
$
 M^n{}_{ij},
$
using the coefficients
$$
N^{n}{}_j
=
-l^n\D{F}{x^j}
-
y^n_h
F^H
 \lf( \D{U^h}{x^j}+ L^h{}_{kj}U^k\rg)
$$
indicated in (I.2.16).


We directly obtain
\ses\\
$$
 M^n{}_{ij}
 =
 -
 \fr1FN^n{}_i\D{F}{x^j}
+
 \fr1FN^n{}_j\D{F}{x^i}
 -
 l^nN^s{}_i\D{l_s}{x^j}
+
l^nN^s{}_j\D{l_s}{x^i}
 $$

\ses

\ses

$$
-
d_i \lf(y^n_hF^H\rg)
 \lf( \D{U^h}{x^j}+ L^h{}_{kj}U^k\rg)
+
d_j \lf(y^n_hF^H\rg)
 \lf( \D{U^h}{x^i}+ L^h{}_{ki}U^k\rg)
$$

\ses

\ses

$$
-
y^n_h
F^H
\lf[
\lf(
\D{ L^h{}_{kj}}{x^i}
-
\D{ L^h{}_{ki}}{x^j}
\rg)
U^k
+
 L^h{}_{kj}\D{U^k}{x^i}
-
 L^h{}_{ki}\D{U^k}{x^j}
 \rg]
$$

\ses

\ses

\ses

\ses

$$
-
y^n_h
F^H
\lf[
N^s{}_i
 \lf( \D{U^h_s}{x^j}+ L^h{}_{kj}U^k_s\rg)
-
N^s{}_j
 \lf( \D{U^h_s}{x^i}+ L^h{}_{ki}U^k_s\rg)
 \rg].
$$

\ses

Using here    the     Riemannian curvature tensor
\be
a_k{}^h{}_{ij}=
\D{ L^h{}_{kj}}{x^i}-\D{ L^h{}_{ki}}{x^j}+ L^u{}_{kj} L^h{}_{ui}- L^u{}_{ki} L^h{}_{uj}
\ee

 {

\nin
leads to
 \ses\\
$$
 M^n{}_{ij}
 =
 -
 \fr1FN^n{}_i\D{F}{x^j}
+
 \fr1FN^n{}_j\D{F}{x^i}
 -
 l^nN^s{}_i\D{l_s}{x^j}
+
l^nN^s{}_j\D{l_s}{x^i}
 $$

\ses

\ses

$$
-
d_i \lf(y^n_hF^H\rg)
 \lf( \D{U^h}{x^j}+ L^h{}_{kj}U^k\rg)
+
d_j \lf(y^n_hF^H\rg)
 \lf( \D{U^h}{x^i}+ L^h{}_{ki}U^k\rg)
$$

\ses

\ses

$$
-
y^n_h
F^H
\lf[
\lf(
a_k{}^h{}_{ij}
-
 L^u{}_{kj} L^h{}_{ui}
+
 L^u{}_{ki} L^h{}_{uj} \rg)
U^k
+
 L^h{}_{kj}\D{U^k}{x^i}
-
 L^h{}_{ki}\D{U^k}{x^j}
 \rg]
$$

\ses

\ses

\ses

$$
-
y^n_h
F^H
\lf[
N^s{}_i
 \lf( d_jU^h_s+ L^h{}_{kj}U^k_s\rg)
-
N^s{}_j
 \lf( d_iU^h_s+ L^h{}_{ki}U^k_s\rg)
 \rg].
$$

\ses

\ses

Now, we apply the vanishing
$
d_iU^h+ L^h{}_{is}U^s=0
$
(see (3.10)-(3.11)),
getting
\ses\\
$$
 M^n{}_{ij}
 =
 -
 \fr1FN^n{}_i\D{F}{x^j}
+
 \fr1FN^n{}_j\D{F}{x^i}
 -
 l^nN^s{}_i\D{l_s}{x^j}
+
l^nN^s{}_j\D{l_s}{x^i}
 $$

\ses

\ses

$$
+
d_i \lf(y^n_hF^H\rg)
N^t{}_jU^h_t
-
d_j \lf(y^n_hF^H\rg)
N^t{}_iU^h_t
$$

\ses

\ses

$$
-
y^n_h
F^H
\lf(
a_k{}^h{}_{ij}
U^k
-
 L^h{}_{kj}
N^t{}_iU^k_t
+
 L^h{}_{ki}
N^t{}_jU^k_t
 \rg)
$$

\ses

\ses

\ses

$$
-
y^n_h
F^H
\lf[
N^t{}_i
 \lf( d_jU^h_t+ L^h{}_{kj}U^k_t\rg)
-
N^t{}_j
 \lf( d_iU^h_t+ L^h{}_{ki}U^k_t\rg)
 \rg].
$$

\ses

Taking into account the equality
$$
F^H U^h_ty^k_{h}
 =h^k_t
 $$
 (see (3.6)),
we arrive at
the representation

{

$$
 M^n{}_{ij}
 =
 -
 \fr1FN^n{}_i\D{F}{x^j}
+
 \fr1FN^n{}_j\D{F}{x^i}
 -
 l^nN^s{}_i\D{l_s}{x^j}
+
l^nN^s{}_j\D{l_s}{x^i}
 $$

\ses

\ses

$$
-
y^n_h
F^H
a_k{}^h{}_{ij}
U^k
-
N^t{}_i
d_jh^n_t
+
N^t{}_j
d_jh^n_t,
$$
\ses\\
which can readily be simplified to read
\ses
$$
 M^n{}_{ij}
 =
 \fr1FN^n{}_iN^s{}_jl_s
-
 \fr1FN^n{}_jN^s{}_il_s
 -
 l^nN^s{}_i\D{l_s}{x^j}
+
l^nN^s{}_j\D{l_s}{x^i}
 $$

\ses

\ses

$$
-
y^n_h
F^H
a_k{}^h{}_{ij}
U^k
+
N^s{}_i
d_j(l^nl_s)
-
N^s{}_j
d_j(l^nl_s)
$$

\ses

\ses

 \ses

$$
 =
 -
 l^nN^s{}_i\D{l_s}{x^j}
+
l^nN^s{}_j\D{l_s}{x^i}
-
y^n_h
F^H
a_k{}^h{}_{ij}
U^k
+
l^nN^s{}_i
d_jl_s
-
l^nN^s{}_j
d_jl_s.
$$

{

The eventual result is
\ses\\
\be
M^n{}_{ij}=
-
y^n_t t^h
a_h{}^t{}_{ij}.
\ee

{

Next, we use the equaity
\ses\\
$$
\cT_i\cT_jw^n{}_m=\wt C^n_h C^k_m\nabla_i\nabla_jW^h{}_k
$$
(see  (4.24))
to consider the relation
$$
[\cT_i\cT_j-\cT_j\cT_i]w^n{}_m
=
\wt C^n_h C^k_m
[\nabla_i\nabla_j-\nabla_j\nabla_i]W^h{}_k.
$$

In the commutator
\be
\lf[\nabla_i\nabla_j-\nabla_j\nabla_i\rg] W^n{}_k
=
-t^sa_s{}^h{}_{ij}  \D{W^n{}_k}{t^h}
-
a_k{}^h{}_{ij}W^n{}_h
+
a_h{}^n{}_{ij}W^h{}_k
\ee
the  Riemannian curvature tensor $a_s{}^h{}_{ij}$
is constructed in accordance with the ordinary rule (5.12).
Whence,
\ses\\
$$
M^h{}_{ij}\D{w^n{}_m}{y^h}  -E_m{}^h{}_{ij}w^n{}_h
+E_h{}^n{}_{ij}w^h{}_m
=
\wt C^n_h C^k_m
\lf[
-t^sa_s{}^r{}_{ij}  \D{W^h{}_k}{t^r}
-
a_k{}^r{}_{ij}W^h{}_r
+
a_r{}^h{}_{ij}W^r{}_k
\rg].
$$

Using here the equality
$$
W^h{}_k= C^h_n \wt C^m_kw^n{}_m
$$
(taken from (4.23))
leads to
\ses\\
$$
M^h{}_{ij}\D{w^n{}_m}{y^h}  -E_m{}^h{}_{ij}w^n{}_h
+E_h{}^n{}_{ij}w^h{}_m
=
-
t^sa_s{}^h{}_{ij}  \D{w^n{}_m}{t^h}
+
W^h{}_k
t^sa_s{}^r{}_{ij}  \D{y^n_h t^k_m}{t^r}
$$

\ses

\ses

$$
-
t^k_m a_k{}^r{}_{ij} y^h_r
w^n{}_h
+
y^n_s   a_r{}^s{}_{ij}t^r_h
w^h{}_m
$$

Now we use here the representation (5.13) obtained for the tensor
$M^h{}_{ij}$.
We are left with
\ses\\
$$
 -E_m{}^h{}_{ij}w^n{}_h
+E_h{}^n{}_{ij}w^h{}_m
=
t^h_u y^v_kw^u{}_v
t^sa_s{}^r{}_{ij}  \D{y^n_h t^k_m}{t^r}
-
t^k_m a_k{}^r{}_{ij} y^h_r
w^n{}_h
+
y^n_s   a_r{}^s{}_{ij}t^r_h
w^h{}_m.
$$

{

In this way we obtain the explicit representation
\ses\\
\be
E_k{}^n{}_{ij}
=
y_h^n t^h_{km}M^m{}_{ij}
+
y^n_ma_h{}^m{}_{ij} t^h_k.
\ee

{

From (5.8) and (5.15) it follows that
$$
\Rho_k{}^n{}_{ij}
=
\lf(
y_h^n t^h_{km}
-C^n{}_{mk}
\rg)
M^m{}_{ij}
+
y^n_ma_h{}^m{}_{ij} t^h_k.
$$
Inserting    here the tensor $C^n{}_{mk}$
taken from (2.21)
and noting
 the vanishing $l_mM^m{}_{ij}=0$
(see (5.5)), we get

$$
\Rho_k{}^n{}_{ij}
=
\lf(
y_h^n t^h_{km}
-
(1-H)\fr1{F}(l_k\de^n_m+l^ng_{mk})
-p^2 t^l_mt^h_{rk} a_{lh}g^{nr}
\rg)
M^m{}_{ij}
+
y^n_ma_h{}^m{}_{ij} t^h_k.
$$

Let us lower here the index $n$
and use the equality
$
g_{nm}y^m_i=p^2t^j_n a_{ij}
$
(see the formulas below  (2.17)).
This yields
$$
\Rho_{knij}
=
\lf(
p^2t^l_n t^h_{km}
a_{lh}
-
(1-H)\fr1{F}(l_kg_{mn}+l_ng_{mk})
-p^2 t^l_mt^h_{nk} a_{lh}
\rg)
M^m{}_{ij}
+
p^2
a_{ml}
t^l_n
a_h{}^m{}_{ij} t^h_k.
$$
Next, we use here the skew-symmetry relation (2.20), obtaining
$$
\Rho_{knij}
=
\lf(
(1-H)\fr2{F}(l_ng_{mk}-l_mg_{kn})
-
(1-H)\fr1{F}(l_kg_{mn}+l_ng_{mk})
\rg)
M^m{}_{ij}
+
p^2
a_{ml}
t^l_n
a_h{}^m{}_{ij} t^h_k,
$$
\ses\\
or
\be
\Rho_{knij}
=
-
(1-H)\fr1{F}(l_kM_{nij}-l_nM_{kij})
+
p^2
a_{hlij} t^h_kt^l_n,
\ee
\ses\\
where
$
a_{hlij}=a_{lr}a_h{}^r{}_{ij}.
$
Finally, we return the index $n$ to the upper position, arriving at
\be
\Rho_k{}^n{}_{ij}
=
-
(1-H)\fr1{F}(l_k\de^n_m-l^ng_{mk})
M^m{}_{ij}
 +
y^n_ma_h{}^m{}_{ij} t^h_k.
\ee

{

The totally contravariant components
$$
\Rho^{knij}~:=g^{pk}a^{mi}a^{nj}
\Rho_p{}^n{}_{mn}
$$
read
\ses\\
\be
\Rho^{knij}
=
-
(1-H)\fr1{F}(l^kM^{nij}-l^nM^{kij})
 +
\fr1{p^2}
 y^k_h y^n_ra^{hrij},
\ee
where
$a^{hrij}=
a^{hl}a^{mi}a^{nj}
a_{l}{}^r{}_{mn}
$
and
$
M^{mij}~:=
a^{hi}a^{nj}
M^m{}_{hn}.
$

\ses

Similarly, we can conclude  from (5.13) that
the tensor
$ M_{nij}~:= g_{nm}M^m{}_{ij}$
can be given by means of the representation
\be
M_{nij}=-p^2t^ht^m_na_{hmij}.
\ee
Squaring yields
\ses
\be
 M^{nij} M_{nij}=
p^2
t^l
a_l{}^{nij}
t^h
a_{hnij}.
\ee

{

Now we square the $\rho$-tensor:
$$
\Rho^{knij}
\Rho_{knij}
=
(1-H)^2\fr2{F^2}M^{nij}M_{nij}
-
2(1-H)\fr1{F}(l^kM^{nij}-l^nM^{kij})
p^2
a_{hlij} t^h_kt^l_n
+
a^{knij}   a_{knij}
$$

\ses

\ses

\ses

$$
=
(1-H)^2\fr2{F^2}M^{nij}M_{nij}
-
2(1-H)H\fr1{F^2}
p^2
(a_{hlij} t^ht^l_nM^{nij}-a_{hlij} t^h_kt^lM^{kij})
+
a^{knij}   a_{knij},
$$
\ses
or
$$
\Rho^{knij}
\Rho_{knij}
=
(1-H)^2\fr{2p^2}
{F^2}
t^l
a_l{}^{nij}
t^h
a_{hnij}
+
2(1-H)\fr{Hp^2}{F^2}
(a_{hlij} t^ht^ra_r{}^{lij}-a_{hlij} t^lt^ra_r{}^{hij})
+
a^{knij}   a_{knij},
$$
\ses
which is
\be
\Rho^{knij}\Rho_{knij}=
a^{knij}a_{knij}
+\fr 2{S^2}
\lf(\fr1{H^2}-1\rg)
t^l
a_l{}^{nij}
t^h
a_{hnij}.
\ee

{

Because of the nullifications
$$
\cT_i \lf(\fr1 py^{n}_m\rg)=0, \qquad \cT_i (Hpt^m)=0
$$
(see  (4.3) and (4.10)),
from (5.13) it follows that
\be
\cT_l M^n{}_{ij}
=
-
y^n_tt^h
\lf(\nabla_l -\fr1HH_l\rg)    a_h{}^t{}_{ij}.
\ee

{

From (5.17) we can conclude that
$$
\cT_l \Rho_k{}^n{}_{ij}=
(1-H)\fr1{F}(l_k\de^n_m-l^ng_{mk})
y^m_tt^h
\lf(\nabla_l -\fr1HH_l\rg) a_h{}^t{}_{ij}
 +
y^n_mt^h_k \nabla_la_h{}^m{}_{ij}
$$

\ses

\ses

\be
+
H_l
\fr1{F}(l_k\de^n_m-l^ng_{mk})
M^m{}_{ij}.
\ee

{

The covariant derivatives
\ses\\
\be
\cT_k M^n{}_{ij}=d_k M^n{}_{ij}
+T^n{}_{kt} M^t{}_{ij}
 - a^s{}_{ki}   M^n{}_{sj}  -  a^s{}_{kj}   M^n{}_{is}
\ee
\ses\\
and
\be
\cT_l \Rho_k{}^n{}_{ij}
=
d_l\Rho_k{}^n{}_{ij}
+
T^n{}_{lt} \Rho_k{}^t{}_{ij}
-
T^t{}_{lk} \Rho_t{}^n{}_{ij}
 -  a^s{}_{li}
 \Rho_k{}^n{}_{sj}
 -   a^s{}_{lj}
 \Rho_k{}^n{}_{is}
\ee
\ses\\
have been used.

{

\ses

\ses

\setcounter{equation}{0}

\ses

\ses

\nin
{ \bf Appendix A:
Evaluations for Finsleroid  connection coefficients with ${\mathbf\it {g=g(x)}}$ }

\ses

\ses

Below we present various important   evaluations
which underlined the consideration
performed in Section I.3 of Chapter I.

We shall  use the relations
\ses\\
\be
\D{\fr bq}{y^n}=\fr{2B}{NKgq^2}A_n
\ee
\ses
and
$$
\D{\fr{q^2}B}{y^n}=-\fr{q^2}B \fr{q^2}B
\lf[
\fr{2B}{NKq^2}A_n +2\fr{2B}{NKgq^2}A_n\fr bq
\rg]
=- \fr{q^2}B
\fr{2}{NK}
\lf(
1 +\fr2g\fr bq
\rg)
A_n,
$$
\ses
so that
\be
\D{\fr{q^2}B}{y^n}=- \fr{q}B
\fr{2}{gNK}
(2b+gq)
A_n.
\ee
\ses\\
Moreover,
$$
\D{\fr{bq}B}{y^n}=
-
\fr{q}B
\fr{2}{gNK}
(2b+gq)
A_n
\fr bq
+
\fr{q^2}B
\fr{2B}{NKgq^2}A_n
$$
and
\ses\\
\be
\D{\fr{bq}B}{y^n}=
-
\fr{b}B
(2b+gq)
\fr{2}{gNK}
A_n
+
\fr{2}{gNK}
A_n.
\ee

We shall  also meet the convenience to apply the identity
\be
(2b+gq)
\lf(q+\fr12 gb\rg)
=
2h^2bq
+
gB.
\ee

{

The equality
\be
\D{\bar M} {y^n}
=
2
\fr{q^2}B
\fr{2}{gNK}
A_n
\ee
\ses\\
 can  be obtained from the relation
\be
\D{y_n}g=\bar M y_n
+
\fr12K^2
\D{\bar M} {y^n}.
\ee

It follows that
$$
\D{g_{mn}}g=\bar M g_{mn}
+
2
\fr{q^2}B
\fr{2}{gNK}
A_m y_n
+
y_m
\fr{q^2}B
\fr{2}{gNK}
A_n
$$

\ses

\ses

$$
- \fr{q}B
\fr{2}{gN}
(2b+gq)
A_m
\fr{2}{gN}
A_n
+
\fr{q^2}B
\fr{2}{gN}
\lf[
-A_ml_n
+
\fr2NA_mA_n
-
\fr{gN}2\fr b{q}{\cal H}_{mn}
\rg],
$$
\ses\\
which is
\be
\D{g_{mn}}g=\bar M g_{mn}
+
\fr{q^2}B
\fr{2}{gN}
(A_m l_n+A_n l_m)
-
 \fr{bq}B
\fr{2}{gN}
\fr{2}{gN}
A_mA_n
-
\fr{bq}B
h_{mn},
\ee
\ses\\
entailing
\be
\D{h_{mn}}g=
-
 \fr{bq}B
\fr{2}{gN}
\fr{2}{gN}
A_mA_n
-
\fr{bq}B
h_{mn}.
\ee

We can also obtain

\ses

\ses

$$
\D{A_{mnj}}g=
\fr32
\bar M
A_{mnj}
+
\lf(\fr1g-\fr{bq}B\rg)
A_{mnj}
-\fr{gbq}B
\fr{1}{gN}
(A_m h_{nj}+A_n h_{mj}+  A_j h_{mn})
$$

\ses

\ses

$$
-
\fr{gbq}B
\fr{1}{gN}
\fr{2}{gN}
\fr{2}{gN}
A_jA_mA_n,
$$
\ses\\
or
\be
\D{A_{mnj}}g=
\fr32
\bar M
A_{mnj}
+
\lf(\fr1g-2\fr{bq}B\rg)
A_{mnj}
-
2
\fr{gbq}B
\fr{1}{gN}
\fr{2}{gN}
\fr{2}{gN}
A_jA_mA_n.
\ee

{

Evaluations frequently involve the vector
$m_i=(2/Ng)A_i$ which possesses the properties
$$
g^{ij}m_im_j=1, \qquad y^im_i=0.
$$
From  (A.24) of [7]
it follows that
\be
m_i=K\fr1{q}(b_i-\fr b{K^2}y_i).
\ee
The equality
\be
K\D{m_i}{y^n}
=-m_nl_i
+
gm_nm_i
-\fr b{q}{\cal H}_{in}
\ee
holds, where  ${\cal H}_{in}=h_{in}-m_im_n$.

\ses

The contravariant components $m^i$ can be taken from
(A.27)  of [7]:
\be
m^i=\fr 1{qK}
\Bigl[q^2b^i-(b+gq)v^i\Bigr],
\ee
entailing
\be
K\D{m^i}{y^n}=-m_nl^i
-gm^im_n
-
\fr 1{q}
(b+gq){\cal H}^i_n.
\ee

\ses

With the representation
\be
A_{ijk}= \fr1N
 \lf[
A_ih_{jk}  +A_jh_{ik}  +A_kh_{ij}
-\fr 4{N^2g^2}
A_iA_jA_k
\rg]
\ee
\ses\\
(see (A.8)  in [7]),
we find that
\be
\D{A_mA^k}   {y^n}
=
\fr1KA^k
\lf[
 -A_nl_m
-\fr{Ng}2
\fr b{q}{\cal H}_{mn}
 \rg]
+
\fr1KA_m
\lf[
-A_nl^k
-
\fr{Ng}2\fr 1{q}
(b+gq){\cal H}^k_n
  \rg].
\ee

{

\ses

Recollecting the scalar
$h(x)=\sqrt{1- (g^2(x)/4}$
and introducing the scalar
 $G=g/h$,
\ses
we get
 \be
  \D hg= -\fr14 G, \qquad  \D Gg= \fr1{h^3},
  \ee
\ses\\
so that
\be
 \D {K^2}g=\bar M
K^2,
\ee
where
$$
\bar M=
 \fr{b q}B - \fr1{h^3}f+
\fr12\fr{G}{hB} ( q^2+\frac12 gb q),
$$
\ses\\
or
\be
\bar M=
 - \fr1{h^3}f+
\fr12\fr{G}{hB}  q^2+  \frac1{h^2B} bq.
\ee

\ses

  The  function $K(x,y)$ is given by the formulas
\be
K(x,y)=
\sqrt{B(x,y)}\,J(x,y), \qquad
J(x,y)=\e^{-\frac12G(x)f(x,y)},
\ee
\ses\\
entailing
$$
\D{\ln J}{y^n} =  \fr1N C_n
$$
\ses\\
and
$$
\D {\bar M}{y^n}
=
 \fr1{h^2}
\fr{2}{gNK}
A_n
-
\fr12\fr{g}{h^2}
 \fr{q}B
\fr{2}{gNK}
(2b+gq)
A_n
+
  \frac1{h^2}
\lf(
-
\fr{b}B
(2b+gq)
\fr{2}{gNK}
A_n
+
\fr{2}{gNK}
A_n
\rg),
$$
\ses\\
or
$$
\D{\bar  M}{y^n}
=
\lf(
-
g
q\lf(b+\fr12gq\rg)
-
b(2b+gq)
+
2B
\rg)
\fr1{h^2}
\fr1B
\fr{2}{gNK}
A_n
$$
\ses
which is equivalent to  (A.5).

{

Starting with (I.3.22)-(I.3.23), we get
\ses\\
$$
\breve N^k{}_{im}
=
\fr1{h^2}g_i
\fr{q^2}{2B}
\lf(1+\fr12 g\fr bq\rg)
\fr2{Ng}
A_m
 l^k
+
\fr1{h^2}g_i
\fr{q^2}{2B}
\lf(1+\fr12 g\fr bq\rg)
\fr 1{q}(b+gq){\cal H}^k_m
$$

\ses

\ses

$$
-
\fr1{h^2}g_i
\lf[
\fr1N
\fr1{Ng}
A_m  A^k
+
l_m
\fr{q^2}B\lf(1+\fr12 g\fr bq\rg)
\fr1{Ng}
A^k
-
\fr2{q}  (b+gq)
\fr2{Ng}
A_m
\fr{q^2}B\lf(1+\fr12 g\fr bq\rg)
\fr1{Ng}
A^k
\rg]
$$

\ses

\ses

$$
-
g_i
\lf(
\fr{q^2}B
\fr{2  }{Ng} A_m
l^k
+
\fr12
\bar M\de^k_m\rg)
$$

\ses

\ses

\ses

\ses

$$
=
\fr1{h^2}g_i
\fr{q^2}{2B}\lf(1+\fr12 g\fr bq-2h^2 \rg)
\fr2{Ng}
A_m
 l^k
+
\fr1{h^2}g_i
\fr{1}{2B}
\lf(q+\fr12 gb\rg)
(b+gq)h^k_m
$$

\ses

\ses

$$
-
\fr1{h^2}g_i
\fr1B
\lf[
gB
-
  (2b+2gq)
\lf(q+ \fr12gb\rg)
\rg]
\fr1{Ng}
\fr1{Ng}
A_m  A^k
-
\fr12
g_i
\bar Mh^k_m
+
\fr1K
l_m\breve N^k{}_{i}.
$$
\ses\\
Using here the equality
(A.4)
leads to
\ses\\
$$
\breve N^k{}_{im}
=
\fr1{h^2}g_i
\fr{q^2}{2B}\lf(1+\fr12 g\fr bq-2h^2 \rg)
\fr2{Ng}
A_m
 l^k
+
\fr1{h^2}g_i
\fr{1}{2B}
\lf(q+\fr12 gb\rg)
(b+gq)h^k_m
$$

\ses

\ses

$$
+
\fr1{h^2}g_i
\fr1B
\lf[
2h^2bq
+
gq
\lf(q+ \fr12gb\rg)
\rg]
\fr1{Ng}
\fr1{Ng}
A_m  A^k
-
\fr12
g_i
\bar Mh^k_m
+
\fr1K
l_m\breve N^k{}_{i}.
$$
\ses\\
Eventually we obtain
\ses\\
$$
\breve N^k{}_{im}
=
\fr1{h^2}g_i
\fr{q^2}{2B}\lf(1+\fr12 g\fr bq-2h^2 \rg)
\fr2{Ng}
A_m
 l^k
+
\fr1{h^2}g_i
\fr{q^2}{2B}
\lf(1+\fr12 g\fr bq\rg)
\lf(\fr bq+g\rg)
h^k_m
$$

\ses

\ses

\be
+
\fr1{h^2}g_i
\fr{q^2}{2B}
\lf(
\fr bq+\fr12g
\rg)
\fr2{Ng}
\fr2{Ng}
A_m  A^k
-
\fr12
g_i
\bar Mh^k_m
+
\fr1K
l_m\breve N^k{}_{i}.
\ee
\ses\\
Thus the representation (I.3.28) is valid.

{

Next, we find that
\ses\\
$$
\breve N^k{}_{imn}=\D{\breve N^k{}_{im}}{y^n}
=
-
\fr1{h^2}g_i
\fr{q}B
\fr{1}{gNK}
(2b+gq)
A_n
\lf(1+\fr12 g\fr bq-2h^2 \rg)
\fr2{Ng}
A_m
 l^k
$$

\ses

\ses

$$
 +
\fr1{h^2}g_i
\fr12 g
 \fr{1}{NKg}A_n
\fr2{Ng}
A_m
 l^k
$$

\ses

\ses

$$
+
\fr1{h^2}g_i
\fr{q^2}{2B}\lf(1+\fr12 g\fr bq-2h^2 \rg)
\fr2{Ng}
\fr1K
\lf[
-
A_nl_m
+
\fr2N
A_nA_m
-
\fr b{q}
\fr{Ng}2
{\cal H}_{mn}
\rg]
 l^k
$$

\ses

\ses

$$
+
\fr1{h^2}g_i
\fr{q^2}{2B}\lf(1+\fr12 g\fr bq-2h^2 \rg)
\fr2{Ng}
A_m
\fr1K h^k_n
 $$

\ses

\ses

\ses

\ses

\ses

\ses

$$
-
\fr1{h^2}g_i
\fr{q}B
\fr{1}{gNK}
(2b+gq)
A_n
\lf(1+\fr12 g\fr bq\rg)
\lf(\fr bq+g\rg)
h^k_m
$$

\ses

\ses

$$
+
\fr1{h^2}g_i
\fr{q^2}{2B}
\fr12 g
\fr{2B}{NKgq^2}A_n
\lf(\fr bq+g\rg)
h^k_m
+
\fr1{h^2}g_i
\fr{q^2}{2B}
\lf(1+\fr12 g\fr bq\rg)
\fr{2B}{NKgq^2}A_n
h^k_m
$$

\ses

\ses

$$
-
\fr1{h^2}g_i
\fr{q^2}{2B}
\lf(1+\fr12 g\fr bq\rg)
\lf(\fr bq+g\rg)
        \fr1K(l^kh_{mn}+l_mh^k_n)
$$

\ses

\ses

{

\ses

\ses

$$
-
\fr1{h^2}g_i
\fr{q}B
\fr{1}{gNK}
(2b+gq)
A_n
\lf(
\fr bq+\fr12g
\rg)
\fr4{N^2g^2}
A_m  A^k
+
\fr1{h^2}g_i
\fr{q^2}{2B}
\fr{2B}{NKgq^2}A_n
\fr2{Ng}
\fr2{Ng}
A_m  A^k
$$

\ses

\ses

$$
-
\fr1{h^2}g_i
\fr{q^2}{2B}
\lf(
\fr bq+\fr12g
\rg)
\fr4{N^2g^2}
\fr1K
\lf[
A^k
\lf(
A_nl_m
+\fr{Ng}2
\fr b{q}{\cal H}_{mn}
 \rg)
+
A_m
\lf(
A_nl^k
+
\fr{Ng}{2q}
(b+gq){\cal H}^k_n
  \rg)
\rg]
$$

\ses

\ses

\ses

\ses

$$
-
g_i
\fr{q^2}B
\fr{2}{gNK}
A_n
h^k_m
+
\fr12
g_i
\bar M
\fr1K(l^kh_{mn}+l_mh^k_n)
$$

\ses

\ses

\ses

$$
+
\fr1{K^2}
h_{mn}
\lf[
-
\fr1{h^2}g_i
\fr{q}B\lf(q+\fr12 gb\rg)
\fr K{Ng}
A^k
-
\fr12
g_i
\bar My^k
\rg]
$$

\ses

\ses

$$
+
\fr1K
l_m
\fr1{h^2}g_i
\fr{q^2}{2B}\lf(1+\fr12 g\fr bq-2h^2 \rg)
\fr2{Ng}
A_n
 l^k
$$

\ses

\ses

$$
        +
\fr1K
l_m
\lf[
\fr1{h^2}g_i
\fr{q^2}{2B}
\lf(1+\fr12 g\fr bq\rg)
\lf(\fr bq+g\rg)
h^k_n
+
\fr1{h^2}g_i
\fr{q^2}{2B}
\lf(
\fr bq+\fr12g
\rg)
\fr4{N^2g^2}
A_n  A^k
-
\fr12
g_i
\bar Mh^k_n
\rg].
$$

\ses

\ses

{

Simplifying yields
\ses\\
$$
\breve N^k{}_{imn}
=
-
\fr1{h^2}g_i
\fr{bq}B
A_n
\lf(1+\fr12 g\fr bq-2h^2 \rg)
\fr4{N^2g^2}
\fr1K
A_m
 l^k
 +
\fr1{h^2}g_i
\fr12 g
 \fr{1}{NKg}A_n
\fr2{Ng}
A_m
 l^k
$$

\ses

\ses

$$
-
\fr1{h^2}g_i
\fr{bq}{2B}\lf(1+\fr12 g\fr bq-2h^2 \rg)
\fr1K
h_{mn}
 l^k
+
\fr1{h^2}g_i
\fr{bq}{2B}\lf(1+\fr12 g\fr bq-2h^2 \rg)
\fr1K
\fr4{N^2g^2}
A_mA_n
 l^k
$$

\ses

\ses

$$
+
\fr1{h^2}g_i
\fr{q^2}{2B}\lf(1+\fr12 g\fr bq-2h^2 \rg)
\fr2{Ng}
A_m
\fr1K h^k_n
 $$

\ses

\ses

$$
-
\fr1{h^2}g_i
\fr{1}B
(2b+gq)
\lf(q+\fr12 gb\rg)
\lf(\fr bq+g\rg)
\fr{1}{gNK}
A_n
h^k_m
$$

\ses

\ses

$$
+
\fr1{h^2}g_i
\fr12 g
\fr{1}{NKg}A_n
\lf(\fr bq+g\rg)
h^k_m
+
\fr1{h^2}g_i
\lf(1+\fr12 g\fr bq\rg)
\fr{1}{NKg}A_n
h^k_m
$$

\ses

\ses

$$
-
\fr1{h^2}g_i
\fr{q^2}{2B}
\lf(1+\fr12 g\fr bq\rg)
\lf(\fr bq+g\rg)
        \fr1Kl^kh_{mn}
$$

\ses

\ses

\ses

\ses

$$
-
\fr1{h^2}g_i
\fr{q}{2B}
(2b+gq)
\lf(
\fr bq+\fr12g
\rg)
\fr2{Ng}
\fr2{Ng}
\fr2{Ng}
\fr1K
A_nA_m  A^k
$$

\ses

\ses

$$
+
\fr1{h^2}g_i
\fr{1}{NKg}A_n
\fr2{Ng}
\fr2{Ng}
A_m  A^k
$$

\ses

\ses

$$
-
\fr1{h^2}g_i
\fr{q^2}{2B}
\lf(
\fr bq+\fr12g
\rg)
\fr2{Ng}
\fr2{Ng}
\fr1K
\lf[
A^k
\fr{Ng}2
\fr b{q}{\cal H}_{mn}
+
A_m
A_nl^k
\rg]
$$

\ses

\ses

$$
-
\fr1{h^2}g_i
\fr{q}{2B}
\lf(
\fr bq+\fr12g
\rg)
\fr2{Ng}
\fr1K
(b+gq)A_m{\cal H}^k_n
$$

\ses

\ses

\ses

\ses

$$
-
g_i
\fr{q^2}B
\fr{2}{gNK}
A_n
h^k_m
-
\fr1{K^2}
h_{mn}
\fr1{h^2}g_i
\fr{q}B\lf(q+\fr12 gb\rg)
\fr K{Ng}
A^k,
$$

\ses

\ses

{

\nin
or
\ses\\
$$
\breve N^k{}_{imn}
=
\fr1{h^2}g_i
\fr12 g
 \fr{1}{NKg}A_n
\fr2{Ng}
A_m
 l^k
$$

\ses

\ses

$$
-
\fr1{h^2}g_i
\fr{bq}{2B}\lf(2+ g\fr bq-2h^2 \rg)
\fr1K
h_{mn}
 l^k
-
\fr1{h^2}g_i
\fr{bq}{2B}\lf(1+\fr12 g\fr bq-2h^2 \rg)
\fr1K
\fr2{Ng}
\fr2{Ng}
A_mA_n
 l^k
$$

\ses

\ses

$$
+
\fr1{h^2}g_i
\fr{q^2}{2B}\lf(1+\fr g2 \fr bq-2h^2 \rg)
\fr2{Ng}
A_m
\fr1K h^k_n
-
\fr1{h^2}g_i
\fr{1}B
(2h^2bq+gB)
\lf(\fr bq+g\rg)
\fr{1}{gNK}
A_n
h^k_m
$$

\ses

\ses

$$
+
\fr1{h^2}g_i
\fr12 g
\fr{1}{NKg}A_n
\lf(\fr bq+g\rg)
h^k_m
+
\fr1{h^2}g_i
\lf(1+\fr12 g\fr bq\rg)
\fr{1}{NKg}A_n
h^k_m
$$

\ses

\ses

$$
-
\fr1{h^2}g_i
\fr{q^2}{2B}
\lf(1+\fr12 g\fr bq\rg)
g
        \fr1Kl^kh_{mn}
-
\fr1{h^2}g_i
\fr{bq}{2B}
\lf(
\fr bq+\fr12g
\rg)
\fr2{Ng}
\fr2{Ng}
\fr2{Ng}
\fr1K
A_nA_m  A^k
$$

\ses

\ses

$$
+
\fr1{h^2}g_i
\fr{1}{NKg}A_n
\fr2{Ng}
\fr2{Ng}
A_m  A^k
$$

\ses

\ses

$$
-
\fr1{h^2}g_i
\fr{q^2}{2B}
\lf(
\fr bq+\fr12g
\rg)
\fr2{Ng}
\fr1K
A^k
\fr b{q}{\cal H}_{mn}
-
\fr1{h^2}g_i
\fr{q^2}{2B}
\lf(
\fr bq+\fr12g
\rg)
\fr2{Ng}
\fr2{Ng}
\fr1K
A_m
A_nl^k
$$

\ses

\ses

$$
-
\fr1{h^2}g_i
\fr{q}{2B}
\lf(
\fr bq+\fr12g
\rg)
\fr2{Ng}
\fr1K
(b+gq)A_mh^k_n
$$

\ses

\ses

\ses

\ses

$$
-
g_i
\fr{q^2}B
\fr{2}{gNK}
A_n
h^k_m
-
\fr1{K^2}
h_{mn}
\fr1{h^2}g_i
\fr{q}B\lf(q+\fr12 gb\rg)
\fr K{Ng}
A^k.
$$

\ses

\ses

{

Additional reductions are possible,
leading to
\ses\\
$$
\breve N^k{}_{imn}
=
\fr1{h^2}g_i
\fr12 g
 \fr{1}{NKg}A_n
\fr2{Ng}
A_m
 l^k
$$

\ses

\ses

$$
-
\fr g{2h^2}g_i
\fr1K
h_{mn}
 l^k
-
\fr1{h^2}g_i
\fr{bq}{2B}\lf(2+\fr12 g\fr bq-2h^2 \rg)
\fr1K
\fr2{Ng}
\fr2{Ng}
A_mA_n
 l^k
$$

\ses

\ses

$$
+
\fr1{h^2}g_i
\fr{q^2}{2B}\lf(1+\fr12 g\fr bq-2h^2 \rg)
\fr2{Ng}
A_m
\fr1K h^k_n
 $$

\ses

\ses

$$
-
\fr1{h^2}g_i
\lf(2h^2+g\lf(\fr bq+g\rg)\rg)
\fr{1}{gNK}
A_n
h^k_m
$$

\ses

\ses

$$
+
\fr1{h^2}g_i
\fr12 g
\fr{1}{NKg}A_n
\lf(\fr bq+g\rg)
h^k_m
+
\fr1{h^2}g_i
\lf(1+\fr12 g\fr bq\rg)
\fr{1}{NKg}A_n
h^k_m
$$

\ses

\ses

\ses

\ses

$$
+
\fr1{h^2}g_i
\fr{1}{NKg}A_n
\fr2{Ng}
\fr2{Ng}
A_m  A^k
$$

\ses

\ses

$$
-
\fr1{h^2}g_i
\fr{bq}{2B}
\lf(
\fr bq+\fr12g
\rg)
\fr2{Ng}
\fr1K
A^k
h_{mn}
-
\fr1{h^2}g_i
\fr{q^2}{2B}
\fr12g
\fr2{Ng}
\fr2{Ng}
\fr1K
A_m
A_nl^k
$$

\ses

\ses

$$
-
\fr1{h^2}g_i
\fr{q}{2B}
\lf(
\fr bq+\fr12g
\rg)
\fr2{Ng}
\fr1K
(b+gq)A_mh^k_n
-
\fr1{K^2}
h_{mn}
\fr1{h^2}g_i
\fr{q}B\lf(q+\fr12 gb\rg)
\fr K{Ng}
A^k,
$$

{

\nin
or simply
$$
\breve N^k{}_{imn}
=
-
\fr g{2h^2}g_i
\fr1K
h_{mn}
 l^k
+
\fr1{h^2}g_i
\fr{1}{2B}\lf(-q^2+\fr12 g bq+\fr12g^2q^2 \rg)
\fr2{Ng}
\fr1K h^k_n
A_m
 $$

\ses

\ses

$$
-
\fr1{h^2}g_i
\lf(2+\fr{g^2}2+g\fr bq\rg)
\fr{1}{gNK}
A_n
h^k_m
+
\fr1{h^2}g_i
\fr{g}{2NK}A_n
h^k_m
+
\fr1{h^2}g_i
\lf(1+ g\fr bq\rg)
\fr{1}{NKg}A_n
h^k_m
$$

\ses

\ses

$$
+
\fr1{h^2}g_i
\fr{1}{NKg}A_n
\fr2{Ng}
\fr2{Ng}
A_m  A^k
$$

\ses

\ses

\be
-
\fr1{h^2}g_i
\fr{1}{2B}
\lf(
b+\fr12gq
\rg)
(b+gq)
\fr2{Ng}
\fr1K
A_mh^k_n
-
\fr1{h^2}g_i
\fr 1{Ng}
\fr1K
h_{mn}
A^k.
\ee
\ses\\
Using here the representation (A.14) of the tensor $A_{ijk}$,
we are
coming to
\be
y_k\Dd{\breve N^k{}_i}{y^m}{y^n}=
\fr2h
h_i
h_{mn}.
\ee

{

Let us verify the validity of the equality
\be
\cD_i h_{nm}
=
-\fr2{h}h_i
h_{nm}.
\ee
\ses\\
To this end
 we find
\ses\\
$$
2
\breve N^k{}_i
C_{kmn}
=
-
\fr1{h^2}g_i
\fr{q}B\lf(q+\fr12 gb\rg)
\fr 2{Ng}
A^k
A_{kmn}
$$

\ses

\ses

$$
=
-
\fr1{h^2}g_i
\fr{q}B\lf(q+\fr12 gb\rg)
\fr 2{Ng}
A^k
\fr1N
 \lf[
A_nh_{nk}  +A_nh_{mk}  +A_kh_{mn}
-\fr 4{N^2g^2}
A_mA_nA_k
\rg],
$$
\ses
so that
\be
2
\breve N^k{}_i
C_{kmn}
=
-
\fr1{h^2}g_i
\fr{q}B\lf(q+\fr12 gb\rg)
\fr 2{Ng}
\fr1N
 \lf[
A_mA_n  + \fr {N^2g^2}4h_{mn}
\rg].
\ee

{

We can also observe that
\ses\\
$$
g_i\D{g_{mn}}g
+
2
\breve N^k{}_i
C_{kmn}
=
g_i
\bar M g_{mn}
+
g_i
\fr{q^2}B
\fr{2}{gN}
(A_m l_n+A_n l_m)
-
g_i
 \fr{bq}B
\fr4{N^2g^2}
A_mA_n
-
g_i
\fr{bq}B
h_{mn}
$$
\ses

\ses

$$
 -
\fr1{h^2}g_i
\fr{q}B\lf(q+\fr12 gb\rg)
 \lf[
\fr 2{Ng}
\fr1N
A_mA_n  + \fr  g2h_{mn}
\rg].
$$
\ses\\
Simultaneously,
$$
g_{kn}\breve N^k{}_{im}
=
\fr1{h^2}g_i
\fr{q^2}{2B}\lf(1+\fr12 g\fr bq-2h^2 \rg)
\fr2{Ng}
A_m
 l_n
+
\fr1{h^2}g_i
\fr{q^2}{2B}
\lf(1+\fr12 g\fr bq\rg)
\lf(\fr bq+g\rg)
h_{nm}
$$

\ses

\ses

$$
+
\fr1{h^2}g_i
\fr{q^2}{2B}
\lf(
\fr bq+\fr12g
\rg)
\fr2{Ng}
\fr2{Ng}
A_m  A_n
-
\fr12
g_i
\bar Mh_{mn}
$$

\ses

\ses

$$
+
\fr1K
l_m
\lf[
-
\fr1{h^2}g_i
\fr{q}B\lf(q+\fr12 gb\rg)
\fr K{Ng}
A_n
-
\fr12
g_i
\bar My_n
\rg].
$$

{

In this way we obtain
$$
g_i\D{g_{mn}}g
+
2
\breve N^k{}_i
C_{kmn}
+
g_{km}\breve N^k{}_{in}
+
g_{kn}\breve N^k{}_{im}
$$

\ses

\ses

  \ses

$$
=
g_i
\fr{q^2}B
\fr{2}{gN}
(A_m l_n+A_n l_m)
-
g_i
 \fr{bq}B
\fr{2}{gN}
\fr{2}{gN}
A_mA_n
-
g_i
\fr{bq}B
h_{mn}
$$
\ses

\ses

$$
 -
\fr1{h^2}g_i
\fr{q}B\lf(q+\fr12 gb\rg)
 \lf[
\fr 2{Ng}
\fr1N
A_mA_n  + \fr  g2h_{mn}
\rg]
$$

  \ses

  \ses

  \ses

$$
+
  \fr1{h^2}g_i
\fr{q^2}{2B}\lf(1+\fr12 g\fr bq-2h^2 \rg)
\fr2{Ng}
(A_m l_n+A_n l_m)
+
\fr1{h^2}g_i
\fr{q^2}{B}
\lf(1+\fr12 g\fr bq\rg)
\lf(\fr bq+g\rg)
h_{nm}
$$

\ses

\ses

$$
+
\fr1{h^2}g_i
\fr{q^2}{B}
\lf(
\fr bq+\fr12g
\rg)
\fr2{Ng}
\fr2{Ng}
A_m  A_n
-
\fr1{h^2}g_i
\fr{q}B\lf(q+\fr12 gb\rg)
\fr 1{Ng}
(l_mA_n+l_nA_m).
$$
\ses\\
Reducing similar terms leads to
$$
g_i\D{g_{mn}}g
+
2
\breve N^k{}_i
C_{kmn}
+
g_{km}\breve N^k{}_{in}
+
g_{kn}\breve N^k{}_{im}
$$

\ses

\ses

  \ses

$$
=
-
g_i
 \fr{bq}B
\fr{2}{gN}
\fr{2}{gN}
A_mA_n
-
g_i
\fr{bq}B
h_{mn}
 -
\fr1{h^2}g_i
\fr{q}B\lf(q+\fr12 gb\rg)
 \lf[
\fr 2{Ng}
\fr1N
A_mA_n  + \fr  g2h_{mn}
\rg]
$$
\ses

\ses

$$
+
\fr1{h^2}g_i
\fr{q^2}{B}
\lf(1+\fr12 g\fr bq\rg)
\lf(\fr bq+g\rg)
h_{nm}
+
\fr1{h^2}g_i
\fr{q^2}{B}
\lf(
\fr bq+\fr12g
\rg)
\fr2{Ng}
\fr2{Ng}
A_m  A_n
$$

\ses

\ses

\ses

$$
=
-
g_i
\fr{bq}B
h_{mn}
+
\fr1{h^2}g_i
\fr1{2B}
(2h^2bq+gB)
h_{nm}
=
\fr1{2h^2}gg_i
h_{nm}.
$$
\ses\\
We get
\be
\cD_i g_{nm}
=
\fr1{2h^2}gg_i
h_{nm}.
\ee
Thus the equality (A.23) is valid.

{

Now we want
to verify the validity of the equality (A.9).
Differentiating  (A.14) with respect to  $y^j$ yields

\ses

$$
2\D{C_{mnj}}g=
2\bar M C_{mnj}
+
2
\fr{q^2}B
\fr{2}{gNK}
A_j
 g_{mn}
$$

\ses

\ses

$$
-
\fr{q}B
\fr{2}{gNK}
(2b+gq)
A_j
\fr{2}{gN}
(A_m l_n+A_n l_m)
+
\fr1K
\fr{q^2}B
\fr{2}{gN}
(A_m h_{nj}+A_n h_{mj})
$$

\ses

\ses

$$
+
\fr{q^2}B
\fr{2}{gN}
\fr1K
\lf[
\lf(
-A_jl_m
+
\fr2NA_mA_j
-\fr b{q}
\fr{Ng}2
{\cal H}_{mj} \rg)
 l_n
 +
\lf(
-A_jl_n
+
\fr2NA_nA_j
-\fr b{q}
\fr{Ng}2
{\cal H}_{nj} \rg)
   l_m
\rg]
$$

\ses

\ses

$$
+
\lf(\fr{b}B
(2b+gq)
\fr{2}{gNK}
A_j
-
\fr{2}{gNK}
A_j
\rg)
\fr{2}{gN}
\fr{2}{gN}
A_mA_n
$$

\ses

\ses

$$
-
 \fr{bq}B
\fr{4}{g^2N^2}
\fr1K
\lf[
\lf(
-A_jl_m
+
\fr2NA_mA_j
-\fr b{q}
\fr{Ng}2
{\cal H}_{mj} \rg)
 A_n
 +
\lf(
-A_jl_n
+
\fr2NA_nA_j
-\fr b{q}
\fr{Ng}2
{\cal H}_{nj} \rg)
A_m
\rg]
$$

\ses

\ses

\ses

$$
+
\lf(
\fr{b}B
(2b+gq)
\fr{2}{gNK}
A_j
-
\fr{2}{gNK}
A_j
\rg)
h_{mn}
+
\fr{bq}B
\fr1K(l_mh_{jn}+l_nh_{mj})
-
\fr{2bq}BC_{mnj}.
$$

{

We may reduce as follows:
\ses\\
$$
2\D{C_{mnj}}g=
2\bar M C_{mnj}
+
2
\fr{q^2}B
\fr{2}{gNK}
A_j
h_{mn}
+
\fr1K
\fr{q^2}B
\fr{2}{gN}
(A_m h_{nj}+A_n h_{mj})
$$

\ses

\ses

\ses

$$
+
\lf(
\fr{B-q^2+b^2}B
\fr{2}{gNK}
-
\fr{2}{gNK}
\rg)
\fr{2}{gN}
\fr{2}{gN}
A_jA_mA_n
$$

\ses

\ses

$$
-
 \fr{bq}B
\fr{2}{gN}
\fr{2}{gN}
\fr1K
\Bigl[
\lf(
\fr2NA_mA_j
-
\fr b{q}
\fr{Ng}2
h_{mj}
+
\fr b{q}
\fr2{Ng}
A_mA_j
\rg)
 A_n
$$

\ses

\ses

$$
 +
\lf(
\fr2NA_nA_j
-
\fr b{q}
\fr{Ng}2
h_{nj}
+
\fr b{q}
\fr2{Ng}
A_nA_j
 \rg)
A_m
\Bigr]
$$

\ses

\ses

\ses

$$
+
\lf(
\fr{B-q^2+b^2}B
\fr{2}{gNK}
-
\fr{2}{gNK}
\rg)
A_j
h_{mn}
-
\fr{2bq}BC_{mnj}
$$

\ses

\ses

\ses

         \ses

$$
=
2\bar M C_{mnj}
+
\fr1K
\fr{B-gbq}B
\fr{2}{gN}
(A_m h_{nj}+A_n h_{mj}+  A_j h_{mn})
-
\fr{2bq}BC_{mnj}
 $$

\ses

\ses

\ses

$$
+
\fr{b^2-q^2}B
\fr{2}{gNK}
\fr{2}{gN}
\fr{2}{gN}
A_jA_mA_n
-
 \fr{2bq}B
\fr{2}{gN}
\fr{2}{gN}
\fr1K
\lf(
\fr2NA_mA_j
+
\fr b{q}
\fr2{Ng}
A_mA_j
\rg)
 A_n
$$

\ses

\ses

\ses

         \ses

$$
=
2\bar M
\fr1K
A_{mnj}
+
\fr1K
\fr{B-gbq}B
\fr{2}{gN}
(A_m h_{nj}+A_n h_{mj}+  A_j h_{mn})
$$

\ses

\ses

$$
-
\fr{B+gbq}B
\fr{2}{gNK}
\fr{2}{gN}
\fr{2}{gN}
A_jA_mA_n
-
\fr{2bq}BC_{mnj}
$$

\ses

\ses

\ses

         \ses

$$
=
2\bar M
\fr1K
A_{mnj}
+
\fr2K
\lf(\fr1g-\fr{bq}B\rg)
A_{mnj}
-
\fr1K
\fr{gbq}B
\fr{2}{gN}
(A_m h_{nj}+A_n h_{mj}+  A_j h_{mn})
$$

\ses

\ses

$$
-
\fr{gbq}B
\fr{2}{gNK}
\fr{2}{gN}
\fr{2}{gN}
A_jA_mA_n.
$$
Thus (A.9) is valid.

{

Next, we evaluate the term
\ses\\
$$
\breve N^k{}_i
\D{A_{jmn}}{y^k}
=
-
\fr1{h^2}g_i
\fr{q}B\lf(q+\fr12 gb\rg)
\fr K{Ng}
A^k
\D{A_{jmn}}{y^k}.
$$
\ses
With the representation
\ses\\
$$
\D{A_{ijk}}{y^n}
= \fr1K    \fr2N\lf(A_{jkn}A_i+A_{ikn}A_j+A_{ijn}A_k  \rg)
- \fr1K \lf(l_jA_{kni}+l_iA_{knj}+l_kA_{ijn}\rg)
$$

\ses

\ses

\be
+  \fr1K \fr1N \fr2N \lf({\cal H}_{jk} A_iA_n  +{\cal H}_{ik} A_jA_n+{\cal H}_{ij} A_kA_n\rg)
 -  \fr{gb}{2K q}
\lf( {\cal H}_{jk}{\cal H}_{in} + {\cal H}_{ik}{\cal H}_{jn}+ {\cal H}_{ji}{\cal H}_{kn}\rg)
\ee
\ses\\
we obtain
\ses\\
$$
KA^n
\D{A_{ijk}}{y^n}
=     \fr2N\lf(A_{jkn}A_i+A_{ikn}A_j+A_{ijn}A_k  \rg)
A^n
-  \lf(l_jA_{kni}+l_iA_{knj}+l_kA_{ijn}\rg)
A^n
$$

\ses

$$
+  \fr {g^2}2
 \lf({\cal H}_{jk} A_i  +{\cal H}_{ik} A_j+{\cal H}_{ij} A_k\rg).
$$

Using here the equality
\ses\\
$$
A_{jkn}A^n=\fr1N\lf(\fr{Ng}2\fr{Ng}2h_{jk}+A_jA_k\rg)
$$
leads to
\ses\\
$$
KA^n
\D{A_{ijk}}{y^n}
=
    \fr2N\fr1N
\lf(\fr{Ng}2\fr{Ng}2
(h_{jk}A_i+h_{ik}A_j+h_{ij}A_k)
+
3A_iA_jA_k
 \rg)
$$

\ses

\ses

$$
-
\fr1N
\fr{Ng}2\fr{Ng}2
(h_{jk}l_i+h_{ik}l_j+h_{ij}l_k)
-
\fr1N
(A_jA_kl_i+A_iA_kl_j+A_iA_jl_k)
$$

\ses

\ses

$$
+  \fr {g^2}2
 \lf(h_{jk} A_i  +h_{ik} A_j+h_{ij} A_k\rg)
-  3
\fr {g^2}2
\fr2{Ng}\fr2{Ng} A_iA_jA_k.
 $$

So we may write
\ses\\
$$
KA^n
\D{A_{ijk}}{y^n}
=
-
\fr1N
\fr{Ng}2\fr{Ng}2
(h_{jk}l_i+h_{ik}l_j+h_{ij}l_k)
-
\fr1N
(A_jA_kl_i+A_iA_kl_j+A_iA_jl_k)
$$

\ses

\ses

\be
+  g^2
 \lf(h_{jk} A_i  +h_{ik} A_j+h_{ij} A_k\rg).
\ee

{

We also need the term
\ses\\
$$
A_{jkn}
\breve N^k{}_{im}
=
\fr1{h^2}g_i
\fr{q^2}{2B}
\lf(1+\fr12 g\fr bq\rg)
\lf(\fr bq+g\rg)
A_{jmn}
$$

\ses

\ses

$$
+
\fr1{h^2}g_i
\fr{q^2}{2B}
\lf(
\fr bq+\fr12g
\rg)
\fr2{Ng}
\fr2{Ng}
A_m
\fr1N\lf(\fr{Ng}2\fr{Ng}2h_{jn}+A_jA_n\rg)
-
\fr12
g_i
\bar M
A_{jmn}
$$

\ses

\ses

$$
-
l_m
\fr1{h^2}g_i
\fr{q}B\lf(q+\fr12 gb\rg)
\fr 1{Ng}
\fr1N
 \lf(
A_jA_n  + \fr {N^2g^2}4h_{jn}
\rg).
$$

{

\ses

Summing all the addends yields
$$
g_i\D{A_{jmn}}g
+
\breve N^k{}_i
\D{A_{jmn}}{y^k}
+
A_{jkm}\breve N^k{}_{in}
+
A_{jkn}\breve N^k{}_{im}
+
A_{mnk}\breve N^k{}_{ij}
$$

\ses

\ses

$$
=
g_i
\lf(\fr1g-2\fr{bq}B\rg)
A_{mnj}
-
g_i
\fr{2bq}B
\fr{1}{N}
\fr{2}{gN}
\fr{2}{gN}
A_jA_mA_n
$$

\ses

          \ses

$$
-
\fr1{h^2}g_i
\fr{q}B\lf(q+\fr12 gb\rg)
\fr 1{Ng}
\Biggl[
-
\fr{Ng^2}4
(h_{jm}l_n+h_{mn}l_j+h_{nj}l_m)
-
\fr1N
(A_jA_ml_n+A_mA_jl_n+A_nA_jl_m)
$$

\ses

\ses

$$
+  g^2
 \lf(h_{mn} A_j  +h_{mj} A_n+h_{nj} A_m\rg)
\Biggr]
 $$

\ses

\ses

             \ses

$$
+
3
\fr1{h^2}g_i
\fr{1}{2B}
\lf(q+\fr12 gb\rg)
( b+gq)
A_{jmn}
$$

\ses

\ses

$$
+
\fr1{h^2}g_i
\fr{q^2}{2B}
\lf(
\fr bq+\fr12g
\rg)
\fr1N
 \lf(h_{mn} A_j  +h_{mj} A_n+h_{nj} A_m\rg)
+
\fr3{h^2}g_i
\fr{q^2}{2B}
\lf(
\fr bq+\fr12g
\rg)
\fr4{N^3g^2}
A_mA_jA_n
$$

\ses

\ses

$$
-
\fr1{h^2}g_i
\fr{q}B\lf(q+\fr12 gb\rg)
\fr 1{Ng}
\fr1N
(A_jA_ml_n+A_mA_jl_n+A_nA_jl_m)
$$

\ses

\ses

$$
-
\fr1{h^2}g_i
\fr{q}B\lf(q+\fr12 gb\rg)
\fr g4
(h_{jm}l_n+h_{mn}l_j+h_{nj}l_m),
$$

{

\nin
or
\ses\\
$$
g_i\D{A_{jmn}}g
+
\breve N^k{}_i
\D{A_{jmn}}{y^k}
+
A_{jkm}\breve N^k{}_{in}
+
A_{jkn}\breve N^k{}_{im}
+
A_{mnk}\breve N^k{}_{ij}
$$

\ses

\ses

$$
=
g_i
\lf(\fr1g-2\fr{bq}B\rg)
A_{mnj}
-
g_i
\fr{2bq}B
\fr{1}{N}
\fr{2}{gN}
\fr{2}{gN}
A_jA_mA_n
$$

\ses

          \ses

\ses

$$
 -
\fr1{h^2}g_i
\fr{g^2q}B\lf(q+\fr12 gb\rg)
\fr 1{Ng}
 \lf(h_{mn} A_j  +h_{mj} A_n+h_{nj} A_m\rg)
 $$

\ses

\ses

             \ses

$$
+
\fr1{h^2}g_i
\fr{1}{2B}
\lf[
6h^2bq+3gB   -b\lf(3q+\fr32 gb\rg)
\rg]
A_{jmn}
$$

\ses

\ses

$$
+
\fr1{h^2}g_i
\fr{gq}{2B}
\lf(
 b+\fr g2q
\rg)
\fr1{Ng}
 \lf(h_{mn} A_j  +h_{mj} A_n+h_{nj} A_m\rg)
+
\fr3{h^2}g_i
\fr{gq}{2B}
\lf(
 b+\fr g2q
\rg)
\fr4{N^3g^3}
A_mA_jA_n
$$

\ses

\ses

\ses

         \ses

$$
=
g_i
\fr1g
A_{jmn}
+
\fr1{h^2}g_i
\fr{1}{2B}
\lf[
4h^2bq+3gB   -b\lf(3q+\fr32 gb\rg)
\rg]
A_{jmn}
$$

\ses

\ses

$$
+
\fr1{h^2}g_i
\fr{q}{2B}
\lf(
 b-2h^2b
 -\fr32gq-g^2b
\rg)
A_{jmn}
=
g_i
\fr1g
A_{jmn}
+
\fr3{4h^2}g_i
A_{jmn}.
$$

{

By the help of such evaluations we eventually obtain
\ses\\
\be
g_i\D{A_{jmn}}g
+
\breve N^k{}_i
\D{A_{jmn}}{y^k}
+
A_{jkm}\breve N^k{}_{in}
+
A_{jkn}\breve N^k{}_{im}
+
A_{mnk}\breve N^k{}_{ij}
=
g_i
\fr1{gh^2}
A_{jmn}
+
\fr1{h^2}g_i
A_{jmn},
\ee
\ses\\
which shows that  the representation
\be
\breve N^k{}_{imn}
=
-
\fr g{2h^2}g_i
\fr1K
h_{mn}
 l^k
-
\fr1{gh^2}g_i
\fr1K
A^k{}_{mn}
\ee
\ses\\
indicated in (I.3.30) is valid.

The full coefficients
 $
 N^k{}_{imn}    =   N^{{{\rm I}}k}{}_{imn}+\breve N^k{}_{imn}
 $
 can be obtained on taking into account the
components
(A.29) together with
the representation
$
 N^{{{\rm I}}k}{}_{imn}
 =
-
(1/K)
\cD_i
A^k{}_{mn}
$
obtainable in the $(g=\const)$-case
(see [10,11]).
The result reads
\be
N^k{}_{imn}
=
\fr2hh_i
\fr1K
l^kh_{mn}
-
\fr1K
\cD_i
A^k{}_{mn}.
\ee
Thus the representation indicated in (I.3.31) is also valid.

{

\ses


\ses

\ses

\setcounter{equation}{0}

\ses

\ses

\nin {\bf Appendix B: Conformal property of the tangent Riemannian space}

\ses

\ses

Given an arbitrary Finsler space of  any dimension
$N\ge3$.
At any fixed point $x$,
the Riemannian curvature tensor
$
\wh R_{\{x\}}=\{\wh R_n{}^m{}_{ij}(x,y)\}
$
of the tangent Riemannian space
$\cR_{\{x\}}$
is
given by means  of  the
components
\be
\wh R_n{}^m{}_{ij}=
\fr1{F^2}
S_n{}^m{}_{ij},
\ee
\ses\\
where
\be
S_n{}^m{}_{ij}
=
\lf(
C^h{}_{nj}    C^m{}_{hi}
-
C^h{}_{ni}    C^m{}_{hj}
\rg)
F^2.
\ee

Let us construct the
 Weyl tensor $W_{ijmn}$ in the space
$\cR_{\{x\}}$,
so that
$$
F^2W_{ijmn}=S_{ijmn}
$$
\ses\\
\be
-
\fr1{N-2}(S_{im}g_{jn}  +  S_{jn}g_{im}
-
S_{in}g_{jm}  -  S_{jm}g_{in})
+
\fr 1{(N-1)(N-2)}
\breve S
(g_{im}g_{jn}-g_{in}g_{jm}),
\ee
\ses\\
where
$S_{ijmn}=g_{jh}S_{i}{}^j{}_{mn}$,
$
S_{im}=g^{jn}S_{ijmn}
$
and
$
\breve S=g^{im}S_{im}.
$
Contracting the tensor two times by the unit vector $l^n=(1/F)y^n$ yields
directly
$$
(N-2)
F^2
W_{ijmn}l^nl^j=
-
S_{im}
+
\fr 1{N-1}
\breve S
h_{im},
$$
\ses\\
where
$h_{im}=g_{im}-(1/F^2)y_iy_m$.
Therefore,
in any dimension $N\ge4$
the vanishing
$
W_{ijmn}=0
$
is tantamount to  the representation
\be
 S_{nmij} =C(h_{nj}h_{mi}-h_{ni}h_{mj}).
\ee
 It is known
 (see Section 5.8 in [1])
 that the indicatrix is a space of constant curvature if and only if the tensor
(B.2)
  fulfills the representation
(B.4),
in which case
$C=C(x)$
(that is, the  factor $C$  is independent of $y$).
The respective  indicatrix  curvature value
$
{\mathcal C}_{\text{Ind.}}
$
is given by
\be
{\mathcal C}_{\text{Ind.}}=1-C.
\ee

{

Next, in the dimension
$
N=3
$
the tensor $W_{ijmn}$
vanishes identically and, therefore,
the equality
\be
S_{ijmn}
=
L
(h_{im}h_{jn}-h_{in}h_{jm}) \quad \text{with} \quad L=\fr 12 \breve S
\ee
\ses\\
holds,
where $L$ may depend on $y$.
Taking
$
S_{im}
=
Lh_{im},
$
\ses
we should examine  the tensor
\be
C_{im}~:=\fr1{F^2}\lf(S_{im}-\fr14\breve Sg_{im}\rg)
\ee
\ses\\
of the Cotton-York type.
Let us use the Riemannian covariant derivative
$\cS$ operative in the space
$\cR_{\{x\}}$ under consideration.
Denoting $L_n=\partial L/\partial y^n$
and taking into account the vanishing
$\cS_ng_{im}=0$,
we have
\ses\\
$$
\cS_nC_{im}-\cS_mC_{in}
=
\fr1{F^2}
\lf(L_nh_{im}-L_mh_{in}
-
L\fr1F
(l_nh_{mi}-l_mh_{ni})
-\fr12
(L_ng_{im}-L_mg_{in})
\rg)
$$

\ses

\ses

$$
+
L
\fr1{F^3}
(l_n h_{im}-l_m h_{in}),
$$
\ses\\
which is
\ses
$$
\cS_nC_{im}-\cS_mC_{in}
=
\fr1{F^2}
\lf(
L_n \lf(h_{im}  -  \fr12  g_{im}  \rg)-L_m \lf(h_{in}  -  \fr12  g_{in}  \rg)
\rg),
$$
\ses\\
so that
\be
\cS_nC_{im}-\cS_mC_{in}=0
\ee
holds iff $L_n=0$,
that is when
 $\breve S=\breve S(x).$
The  vanishing (B.8) means the conformal flatness of the
three-dimensional
space
$\cR_{\{x\}}$.

{

Thus we are entitled to set forth the validity of the following proposition.

\ses

\ses

{\bf Proposition}.
{\it
Given an arbitrary Finsler space of  any dimension
$N\ge3$.
The  tangent Riemannian  space
$\cR_{\{x\}}$
is conformally flat if and only if the  indicatrix
is a space of constant curvature.

}

\ses

\ses

{

The question arises: What is the form of the conformal multiplier
of the  space
$\cR_{\{x\}}$
under study?
See the next appendix.

{


\ses

\ses

\setcounter{equation}{0}

\ses

\ses

\nin{\bf Appendix C:  Multiplier for the tangent Riemannian space}

\ses

\ses

To find the form of the conformal multiplier
of the  space
$\cR_{\{x\}}$
under study,
we can  start with
the conformal tensor
\be
u_{ij}
=
\fr{z(x,y)}
{(c_1(x))^2 F^{2a(x)}}g_{ij}
\ee
(cf. (II.2.3)),
where $z$ is a test smooth positive function.

Denoting
$
u_{ijk}=\partial u_{ij}/\partial{y^k},
$
\ses
we get
$$
(c_1)^2
u_{ijk}=
\fr1{F^{2a}}
     \lf(-  2z\fr aF l_k    +z_k\rg)
g_{ij}
+
2z
\fr1{F^{2a}}
C_{ijk},
$$
\ses\\
where
$C_{ijk}=(1/2)\partial g_{ij}/\partial y^k$
and $z_k=\partial z/\partial y^k$.

{

Constructing the coefficients
$$
Z_{ijk}~:=\fr12(u_{kji}+u_{iki}-u_{ijk})
$$
leads to
$$
F^{2a}
(c_1)^2 Z_{ijk}=
\lf(-  \fr {za}F l_i    + \fr12 z_i\rg)
g_{kj}
+
\lf(-  \fr {za}F l_j    + \fr12 z_j\rg)
g_{ik}
-
\lf(-  \fr {za}F l_k    + \fr12 z_k\rg)
g_{ij}
+
zC_{ijk}.
$$

{

Since the components $u^{ij}$ reciprocal to the components (C.1) are of the form
$
u^{ij}=(1/z)F^{2a}g^{ij}(c_1)^2 ,
$
the coefficients
$
Z^m{}_{ij}=u^{mh}Z_{ijh}
$
read merely
\ses\\
\be
Z^m{}_{ij}=
\lf(-  \fr aF l_i    + \fr1{2z} z_i\rg)
\de^m_j
+
\lf(-  \fr aF l_j    + \fr1{2z} z_j\rg)
\de^m_i
-
\lf(- \fr aF l^m    + \fr1{2z} g^{mk} z_k\rg)
g_{ij}
+
C^m{}_{ij}.
\ee

{

We straightforwardly obtain
$$
\D{Z^m{}_{ni}}{y^j}
=
\fr {a}{F^2}
l_j(l_n\de^m_i+l_i\de^m_n-l^mg_{ni})
$$

\ses

\ses

$$
-
\fr a{F^2}
\lf[h_{ij}\de^m_n+h_{nj}\de^m_i-h^m_jg_{in}-2
\lf(l^m -  \fr1{2z} \fr Fag^{mk} z_k \rg)
FC_{inj}
\rg]
+
\D{C^m{}_{ni}}{y^j}
$$

      \ses

\ses

$$
+
\fr12
\lf[
(\ln z)_{nj}  \de^m_i  +   (\ln z)_{ij}  \de^m_n
-
\D{\lf( g^{mk} (\ln z)_k\rg)}{y^j}
g_{ni}
\rg]
$$
\ses\\
and
\ses\\
$$
\D{Z^m{}_{ni}}{y^j}-
\D{Z^m{}_{nj}}{y^i}
=
\fr {a}{F^2}
[l_n(l_j\de^m_i-l_i\de^m_j) -l^m(l_jg_{ni}-l_ig_{nj})]
$$

\ses

\ses

\ses

$$
-
\fr a{F^2}
[(h_{nj}\de^m_i-h_{ni}\de^m_j)  -(h^m_jg_{in}-h^m_ig_{jn})]
+
\D{C^m{}_{ni}}{y^j}
-
\D{C^m{}_{nj}}{y^i}
$$

\ses

      \ses

\ses

$$
+
\fr12
\lf((\ln z)_{nj}  \de^m_i -(\ln z)_{ni}  \de^m_j \rg)
-
\fr12
\lf(
\D{\lf( g^{mk} (\ln z)_k\rg)}{y^j}  g_{ni}  - \D{\lf( g^{mk} (\ln z)_k \rg)}{y^i}  g_{nj} \rg),
$$

{

\ses

\nin
so that
$$
\D{Z^m{}_{ni}}{y^j}-    \D{Z^m{}_{nj}}{y^i}
=
-
\fr {2a}{F^2}
(h_{nj}h^m_i-h_{ni}h^m_j)
+
\D{C^m{}_{ni}}{y^j}
-
\D{C^m{}_{nj}}{y^i}
$$

\ses

\ses

      \ses

\ses

\be
+
\fr12
\lf((\ln z)_{nj}  \de^m_i -(\ln z)_{ni}  \de^m_j \rg)
-
\fr12
 g^{mh} \lf[(\ln z)_{hj}  g_{ni}    -    (\ln z)_{hi}  g_{nj}\rg]
+
(\ln z)_h
\lf[
C^{mh}{}_j  g_{ni}  - C^{mh}{}_i  g_{nj} \rg].
\ee

{

\ses

Also,
\ses\\
$$
Z^h{}_{ni}    Z^m{}_{hj}
-
Z^h{}_{nj}    Z^m{}_{hi}
=
-
 \fr a{F}
\Bigl[
-
\fr a{F}[l_n(l_i\de^m_j-l_j\de^m_i) +l^m(l_jg_{ni}-l_ig_{nj})]
+(l_iC^m{}_{nj}-l_jC^m{}_{ni})
\Bigr]
$$

\ses

\ses

 \ses

$$
-  \fr aF l_i
\de^h_n
\Biggl[
 \fr1{2z} z_h
\de^m_j
-
 \fr1{2z} g^{mk} z_k
g_{hj}
\Biggr]
-[ij]
$$

\ses

 \ses

\ses

$$
+
 \fr1{2z} z_i
\de^h_n
\Biggl[
\lf(-  \fr aF l_h    + \fr1{2z} z_h\rg)
\de^m_j
-
\lf(- \fr aF l^m    + \fr1{2z} g^{mk} z_k\rg)
g_{hj}
+
C^m{}_{hj}
\Biggr]
-[ij]
$$

 \ses

\ses

\ses

\ses

$$
-
\lf(\fr a{F}\rg)^2
(g_{in}\de^m_j-g_{jn}\de^m_i)
$$

\ses

\ses

\ses

$$
+ \fr aF l^h
g_{in}
\Biggl[
 \fr1{2z} z_h
\de^m_j
+
\fr1{2z} z_j
\de^m_h
-
 \fr1{2z} g^{mk} z_k
g_{hj}
\Biggr]
-[ij]
$$

 \ses

\ses

\ses

$$
-
 \fr1{2z} g^{hk} z_k
g_{in}
\Biggl[
\lf(-  \fr aF l_h    + \fr1{2z} z_h\rg)
\de^m_j
+
\lf(-  \fr aF l_j    + \fr1{2z} z_j\rg)
\de^m_h
-
\lf(- \fr aF l^m    + \fr1{2z} g^{mk} z_k\rg)
g_{hj}
\Biggr]
-[ij]
$$

 \ses

\ses

\ses

$$
-
 \fr1{2z} g^{hk} z_k
(g_{in}C^m{}_{hj}-g_{jn}C^m{}_{hi})
$$

\ses

\ses

\ses

$$
-
\fr a{F}
l_jC^m{}_{in}
+
C^h{}_{in}
\Biggl[
 \fr1{2z} z_h
\de^m_j
+
 \fr1{2z} z_j
\de^m_h
-
 \fr1{2z} g^{mk} z_k
g_{hj}
\Biggr]
-[ij]
$$

\ses

\ses

\ses

$$
+
C^h{}_{ni}    C^m{}_{hj}
-
C^h{}_{nj}    C^m{}_{hi},
$$

{

\ses

\nin
or
$$
Z^h{}_{ni}    Z^m{}_{hj}
-
Z^h{}_{nj}    Z^m{}_{hi}
=
-
\lf(\fr a{F}\rg)^2
(h_{in}h^m_j-h_{jn}h^m_i)
+
C^h{}_{ni}    C^m{}_{hj}
-
C^h{}_{nj}    C^m{}_{hi}
$$

\ses

\ses

\ses

$$
-  \fr aF l_i
\Biggl[
 \fr1{2z} z_n
\de^m_j
-
 \fr1{2z} g^{mk} z_k
g_{nj}
\Biggr]
+
 \fr1{2z} z_i
\lf(
-  \fr aF l_n    + \fr1{2z} z_n
\rg)
\de^m_j
-[ij]
$$

 \ses

\ses

\ses

$$
+ \fr aF l^h
g_{in}
\Biggl[
 \fr1{2z} z_h
\de^m_j
-
 \fr1{2z} g^{mk} z_k
g_{hj}
\Biggr]
-[ij]
$$

 \ses

\ses

\ses

$$
-
 \fr1{2z} g^{hs} z_s
g_{in}
\Biggl[
\lf(
-
  \fr aF l_h    + \fr1{2z} z_h\rg)
\de^m_j
-
\fr aF l_j
\de^m_h
+
 \fr aF l^m
 g_{hj}
-
\fr1{2z} g^{mk} z_k
g_{hj}
\Biggr]
-[ij]
$$

 \ses

\ses

\ses

$$
+
 \fr1{2z}  z_h
\lf[
-
C^{hm}{}_{j}   g_{in}  +  C^{hm}{}_{i}   g_{jn}
+
 C^h{}_{in}
\de^m_j
-
 C^h{}_{jn}
\de^m_i
\rg].
$$

{

\ses

In this way we come to
$$
Z^h{}_{ni}    Z^m{}_{hj}
-
Z^h{}_{nj}    Z^m{}_{hi}
=
-
\lf(\fr a{F}\rg)^2
(h_{in}h^m_j-h_{jn}h^m_i)
+
C^h{}_{ni}    C^m{}_{hj}
-
C^h{}_{nj}    C^m{}_{hi}
$$

\ses

\ses

\ses

$$
-  \fr aF
 \fr1{2z} z_n
(l_i\de^m_j-l_j\de^m_i)
+
  \fr aF
 \fr1{2z} g^{mk} z_k
(l_ig_{nj}-l_jg_{ni})
$$

\ses

 \ses

$$
+
 \fr1{2z}
\lf(
-  \fr aF l_n    + \fr1{2z} z_n
\rg)
(z_i\de^m_j-z_j\de^m_i)
$$

 \ses

\ses

\ses

$$
+ \fr aF
 \fr1{z} (l^sz_s)
(g_{in} \de^m_j-g_{jn} \de^m_i)
-
 \fr1{4z^2}
  \bigl(z_h g^{hs} z_s\bigr)
(
g_{in}\de^m_j-g_{jn}\de^m_i
)
$$

\ses

\ses

$$
+
 \fr1{2z}
\lf(
- \fr aF l^m
+
\fr1{2z} g^{mk} z_k
\rg)
(z_jg_{in}  -  z_ig_{jn})
$$

 \ses

\ses

\ses

\be
+
 \fr1{2z}  z_h
\lf[
-
C^{hm}{}_{j}   g_{in}  +  C^{hm}{}_{i}   g_{jn}
+
 C^h{}_{in}
\de^m_j
-
 C^h{}_{jn}
\de^m_i
\rg].
\ee

{

\ses

Thus we are able to evaluate
the curvature tensor
$$
\wt R_n{}^m{}_{ij} ~:
=
\D{Z^m{}_{ni}}{y^j}-
\D{Z^m{}_{nj}}{y^i}
+Z^h{}_{ni}    Z^m{}_{hj}
-
Z^h{}_{nj}    Z^m{}_{hi}.
$$
\ses\\
By  lowering  the index
$$
\wt R_{nmij} ~:
=
u_{mt}\wt R_n{}^t{}_{ij},
$$

\ses

{

\nin
we obtain the representation
\ses\\
$$
\fr1z
(c_1(x))^2 F^{2a(x)}
\wt R_{nmij}
=
\fr 1{F^2}
(a^2-2a)
(h_{nj}h_{mi}-h_{ni}h_{mj})
+
\fr 1{F^2}
S_{nmij}
$$

\ses

\ses

\ses

$$
+
\fr12
\Bigl((\ln z)_{nj} g_{mi} -(\ln z)_{ni} g_{mj}
-
(\ln z)_{mj}  g_{ni}    +    (\ln z)_{mi}  g_{nj}\Bigr)
$$

\ses

\ses

\ses

$$
-  \fr aF
 \fr1{2z}
\Bigl(
z_n(
l_ih_{mj}-l_jh_{mi})
- z_m
(l_ih_{nj}-l_jh_{ni})
\Bigr)
$$

\ses

 \ses

\ses

$$
+
 \fr1{2z}
\lf(
-  \fr aF l_n    + \fr1{2z} z_n
\rg)
(z_ig_{mj}-z_jg_{mi})
-
 \fr1{2z}
\lf(
- \fr aF l_m
+
\fr1{2z}  z_m
\rg)
(  z_ig_{jn} -z_jg_{in})
$$

 \ses

\ses

\ses

$$
+ \fr aF
 \fr1{z} (l^sz_s)
(g_{in} g_{mj}-g_{jn} g_{mi})
-
 \fr1{4z^2}
  \bigl(z_h g^{hs} z_s\bigr)
(g_{in} g_{mj}-g_{jn} g_{mi})
$$

\ses

\ses

\ses

\be
+
 \fr1{2z}  z_h
\Bigl(
C^{h}{}_{mj}   g_{in}  -  C^{h}{}_{mi}   g_{jn}
-
 C^h{}_{nj}
g_{mi}
+
 C^h{}_{in}
g_{mj}
\Bigr),
\ee
\ses\\
where
\ses\\
$$
S_n{}^m{}_{ij}
=
\lf(
\D{C^m{}_{ni}}{y^j}
-
\D{C^m{}_{nj}}{y^i}
+
C^h{}_{ni}    C^m{}_{hj}
-
C^h{}_{nj}    C^m{}_{hi}
\rg)
F^2
$$
\ses
and
$$
S_{nmij}
=
g_{mt}S_n{}^t{}_{ij}.
$$

{

\ses

Henceforth we assume the zero-degree homogeneity of the function
$z(x,y)$ with respect to the argument $y$,
having the identities
\ses\\
\be
(\ln z)_{ni} l^i=
-
\fr1F
(\ln z)_{n}, \qquad
(\ln z)_{i} l^i=0.
\ee
\ses\\
By performing the contraction in (C.5), we get
\ses\\
$$
\fr1z
(c_1(x))^2 F^{2a(x)}
\wt R_{nmij}l^ml^j
=
\fr1{2F}
\Bigl(-(\ln z)_{n}l_i -(\ln z)_{ni}
-    (\ln z)_{i}  l_n\Bigr)
$$

\ses

\ses

\ses

$$
+
 \fr1{2z}
\lf(
-  \fr aF l_n    + \fr1{2z} z_n
\rg)
z_i
+
 \fr1{2z}
 \fr aF
  z_il_n
-
 \fr1{4z^2}
  \bigl(z_h g^{hs} z_s\bigr)
h_{in}
+
 \fr1{2z}  z_h
 C^h{}_{in},
$$

{

\ses

\nin
or
\ses\\
$$
2
\fr1z
(c_1(x))^2 F^{2a(x)}
\wt R_{nmij}l^ml^j
=
-
\fr1F
(\ln z)_{n}l_i
-(\ln z)_{ni}
-
\fr1F
 (\ln z)_{i}  l_n
+
 \fr1{2z^2}
 z_n
z_i
$$

\ses

\ses

\be
-
 \fr1{2z^2}
  \bigl(z_h g^{hs} z_s\bigr)
h_{in}
+
 \fr1{z}  z_h
 C^h{}_{in}.
\ee

{

\ses

The vanishing
\be
\wt R_{nmij}=0
\ee
\ses\\
holds when
$$
(\ln z)_{ni}=
-
\fr1F
(\ln z)_{n}l_i
-
\fr1F
(\ln z)_{i}  l_n
+
 \fr1{2z^2}
 z_n
z_i
-
 \fr1{2z^2}
  \bigl(z_h g^{hs} z_s\bigr)
h_{in}
+
 \fr1{z}  z_h
 C^h{}_{in},
$$
\ses\\
in which case from (C.5) we get
\ses\\
$$
\fr 1{F^2}
S_{nmij}
=
-
\fr 1{F^2}
(a^2-2a)
(h_{nj}h_{mi}-h_{ni}h_{mj})
$$

\ses

\ses

      \ses

\ses

$$
+
\fr12
\Bigl(
\fr1F(\ln z)_{n}l_j
+    \fr1F(\ln z)_{j}  l_n
-
 \fr1{2z^2}
 z_n
z_j
+
 \fr1{2z^2}
  \bigl(z_h g^{hs} z_s\bigr)
h_{jn}
\Bigr)
 g_{mi}
$$

\ses

\ses

$$
 -
\fr12
\Bigl(
\fr1F(\ln z)_{n}l_i
+    \fr1F(\ln z)_{i}  l_n
-
 \fr1{2z^2}
 z_n
z_i
+
 \fr1{2z^2}
  \bigl(z_h g^{hs} z_s\bigr)
h_{in}
\Bigr)
 g_{mj}
$$

\ses

\ses

$$
-
\fr12
\Bigl(
\fr1F(\ln z)_{m}l_j
+
\fr1F
(\ln z)_{j}  l_m
-
 \fr1{2z^2}
 z_m
z_j
+
 \fr1{2z^2}
  \bigl(z_h g^{hs} z_s\bigr)
h_{jm}
\Bigr)
 g_{ni}
$$

\ses

\ses

$$
+
\fr12
\Bigl(
\fr1F
(\ln z)_{m}l_i
+
\fr1F
(\ln z)_{i}  l_m
-
\fr1{2z^2}
z_m
z_i
+
\fr1{2z^2}
\bigl(z_h g^{hs} z_s\bigr)
h_{im}
\Bigr)
g_{nj}
$$

\ses

\ses

\ses

$$
+  \fr aF
 \fr1{2z}
\Bigl(
z_n(
l_ih_{mj}-l_jh_{mi})
- z_m
(l_ih_{nj}-l_jh_{ni})
\Bigr)
$$

\ses

 \ses

$$
-
 \fr1{2z}
\lf(
-  \fr aF l_n    + \fr1{2z} z_n
\rg)
(z_ig_{mj}-z_jg_{mi})
+
 \fr1{2z}
\lf(
- \fr aF l_m
+
\fr1{2z}  z_m
\rg)
(  z_ig_{jn} -z_jg_{in})
$$

 \ses

\ses

\ses

$$
+ \fr1{4z^2}
  \bigl(z_h g^{hs} z_s\bigr)
(g_{in} g_{mj}-g_{jn} g_{mi}).
$$

{

\ses

Due simplifying yields
\ses\\
$$
\fr 1{F^2}
S_{nmij}
=
-
\fr 1{F^2}
(a^2-2a)
(h_{nj}h_{mi}-h_{ni}h_{mj})
$$

\ses

\ses

      \ses

\ses

$$
+
\fr1{2F}
\Bigl(
(\ln z)_{n}l_j
+    (\ln z)_{j}  l_n
\Bigr)
 g_{mi}
 -
\fr1{2F}
\Bigl(
(\ln z)_{n}l_i
+    (\ln z)_{i}  l_n
\Bigr)
 g_{mj}
$$

\ses

\ses

$$
-
\fr1{2F}
\Bigl(
(\ln z)_{m}l_j
+    (\ln z)_{j}  l_m
\Bigr)
 g_{ni}
+
\fr1{2F}
\Bigl(
(\ln z)_{m}l_i
+
(\ln z)_{i}  l_m
\Bigr)
g_{nj}
$$

\ses

\ses

\ses

$$
+
 \fr1{4z^2}
  \bigl(z_h g^{hs} z_s\bigr)
\Bigl[
h_{jn}   g_{mi}
-
h_{in}
 g_{mj}
-
h_{jm}
 g_{ni}
+
h_{im}
g_{nj}
 \Bigr]
$$

\ses

\ses

 \ses

$$
+  \fr aF
 \fr1{2z}
\Bigl(
z_n(
l_ih_{mj}-l_jh_{mi})
- z_m
(l_ih_{nj}-l_jh_{ni})
\Bigr)
$$

\ses

 \ses

\ses

$$
+
 \fr aF
 \fr1{2z}
\Bigl[
 l_n
(z_ih_{mj}-z_jh_{mi})
-
 l_m
(  z_ih_{jn} -z_jh_{in})
\Bigr]
+ \fr1{4z^2}
  \bigl(z_h g^{hs} z_s\bigr)
(g_{in} g_{mj}-g_{jn} g_{mi})
$$

\ses

\ses

\ses

\ses

\ses

$$
=
-
\fr 1{F^2}
(a^2-2a)
(h_{nj}h_{mi}-h_{ni}h_{mj})
$$

\ses

\ses

      \ses

$$
+  \fr {a-1}
{2zF}
\Bigl(
z_n(
l_ih_{mj}-l_jh_{mi})
- z_m
(l_ih_{nj}-l_jh_{ni})
+
 l_n
(z_ih_{mj}-z_jh_{mi})
-
 l_m
(  z_ih_{jn} -z_jh_{in})
\Bigr)
$$

 \ses

\ses

\ses

$$
+ \fr1{4z^2}
  \bigl(z_h g^{hs} z_s\bigr)
\Bigl(h_{jn}   g_{mi}
-
h_{in}
 g_{mj}
-
h_{jm}
 g_{ni}
+
h_{im}
g_{nj}
+
g_{in} g_{mj}-g_{jn} g_{mi}
\Bigr),
$$

{

\nin
which  is
$$
S_{nmij}
=
a(2-a)
(h_{nj}h_{mi}-h_{ni}h_{mj})
+
F^2
 \fr1{2z^2}
  \bigl(z_h g^{hs} z_s\bigr)
(h_{nj}h_{mi}-h_{ni}h_{mj})
$$

\ses

\ses

      \ses

\be
+  \fr {a-1}
{2z}
\Bigl(
z_n(
l_ih_{mj}-l_jh_{mi})
- z_m
(l_ih_{nj}-l_jh_{ni})
+
 l_n
(z_ih_{mj}-z_jh_{mi})
-
 l_m
(  z_ih_{jn} -z_jh_{in})
\Bigr)
F.
\ee

\ses

Therefore, the known vanishing
$
S_{nmij}l^i=0
$
requires
\be
(a-1)\Bigl((\ln z)_{n} h_{mj} -
(\ln z)_{m} h_{nj}
\Bigr)
=0.
\ee
\ses\\
Whenever $a\ne1$,
we should take
 $z_n=0$, which means that the
 {\it  function
 $z$ is independent of $y$.
}

{

\ses

Lastly,
it is worth noting that the case $a=1$  would mean
\be
 S_{nmij}=h_{nj}h_{mi}-h_{ni}h_{mj}.
\ee
In this case
\be
{\mathcal C}_{\text{Ind.}}=0
\ee

                  {


\ses

\ses

\setcounter{equation}{0}

\ses

\ses

\nin {\bf Appendix D:  Evaluation of  the coefficients   $N^k{}_{mni}$ in the ${\cF}^N$-space}

\ses

\ses

To  evaluate the coefficients
$N^k{}_{mni}=\partial{N^k{}_{mn}}/\partial{y^i}$,
we use  (II.3.32) and obtain
\ses\\
$$
N^k{}_{mni}
=
\fr1{F^2}l_i
h^k_n\D{F}{x^m}
+
\fr1{F^2}
(l^kh_{ni}+l_n h^k_i)\D{F}{x^m}
-
\fr1F
h^k_n\D{l_i}{x^m}
  -
\fr1F h^k_i\D{l_n}{x^m}
  -
l^k\D{\lf(\fr1Fh_{ni}\rg)}{x^m}
$$

\ses

\ses

$$
-
\D{ C^k{}_{ns}}{y^i}
  N^s{}_m
-
 C^k{}_{ns}
  N^s{}_{mi}
-
\fr1{F^2}
l_i
\lf(
l_nh^k_s
-
(1-H)
l^kh_{ns}
\rg)
 N^s{}_m
$$

\ses

\ses

$$
+
\fr1{F}
\lf(
l_nh^k_s
-
(1-H)
l^kh_{ns}
\rg)
 N^s{}_{mi}
+
\fr1{F^2}
\lf(
h_{ni}h^k_s
-
l_nl^kh_{si}-l_nl_sh^k_i
\rg)
 N^s{}_m
$$

\ses

\ses

$$
-
\fr1{F}
(1-H)
\lf(\fr1Fh^k_ih_{ns}
+
2l^kC_{nsi}
\rg)
 N^s{}_m
+
\fr1{F^2}
(1-H)
l^k
l_nh_{si}
 N^s{}_m
+
\fr1{F^2}
(1-H)
l^k
l_sh_{ni}
 N^s{}_m
$$

\ses

\ses

$$
-
\lf(y^k_{hp}t^p_i+\fr HFy^k_hl_i\rg)
F^H
 \lf( \D{ U^h_n}{x^m}+ L^h{}_{ms}U^s_n\rg)
-
y^k_hF^H
 \lf( \D{ U^h_{ni}}{x^m}+ L^h{}_{ms}U^s_{ni}\rg),
  $$

{

\nin
or
\ses\\
$$
N^k{}_{mni}
=
\fr1{F^2}l_i
h^k_n\D{F}{x^m}
+
\fr1{F^2}
l_n h^k_i
\D{F}{x^m}
-
\fr1F
h^k_n\D{l_i}{x^m}
  -
\fr1F h^k_i\D{l_n}{x^m}
+
\fr1{F^2}
l^kh_{ni}\D{F}{x^m}
  -
\fr1Fl^k\D{h_{ni}}{x^m}
$$

\ses

\ses

$$
-
\D{ C^k{}_{ns}}{y^i}
  N^s{}_m
-
 C^k{}_{ns}
  N^s{}_{mi}
-
\fr1{F^2}
l_i
\lf(
l_nh^k_s
-
(1-H)
l^kh_{ns}
\rg)
 N^s{}_m
$$

\ses

\ses

$$
+
\fr1{F}
\lf(
l_nh^k_s
-
(1-H)
l^kh_{ns}
\rg)
 N^s{}_{mi}
+
\fr1{F^2}
\lf(
h_{ni}h^k_s
-
l_nl^kh_{si}-l_nl_sh^k_i
\rg)
 N^s{}_m
$$

\ses

\ses

$$
-
\fr1{F^2}
(1-H)
h^k_ih_{ns}
 N^s{}_m
-
\fr2{F}
(1-H)
l^kC_{nsi}
 N^s{}_m
$$

\ses

\ses

$$
+
\fr1{F^2}
(1-H)
l^k
l_nh_{si}
 N^s{}_m
-
H
\fr1{F^2}
l^k
l_sh_{ni}
 N^s{}_m
$$

\ses

\ses

$$
-
\lf(y^k_{hp}t^p_i+\fr HFy^k_hl_i\rg)
F^H
 \lf(
 - N^s{}_mU^h_{ns}
-N^s{}_{mn}U^h_{s}
\rg)
-
y^k_hF^H
 \lf( \D{ U^h_{ni}}{x^m}+ L^h{}_{ms}U^s_{ni}\rg).
  $$

{

Reducing similar terms leads to
\ses\\
$$
N^k{}_{mni}
=
\fr1{F^2}l_i
h^k_n\D{F}{x^m}
+
\fr1{F^2}
l_n h^k_i
\D{F}{x^m}
-
\fr1F
h^k_n\D{l_i}{x^m}
  -
\fr1F h^k_i\D{l_n}{x^m}
+
\fr1{F^2}
l^kh_{ni}\D{F}{x^m}
  -
\fr1Fl^k\D{h_{ni}}{x^m}
$$

\ses

\ses

$$
-
\D{ C^k{}_{ns}}{y^i}
  N^s{}_m
-
 C^k{}_{ns}
  N^s{}_{mi}
-
\fr1{F^2}
l_i
\lf(
l_nh^k_s
-
(1-H)
l^kh_{ns}
\rg)
 N^s{}_m
$$

\ses

\ses

$$
+
\fr1{F}
\lf(
l_nh^k_s
-
(1-H)
l^kh_{ns}
\rg)
 N^s{}_{mi}
+
\fr1{F^2}
\lf(
h_{ni}h^k_s
-
l_nl^kh_{si}-l_nl_sh^k_i
\rg)
 N^s{}_m
$$

\ses

\ses

$$
-
\fr1{F^2}
(1-H)
h^k_ih_{ns}
 N^s{}_m
-
\fr2{F}
(1-H)
l^kC_{nsi}
 N^s{}_m
$$

\ses

\ses

$$
+
\fr1{F^2}
(1-H)
l^k
l_nh_{si}
 N^s{}_m
-
H
\fr1{F^2}
l^k
l_sh_{ni}
 N^s{}_m
$$

\ses

\ses

$$
- y^k_{p} t^p_{iv}y^v_h
F^H
 N^s{}_mU^h_{ns}
-
h^v_s y^k_{p} t^p_{iv}
N^s{}_{mn}
$$

\ses

\ses

$$
+
\fr HFy^k_hl_i
F^H
 \lf(
 N^s{}_mU^h_{ns}
+
N^s{}_{mn}U^h_{s}
\rg)
-
y^k_hF^H
 \lf( \D{ U^h_{ni}}{x^m}+ L^h{}_{ms}U^s_{ni}\rg).
  $$

\ses

Applying here  (II.3.30) and (II.3.31) yields

{

\ses

$$
N^k{}_{mni}
=
\fr1{F^2}l_i
h^k_n\D{F}{x^m}
+
\fr1{F^2}
l_n h^k_i
\D{F}{x^m}
-
\fr1F
h^k_n\D{l_i}{x^m}
  -
\fr1F h^k_i\D{l_n}{x^m}
+
\fr1{F^2}
l^kh_{ni}\D{F}{x^m}
  -
\fr1Fl^k\D{h_{ni}}{x^m}
$$

\ses

\ses

$$
-
\D{ C^k{}_{ns}}{y^i}
  N^s{}_m
-
 C^k{}_{ns}
  N^s{}_{mi}
-
 C^k{}_{is}
  N^s{}_{mn}
-
\fr1{F^2}
l_i
\lf(
l_nh^k_s
-
(1-H)
l^kh_{ns}
\rg)
 N^s{}_m
$$

\ses

\ses

$$
+
\fr1{F}
\lf(
l_nh^k_s
-
(1-H)
l^kh_{ns}
\rg)
 N^s{}_{mi}
+
\fr1{F^2}
\lf(
h_{ni}h^k_s
-
l_nl^kh_{si}-l_nl_sh^k_i
\rg)
 N^s{}_m
$$

\ses

\ses

$$
-
\fr1{F^2}
(1-H)
h^k_ih_{ns}
 N^s{}_m
-
\fr2{F}
(1-H)
l^kC_{nsi}
 N^s{}_m
$$

\ses

\ses

$$
+
\fr1{F^2}
(1-H)
l^k
l_nh_{si}
 N^s{}_m
-
H
\fr1{F^2}
l^k
l_sh_{ni}
 N^s{}_m
$$

\ses

\ses

$$
    -
\lf(C^k{}_{iv}
-(1-H)\fr1{F}(l_v\de^k_i+l_i\de^k_v-l^kg_{iv})
\rg)
                         y^v_hF^HU^h_{ns}
 N^s{}_m
+
(1-H)\fr1{F}(l_ih^k_s-l^kh_{is})
N^s{}_{mn}
$$

\ses

\ses

$$
+
\fr HFy^k_hl_i
F^H
 \lf(
 N^s{}_mU^h_{ns}
+
N^s{}_{mn}U^h_{s}
\rg)
-
y^k_hF^H
 \lf( \D{ U^h_{ni}}{x^m}+ L^h{}_{ms}U^s_{ni}\rg).
  $$

{

From the representation
$$
 U^h_n
 =
\fr1{F^H} t^h_n-\fr1FH U^hl_n
 $$
(see (II.3.5))
let us evaluate
the coeffcients
$$
 U^h_{ni}=\D{U^h_n}{y^i}.
 $$
We get
$$
 U^h_{ni}=
-\fr HFl_i
\fr1{F^H} t^h_n
+
\fr1{F^H} t^h_{ni}
+\fr1{F^2}l_iH U^hl_n
-\fr1FH U^h_il_n
-\fr1{F^2}H U^hh_{ni}
$$

\ses

\ses

$$
=
-\fr HF(l_iU^h_n+l_nU^h_i)
+
\fr1{F^H} t^h_{ni}
+
\fr1{F^2}(1-H)H U^hl_nl_i
-\fr1{F^2}H U^hh_{ni}.
$$

Here,
$$
 t^h_{ni}
=C^s{}_{ni}t^h_s
-(1-H)\fr1{F}
\lf(l_it^h_n+l_nt^h_i-\fr HF t^hg_{ni}\rg)
$$
(see (II.3.30)). From this it follows that
\ses\\
$$
 U^h_{ni}=
-\fr HF(l_iU^h_n+l_nU^h_i)
+
\fr1{F^H}
\lf[        C^s{}_{ni}t^h_s
-(1-H)\fr1{F}
\lf(l_it^h_n+l_nt^h_i\rg)   \rg]
$$

\ses

\ses

$$
+
\fr1{F^2}
\fr1{F^H}
(1-H)
 H t^hg_{ni}
-
\fr1{F^2}H U^h
(g_{ni}-(2-H)l_nl_i).
$$

\ses

We have here
\ses\\
$$
 U^h_{ni}=
-\fr HF(l_iU^h_n+l_nU^h_i)
+
\fr1{F^H}
\lf[        C^s{}_{ni}t^h_s
-\fr1{F}
\lf(l_it^h_n+l_nt^h_i\rg)   \rg]
+
\fr1{F^H}
H\fr1{F}
\lf(l_it^h_n+l_nt^h_i\rg)
$$

\ses

\ses

$$
-
\fr1{F^2}
\fr1{F^H}
 H^2 t^hg_{ni}
+
\fr{2-H}{F^2}H U^h
l_nl_i
$$

\ses

\ses

$$
=
        C^s{}_{ni}t^h_s  \fr1{ F^H}
-\fr1{F}
\fr1{ F^H}
\lf(l_it^h_n+l_nt^h_i\rg)
-
\fr1{F^2}
 H^2 U^hg_{ni}
+
\fr{2+H}{F^2}H U^h
l_nl_i.
$$

Thus it is valid that
\ses\\
\be
U^h_{ni}=
        C^s{}_{ni}t^h_s  \fr1{ F^H}
-\fr1{F}
\lf(l_iU^h_n+l_nU^h_i\rg)
-
\fr1{F^2}
 H^2 U^hh_{ni}.
\ee

\ses

We may readily deduce
the contraction
\ses\\
\be
 y^v_hF^HU^h_{ns}
=
       C^v{}_{sn}
-\fr1{F}
\lf(l_sh^v_n+l_nh^v_s\rg)
-
\fr1{F}
 Hh_{ns}
l^v.
\ee

{

Inserting yields

\ses

$$
N^k{}_{mni}
=
\fr1{F^2}l_i
h^k_n\D{F}{x^m}
+
\fr1{F^2}
l_n h^k_i
\D{F}{x^m}
-
\fr1F
h^k_n\D{l_i}{x^m}
  -
\fr1F h^k_i\D{l_n}{x^m}
+
\fr1{F^2}
l^kh_{ni}\D{F}{x^m}
  -
\fr1Fl^k\D{h_{ni}}{x^m}
$$

\ses

\ses

$$
-
\D{ C^k{}_{ns}}{y^i}
  N^s{}_m
-
 C^k{}_{ns}
  N^s{}_{mi}
-
 C^k{}_{is}
  N^s{}_{mn}
-
\fr1{F^2}
l_i
\lf(
l_nh^k_s
-
(1-H)
l^kh_{ns}
\rg)
 N^s{}_m
$$

\ses

\ses

$$
+
\fr1{F}
\lf(
l_nh^k_s
-
(1-H)
l^kh_{ns}
\rg)
 N^s{}_{mi}
+
\fr1{F^2}
\lf(
h_{ni}h^k_s
-
l_nl^kh_{si}-l_nl_sh^k_i
\rg)
 N^s{}_m
$$

\ses

\ses

$$
-
\fr1{F^2}
(1-H)
h^k_ih_{ns}
 N^s{}_m
-
\fr2{F}
(1-H)
l^kC_{nsi}
 N^s{}_m
$$

\ses

\ses

$$
+
\fr1{F^2}
(1-H)
l^k
l_nh_{si}
 N^s{}_m
-
H
\fr1{F^2}
l^k
l_sh_{ni}
 N^s{}_m
$$

\ses

\ses

$$
    -
\lf(C^k{}_{iv}
-(1-H)\fr1{F}(l_v\de^k_i+l_i\de^k_v-l^kg_{iv})
\rg)
           \lf[       C^v{}_{sn}
-\fr1{F}
\lf(l_sh^v_n+l_nh^v_s\rg)
-
\fr1{F}
 Hh_{ns}
l^v
               \rg]
 N^s{}_m
$$

\ses

\ses

$$
+
(1-H)\fr1{F}(l_ih^k_s-l^kh_{is})
N^s{}_{mn}
+
\fr HFl_i
            \lf[       C^k{}_{sn}
-\fr1{F}
\lf(l_sh^k_n+l_nh^k_s\rg)
-
\fr H{F}
h_{ns}
l^k
               \rg]
 N^s{}_m
+
\fr HFl_i
N^s{}_{mn}h^k_s
$$

\ses

\ses

   \ses

$$
+
y^k_hF^H
\D{
\fr1{F}
\lf(l_iU^h_n+l_nU^h_i\rg)
}
{x^m}
+
y^k_hF^H
\D{
\fr1{F^2}
 H^2 U^hh_{ni}
}
{x^m}
  $$

\ses

\ses

$$
-
\D{C^k{}_{ni}}{x^m}
-
 C^s{}_{ni}
 y^k_hF^H
\D{      \lf( U^h_s  +\fr1F HU^hl_s\rg)}{x^m}
-
y^k_hF^H
 L^h{}_{ms}U^s_{ni},
  $$

{

\nin
or
\ses\\
$$
N^k{}_{mni}
=
\fr1{F^2}l_i
h^k_n\D{F}{x^m}
+
\fr1{F^2}
l_n h^k_i
\D{F}{x^m}
-
\fr1F
h^k_n\D{l_i}{x^m}
  -
\fr1F h^k_i\D{l_n}{x^m}
+
\fr1{F^2}
l^kh_{ni}\D{F}{x^m}
  -
\fr1Fl^k\D{h_{ni}}{x^m}
$$

\ses

\ses

$$
-
\D{C^k{}_{ni}}{x^m}
-
\D{ C^k{}_{ns}}{y^i}
  N^s{}_m
-
 C^k{}_{ns}
  N^s{}_{mi}
-
 C^k{}_{is}
  N^s{}_{mn}
-
\fr1{F^2}
l_i
\lf(
l_nh^k_s
-
(1-H)
l^kh_{ns}
\rg)
 N^s{}_m
$$

\ses

\ses

$$
+
\fr1{F}
\lf(
l_nh^k_s
-
(1-H)
l^kh_{ns}
\rg)
 N^s{}_{mi}
+
\fr1{F^2}
\lf(
h_{ni}h^k_s
-
l_nl^kh_{si}-l_nl_sh^k_i
\rg)
 N^s{}_m
$$

\ses

\ses

$$
-
\fr1{F^2}
(1-H)
h^k_ih_{ns}
 N^s{}_m
-
\fr2{F}
(1-H)
l^kC_{nsi}
 N^s{}_m
$$

\ses

\ses

$$
+
\fr1{F^2}
(1-H)
l^k
l_nh_{si}
 N^s{}_m
-
H
\fr1{F^2}
l^k
l_sh_{ni}
 N^s{}_m
$$

\ses

\ses

$$
    -
           \lf[       C^v{}_{sn}C^k{}_{iv}
-\fr1{F}
\lf(l_sC^k{}_{in}+l_nC^k{}_{is}\rg)
               \rg]
 N^s{}_m
-
(1-H)\fr1{F^2}\de^k_i
 Hh_{ns}
 N^s{}_m
$$

\ses

\ses

$$
+
(1-H)\fr1{F}(l_i\de^k_v-l^kg_{iv})
           \lf[       C^v{}_{sn}
-\fr1{F}
\lf(l_sh^v_n+l_nh^v_s\rg)
               \rg]
 N^s{}_m
-
(1-H)\fr1{F}l^kh_{is}
N^s{}_{mn}
$$

\ses

\ses

$$
+
\fr HFl_i
            \lf[       C^k{}_{sn}
-\fr1{F}
\lf(l_sh^k_n+l_nh^k_s\rg)
-
\fr H{F}
h_{ns}
l^k
               \rg]
 N^s{}_m
+
\fr 1Fl_i
N^s{}_{mn}h^k_s
$$

\ses

\ses

   \ses

$$
+
y^k_hF^H
\D{
\fr1{F}
\lf(l_iU^h_n+l_nU^h_i\rg)
}
{x^m}
+
y^k_hF^H
\fr1{F^2}
 H^2h_{ni}
 \D{
 U^h
}
{x^m}
+
y^k
H
\D{
\fr1{F^2}
h_{ni}
}
{x^m}
+
2
y^k
H_m
\fr1{F^2}
 h_{ni}
  $$

\ses

\ses

$$
-
 C^s{}_{ni}
 y^k_hF^H
\D{      \lf(\fr1F HU^hl_s\rg)}{x^m}
-
y^k_hF^H
 L^h{}_{sm}U^s_{ni}
-
 C^v{}_{ni}
 y^k_hF^H
\D{       U^h_v}{x^m}.
  $$

{

Finally we apply  here (II.3.32),
obtaining
\ses\\
$$
N^k{}_{mni}
=
\fr1{F^2}l_i
h^k_n\D{F}{x^m}
+
\fr1{F^2}
l_n h^k_i
\D{F}{x^m}
-
\fr1F
h^k_n\D{l_i}{x^m}
  -
\fr1F h^k_i\D{l_n}{x^m}
+
\fr1{F^2}
l^kh_{ni}\D{F}{x^m}
  -
\fr1Fl^k\D{h_{ni}}{x^m}
$$

\ses

\ses

$$
-
\D{C^k{}_{ni}}{x^m}
-
\D{ C^k{}_{ns}}{y^i}
  N^s{}_m
-
 C^k{}_{ns}
  N^s{}_{mi}
-
 C^k{}_{is}
  N^s{}_{mn}
+
 N^k{}_{ms}
  C^s{}_{ni}
   $$

\ses

\ses

$$
-
\fr1{F^2}
l_i
\lf(
l_nh^k_s
-
(1-H)
l^kh_{ns}
\rg)
 N^s{}_m
+
\fr1{F^2}
l_n
(1-H)
l^kh_{is}
 N^s{}_m
$$

\ses

\ses

$$
+
\fr1{F}
h^k_s
(l_n N^s{}_{mi}+l_i N^s{}_{mn})
-
\fr1{F}
(1-H)
l^k
\lf(
h_{ns}N^s{}_{mi}
+
h_{is}N^s{}_{mn}
\rg)
$$

\ses

\ses

$$
+
\fr1{F^2}
\lf(
h_{ni}h^k_s
-
l_nl^kh_{si}-l_nl_sh^k_i
\rg)
 N^s{}_m
$$

\ses

\ses

$$
-
\fr1{F^2}
(1-H)
h^k_ih_{ns}
 N^s{}_m
-
\fr2{F}
(1-H)
l^kC_{nsi}
 N^s{}_m
-
H
\fr1{F^2}
l^k
l_sh_{ni}
 N^s{}_m
$$

\ses

\ses

$$
    -
           \lf[       C^v{}_{sn}C^k{}_{iv}
-\fr1{F}
\lf(l_sC^k{}_{in}+l_nC^k{}_{is}\rg)
               \rg]
 N^s{}_m
-
(1-H)\fr1{F^2}\de^k_i
 Hh_{ns}
 N^s{}_m
$$

\ses

\ses

$$
+
(1-H)\fr1{F}l_i
           \lf[       C^k{}_{sn}
-\fr1{F}
\lf(l_sh^k_n+l_nh^k_s\rg)
               \rg]
 N^s{}_m
-
(1-H)\fr1{F}l^k
           \lf[       C_{isn}
-\fr1{F}
\lf(l_sh_{in}+l_nh_{is}\rg)
               \rg]
 N^s{}_m
$$

\ses

\ses

$$
+
\fr HFl_i
            \lf[       C^k{}_{sn}
-\fr1{F}
\lf(l_sh^k_n+l_nh^k_s\rg)
-
\fr H{F}
h_{ns}
l^k
               \rg]
 N^s{}_m
$$

\ses

\ses

   \ses

$$
+
y^k_hF^H
\D{
\fr1{F}
\lf(l_iU^h_n+l_nU^h_i\rg)
}
{x^m}
+
y^k_hF^H
\fr1{F^2}
 H^2h_{ni}
 \D{
 U^h
}
{x^m}
+
y^k
H
\D{
\fr1{F^2}
h_{ni}
}
{x^m}
+
2
y^k
H_m
\fr1{F^2}
 h_{ni}
  $$
\ses

\ses

$$
-
 C^s{}_{ni}
 l^k
\D{     l_s}{x^m}
-
y^k_hF^H
 L^h{}_{sm}U^s_{ni}
  $$

\ses

\ses

$$
 -
 C^v{}_{ni}
\lf[
-
\fr1F
h^k_v\D{F}{x^m}
-
l^k\D{l_v}{x^m}
-
 C^k{}_{vs}
  N^s{}_m
-
\fr1{F}
(1-H)
l^kh_{vs}
 N^s{}_m
\rg]
+
 C^v{}_{ni}
 y^k_hF^H
 L^h{}_{sm}U^s_v.
  $$

{

Reducing similar terms yields now
\ses\\
$$
N^k{}_{mni}
=
\fr1{F^2}l_i
h^k_n\D{F}{x^m}
+
\fr1{F^2}
l_n h^k_i
\D{F}{x^m}
-
\fr1F
h^k_n\D{l_i}{x^m}
  -
\fr1F h^k_i\D{l_n}{x^m}
+
\fr1{F^2}
l^kh_{ni}\D{F}{x^m}
$$

\ses

\ses

$$
  -
(1-H)\fr1Fl^k
\lf[
\D{h_{ni}}{x^m}
+
2C_{nsi}
 N^s{}_m
+
h_{ns}N^s{}_{mi}
+
h_{is}N^s{}_{mn}
 \rg]
-
\cD_m C^k{}_{ni}
$$

\ses

\ses

$$
-
\fr1{F^2}
l_i
\lf(
l_nh^k_s
-
(1-H)
l^kh_{ns}
\rg)
 N^s{}_m
+
\fr1{F^2}
l_n
(1-H)
l^kh_{is}
 N^s{}_m
+
\fr1{F}
h^k_s
(l_n N^s{}_{mi}+l_i N^s{}_{mn})
$$

\ses

\ses

$$
+
\fr1{F^2}
\lf(
h_{ni}h^k_s
-
l_nl^kh_{si}-l_nl_sh^k_i
\rg)
 N^s{}_m
-
\fr1{F^2}
(1-H)
h^k_ih_{ns}
 N^s{}_m
-
H
\fr1{F^2}
l^k
l_sh_{ni}
 N^s{}_m
$$

\ses

\ses

$$
    -
           \lf[       C^v{}_{sn}C^k{}_{iv}
-\fr1{F}
\lf(l_iC^k{}_{ns}+l_nC^k{}_{is}\rg)
               \rg]
 N^s{}_m
-
(1-H)\fr1{F^2}h^k_i
 Hh_{ns}
 N^s{}_m
-
(1-H)\fr1{F^2}l^kl_i
 Hh_{ns}
 N^s{}_m
$$

\ses

\ses

$$
-
\fr1{F^2}l_i
\lf(l_sh^k_n+l_nh^k_s\rg)
 N^s{}_m
-
(1-H)\fr1{F}l^k
           \lf[       C_{isn}
-\fr1{F}
\lf(l_sh_{in}+l_nh_{is}\rg)
               \rg]
 N^s{}_m
-
\fr H{F^2}
l_il^k
 H
h_{ns}
 N^s{}_m
$$

\ses

\ses

   \ses

$$
+
y^k_hF^H
\D{
\fr1{F}
\lf(l_iU^h_n+l_nU^h_i\rg)
}
{x^m}
+
y^k_hF^H
\fr1{F^2}
 H^2h_{ni}
 \D{
 U^h
}
{x^m}
+
y^k
H
h_{ni}
\D{
\fr1{F^2}
}
{x^m}
+
2
y^k
H_m
\fr1{F^2}
 h_{ni}
  $$
\ses

\ses

$$
-
y^k_hF^H
 L^h{}_{sm}U^s_{ni}
 -
 C^v{}_{ni}
\lf[
-
 C^k{}_{vs}
  N^s{}_m
-
\fr1{F}
(1-H)
l^kh_{vs}
 N^s{}_m
\rg]
+
 C^v{}_{ni}
 y^k_hF^H
 L^h{}_{sm}U^s_v.
  $$

{

The next step is to transform the representation to
\ses\\
$$
N^k{}_{mni}
=
\fr1{F^2}l_i
h^k_n\D{F}{x^m}
+
\fr1{F^2}
l_n h^k_i
\D{F}{x^m}
-
\fr1F
h^k_n\D{l_i}{x^m}
  -
\fr1F h^k_i\D{l_n}{x^m}
+
\fr1{F^2}
l^kh_{ni}\D{F}{x^m}
$$

\ses

\ses

$$
  -
(1-H)\fr1Fl^k
\lf[
\D{h_{ni}}{x^m}
+
2C_{nsi}
 N^s{}_m
+
h_{ns}N^s{}_{mi}
+
h_{is}N^s{}_{mn}
 \rg]
-
\cD_m C^k{}_{ni}
$$

\ses

\ses

$$
-
\fr1{F^2}
l_i
\lf(
l_nh^k_s
-
(1-H)
l^kh_{ns}
\rg)
 N^s{}_m
+
\fr1{F^2}
l_n
(1-H)
l^kh_{is}
 N^s{}_m
+
\fr1{F}
h^k_s
(l_n N^s{}_{mi}+l_i N^s{}_{mn})
$$

\ses

\ses

$$
+
\fr1{F^2}
\lf(
h_{ni}h^k_s
-
l_nl^kh_{si}-l_nl_sh^k_i-l_il_sh^k_n
\rg)
 N^s{}_m
-
\fr1{F^2}
(1-H)
h^k_ih_{ns}
 N^s{}_m
-
H
\fr1{F^2}
l^k
l_sh_{ni}
 N^s{}_m
$$

\ses

\ses

$$
    -
\lf[
\fr1{F^2}
S_n{}^k{}_{si}
-
\fr1{F}
\lf(l_iC^k{}_{ns}+l_nC^k{}_{is}\rg)
               \rg]
 N^s{}_m
-
(1-H)\fr1{F^2}h^k_i
 Hh_{ns}
 N^s{}_m
$$

\ses

\ses

$$
-
\fr1{F^2}l^kl_i
 Hh_{ns}
 N^s{}_m
-
\fr1{F^2}l_i
l_nh^k_s
 N^s{}_m
+
(1-H)\fr1{F^2}l^k
l_sh_{in}
 N^s{}_m
+
(1-H)\fr1{F^2}l^k
l_nh_{is}
 N^s{}_m
$$

\ses

\ses

   \ses

$$
+
y^k_hF^H
\fr1{F}
\lf(
l_i\D{U^h_n}{x^m}
+
l_n\D{U^h_i}{x^m}
\rg)
+
h^k_n
\D{
\fr1{F}
l_i
}
{x^m}
+
h^k_i
\D{
\fr1{F}
l_n
}
{x^m}
 $$

\ses

\ses

$$
+
y^k_hF^H
\fr1{F^2}
 H^2h_{ni}
 \D{
 U^h
}
{x^m}
-
\fr2{F^2}
l^k
H
h_{ni}
\D{
F
}
{x^m}
+
2
\fr1F
l^k
H_m
 h_{ni}
-
y^k_hF^H
 L^h{}_{sm}U^s_{ni}
+
 C^v{}_{ni}
 y^k_hF^H
 L^h{}_{sm}U^s_v,
  $$
\ses
where
$$
S_n{}^k{}_{si}   =F^2(   C^v{}_{in}C^k{}_{sv}-   C^v{}_{sn}C^k{}_{iv}).
$$

{

We come to
\ses\\
$$
N^k{}_{mni}
=
  -
(1-H)\fr1Fl^k
\lf[
\D{h_{ni}}{x^m}
+
2C_{nsi}
 N^s{}_m
+
h_{ns}N^s{}_{mi}
+
h_{is}N^s{}_{mn}
 \rg]
-
\cD_m C^k{}_{ni}
$$

\ses

\ses

$$
-
\fr1{F^2}
l_i
\lf(
l_nh^k_s
-
(1-H)
l^kh_{ns}
\rg)
 N^s{}_m
+
\fr1{F^2}
l_n
(1-H)
l^kh_{is}
 N^s{}_m
$$

\ses

\ses

$$
+
\fr1{F}
h^k_s
(l_n N^s{}_{mi}+l_i N^s{}_{mn})
+
\fr1{F^2}
\lf(
h_{ni}h^k_s
-
l_nl_sh^k_i-l_il_sh^k_n
\rg)
 N^s{}_m
$$

\ses

\ses

$$
    -
\lf[
\fr1{F^2}
S_n{}^k{}_{si}
-\fr1{F}
\lf(l_iC^k{}_{ns}+l_nC^k{}_{is}\rg)
               \rg]
 N^s{}_m
-
(1-H^2)\fr1{F^2}h^k_i
h_{ns}
 N^s{}_m
$$

\ses

\ses

$$
-
\fr1{F^2}l_i
l_nh^k_s
 N^s{}_m
-H\fr1{F^2}l^k
(l_ih_{ns}+l_nh_{is})
 N^s{}_m
+
y^k_hF^H
\fr1{F}
\lf(
l_i\D{U^h_n}{x^m}
+
l_n\D{U^h_i}{x^m}
\rg)
 $$

\ses

\ses

$$
+
y^k_hF^H
\fr1{F^2}
 H^2h_{ni}
 \D{
 U^h
}
{x^m}
+
2
\fr1F
l^k
H_m
 h_{ni}
+
y^k_hF^H
 L^h{}_{sm}
\lf[
\fr1{F}
\lf(l_iU^s_n+l_nU^s_i\rg)
+
\fr{H^2}{F^2}
 U^sh_{ni}
 \rg]
  $$

\ses

\ses

\ses

\ses

$$
=
  -
(1-H)\fr1Fl^k
\cD_mh_{ni}
-
\cD_m C^k{}_{ni}
-
\fr1{F^2}
l_i
l_nh^k_s
 N^s{}_m
+
\fr1{F^2}
\lf(
h_{ni}h^k_s
-
l_nl_sh^k_i-l_il_sh^k_n
\rg)
 N^s{}_m
$$

\ses

\ses

$$
    -
  \lf[
\fr1{F^2}
S_n{}^k{}_{si}
-\fr1{F}
\lf(l_iC^k{}_{ns}+l_nC^k{}_{is}\rg)
               \rg]
 N^s{}_m
-
(1-H^2)\fr1{F^2}h^k_i
h_{ns}
 N^s{}_m
$$

\ses

\ses

$$
-
\fr1{F^2}l_i
l_nh^k_s
 N^s{}_m
-H\fr1{F^2}l^k
(l_ih_{ns}+l_nh_{is})
 N^s{}_m
$$

\ses

\ses

$$
+
y^k_hF^H
\fr1{F}
\lf[
l_i
\lf(\D{U^h_n}{x^m}+U^h_sN^s{}_{mn}+ L^h{}_{sm}U^s_n\rg)
+
l_n\lf(\D{U^h_i}{x^m}+U^h_sN^s{}_{mi}+ L^h{}_{sm}U^s_i\rg)
\rg]
 $$

\ses

\ses

$$
+
y^k_hF^H
\fr1{F^2}
 H^2h_{ni}
\lf( \D{ U^h}{x^m}+ L^h{}_{sm}U^s_m\rg)
+
2
\fr1F
l^k
H_m
 h_{ni}.
  $$

\ses

Use here the equality
$$
\cD_m h_{ni}
=
-\fr 2HH_m
h_{ni}
$$
(see (II.3.21)).

{

The rest is
\ses\\
$$
N^k{}_{mni}
=
\fr2HH_m
\fr1Fl^k
h_{ni}
-
\cD_m C^k{}_{ni}
-
\fr1{F^2}
l_i
l_nh^k_s
 N^s{}_m
+
\fr1{F^2}
\lf(
h_{ni}h^k_s
-
l_nl_sh^k_i-l_il_sh^k_n
\rg)
 N^s{}_m
$$

\ses

\ses

$$
    -
 \lf[
\fr1{F^2}
S_n{}^k{}_{si}
-\fr1{F}
\lf(l_iC^k{}_{ns}+l_nC^k{}_{is}\rg)
               \rg]
 N^s{}_m
-
(1-H^2)\fr1{F^2}h^k_i
h_{ns}
 N^s{}_m
-
\fr1{F^2}l_i
l_nh^k_s
 N^s{}_m
$$

\ses

\ses

$$
-H\fr1{F^2}l^k
(l_ih_{ns}+l_nh_{is})
 N^s{}_m
-
y^k_hF^H
\fr1{F}
\lf(
l_i
U^h_{ns}
+
l_nU^h_{is}
\rg)
N^s{}_m
-
y^k_hF^H
\fr1{F^2}
 H^2h_{ni}
U^h_sN^s{}_m
  $$

\ses

\ses

\ses

\ses

$$
=
\fr2HH_m
\fr1Fl^k
h_{ni}
-
\cD_m C^k{}_{ni}
-
\fr1{F^2}
l_i
l_nh^k_s
 N^s{}_m
-
\fr1{F^2}
\lf(
l_nl_sh^k_i+  l_il_sh^k_n
\rg)
 N^s{}_m
$$

\ses

\ses

$$
    -
 \lf[
\fr1{F^2}
S_n{}^k{}_{si}
-\fr1{F}
\lf(l_iC^k{}_{ns}+l_nC^k{}_{is}\rg)
               \rg]
 N^s{}_m
-
(1-H^2)\fr1{F^2}
\lf(h^k_ih_{ns}-h_{in}h^k_s\rg)
 N^s{}_m
$$

\ses

\ses

$$
-
\fr1{F^2}l_i
l_nh^k_s
 N^s{}_m
-H\fr1{F^2}l^k
(l_ih_{ns}+l_nh_{is})
 N^s{}_m
$$

\ses

\ses

 $$
-
\fr1{F}
\lf[
l_i
 \lf(
       C^k{}_{sn}
-\fr1{F}
\lf(l_sh^k_n+l_nh^k_s\rg)
-
\fr H{F}
h_{ns}
l^k
\rg)
+
l_n
 \lf(
       C^k{}_{si}
-\fr1{F}
\lf(l_sh^k_i+l_ih^k_s\rg)
-
\fr H{F}
h_{is}
l^k
\rg)
\rg]
N^s{}_m.
 $$
Since the indicatrix curvature equals $H^2$, we have
$
S_n{}^k{}_{si}
=
-
(1-H^2)
\lf(h^k_ih_{ns}-h_{in}h^k_s\rg).
$

\ses

The eventual result is the representation
\be
N^k{}_{mni}
=
\fr2HH_m
\fr1Fl^k
h_{ni}
-
\cD_m C^k{}_{ni}.
\ee

Thus,  Proposition II.3.5   is valid.

{

\setcounter{equation}{0}

\ses

\ses

\nin
{ \bf Appendix E:  Implications from angle}

\ses

\ses

Below, the consideration refers
to an {\it arbitrary} Finsler space.
No assumptions concerning the curvature of indicatrix
are made.
We use the  angle $\al=\al_{\{x\}}(y_1,y_2)$
which is
the
geodesic-arc distance on
the indicatrix,
in accordance with the initial definition (I.1.1).

{

In terms of  the function
\be
E~:=\fr12\al^2
\ee
the preservation equation
$
d_i \al+ (1/H)H_i\al=0
$
proposed by (I.1.28) reads
\be
\D{E}{x^i}  +N^k{}_{1i}   \D{E}{y^k_1} +   N^k{}_{2i}   \D{E}{y^k_2}
=
-
\fr 2HH_i
 E,
\ee
\ses\\
where
$
N^k{}_{1i}  =N^k{}_{i}(x,y_1), \,
N^k{}_{2i}  =N^k{}_{i}(x,y_2),
$
and
$
H_i=\partial H/\partial x^i.
$
Henceforth,
\be
H=H(x).
\ee

\ses

We are aimed to extract  the required coincidence limits
from the function $E$,
treating
 the indicatrix naturally to be a particular Riemannian space metricized by the help of
the metric tensor induced by  the Finsler metric tensor.

{

Let a  set of scalars
$
u^a=u^a(x,y)
$
be used to coordinatize the indicatrices;
the indices $a,b,...$ will be specified over the range (1,2,...N-1).
We shall use the derivative objects
$$
u^a_m=  \D{u^a}{y^m},
\qquad
u^a_{mk}=  \D{u^a_m}{y^k}.
$$
The scalars are assumed to be positively homogeneous of degree zero with respect to
the variable $y$:
\be
u^a(x,ky)=u^a(x,y), \qquad k>0, \quad \forall y,
\ee
which directly entails the identities
\be
u^a_my^m=0,
\qquad
u^a_{mk}y^k= -u^a_m.
\ee

\ses

Using the parametrical representation $l^i=t^i(u^a)$
of the indicatrix, where $l^i$ are unit vectors (possessing the property
$F(l)=1$),
 we can construct the induced metric tensor
\be
i_{ab}(u^c)=g_{mn}t^m_at^n_b\equiv
h_{mn}t^m_at^n_b
\ee
on the  indicatrix
by the help of the projection factors
$
t^m_a=\partial t^m/\partial u^a
$
(the method was described in detail in Section
5.8) of [1]).

The validity of the equalities
\be
Fu^b_mt^m_c=\de^b_c, \quad
Fu^b_{m}
t^k_b=h^k_m,   \quad
Fu^c_k=g_{km}t^m_ai^{ac},   \quad
t^n_e   i^{ec}t^i_c=h^{ni}, \quad
\fr1{F^2}h_{mn}=i_{ab}u^a_mu^b_n
\ee
can readily be verified.

{

From the identity $l_mt^m_a=0$ it follows that
\be
l_mt^m_{ab}=-i_{ab},
\ee
where $t^m_{ab}=\partial t^m_{a}/\partial u^b$.

From (E.6) we get

$$
i_{ab,c}=2FC_{mnk}t^m_at^n_bt^k_c
+
g_{nm}
(t^m_{ac}t^n_b+t^m_at^n_{bc}).
$$
With the coefficients  $i_{ac,b}=\partial i_{ac}/\partial u^b$
we obtain
$$
(i_{ac,b}+i_{bc,a}-i_{ab,c})=2FC_{mnk}t^m_at^n_bt^k_c
+
2g_{nm}
t^m_{ab} t^n_c,
$$
which entails
$$
(i_{ae,b}+i_{be,a}-i_{ab,e})i^{ec}=2FC_{mnk}t^m_at^n_bt^k_e   i^{ec}
+
2g_{nm}
t^m_{ab} t^n_e   i^{ec},
$$
so that
$$
t^i_c
\Bigl(i^c{}_{ab}-FC_{mnk}t^m_at^n_bt^k_e  i^{ec}  \Bigr)
=
h^i_m
t^m_{ab}
$$
and
\be
t^i_{ab}=
t^i_c
\Bigl(i^c{}_{ab}-FC_{mnk}t^m_at^n_bt^k_e  i^{ec}  \Bigr)
-l^ii_{ab}
\ee
(this equation is equivalent to (5.8.8) of [1]).

{

The  indicatrix Christoffel symbols
$$
i^c{}_{ab}
=\fr12i^{ce}\lf(\D{i_{ea}}{u^b}+\D{i_{eb}}{u^a}-\D{i_{ab}}{u^e}\rg)
$$
and the indicatrix curvature tensor
\be
I_a{}^e{}_{bd} ~:=
\D{i^e{}_{ab}}{u^d}-
\D{i^e{}_{ad}}{u^b}
+    i^f{}_{ab}   i^e{}_{fd}
-
i^f{}_{ad}   i^e{}_{fb}
\ee
will be used.

\ses

Constructing the tensor
\be
S_{abcd} =   -\fr13(I_{acbd}+I_{adbc}),
\ee
where $I_{acbd}=I_a{}^e{}_{bd}i_{ec},$
we obtain
the useful  identity
\be
S_{abcd} -S_{acbd}=  -I_{adbc}.
\ee

{

It follows that
$$
i_{ab}\lf(u^b_{mj}+i^b{}_{cv}u^v_mu^c_j\rg)
=
h_{pq}
\lf(t^p_at^q_bu^b_{mj}+i^b{}_{cv}t^p_at^q_bu^v_mu^c_j\rg)
=
h_{pq}
t^p_a
\lf(
 \D{\fr1Fh^q_m}{y^j}- u^b_mt^q_{bc}u^c_j
+
i^b{}_{cv}t^q_bu^v_mu^c_j\rg).
$$
Taking  $t^p_{ab}$  from (E.9), we get
\ses\\
$$
i_{ab}\lf(u^b_{mj}+i^b{}_{cv}u^v_mu^c_j\rg)
=
h_{pq}
t^p_a
\lf(
 \D{\fr1Fh^q_m}{y^j}- u^b_m
t^q_f \Bigl(i^f{}_{bc}-FC_{rsk}t^r_bt^s_ct^k_e  i^{ef}  \Bigr)u^c_j
+
i^b{}_{cv}t^q_bu^v_mu^c_j\rg).
$$
\ses\\
Here,
$$
t^k_e  i^{ef}=Fu^f_hg^{hk},
$$
so that
\be
i_{ab}\lf(u^b_{mj}+i^b{}_{cv}u^v_mu^c_j\rg)
=
h_{pq}
t^p_a
\lf(
 \D{\fr1Fh^q_m}{y^j}
 +
\fr1F
 C_{mj}{}^q
\rg),
\ee
\ses
from which it follows that
\be
u^b_{mj}+i^b{}_{cv}u^v_mu^c_j
=
Fu^b_t
\lf(
 \D{\fr1Fh^t_m}{y^j}
 +
\fr1F
 C_{mj}{}^t
\rg)
\ee
\ses\\
and
\be
u^a_{mk}u^b_ni_{ab}+u^a_mu^f_ku^b_ni^c{}_{af}i_{cb}
=
-
\fr1{F^3}(h_{nk}l_m+h_{nm}l_k)
+
\fr1{F^2}
  C_{kmn}.
  \ee

{

Now we consider the quantity (E.1) on the indicatrix:
\be
E=M(x,u_1,u_2),
\ee
where $M$ is a scalar.

\ses

There arise the objects
\be
M_{1a2c} \eqdef
\D{M_{1a}}{u^c_2}, \qquad
M_{1a2c2d} \eqdef
\D{M_{1a2c}}{u^d_2}
-
i^{2f}{}_{2d2c}M_{1a2f},
\ee
\ses
together with
\ses\\
\be
M_{1a1c}
\eqdef
\D{M_{1a}}{u^c_1}
-
i^{1f}{}_{1a1c}M_{1f}, \qquad
M_{1a1b2c} \eqdef
\D{M_{1a1b}}{u^c_2}
=
{\displaystyle\frac{\partial^3{M}}
{\partial{u^a_1  }    \partial{u^b_1}    \partial{u^c_2} }}
 -
i^{1f}{}_{1a1b} M_{1f2c}
\ee

\ses

\be
M_{1a1b2c2d} \eqdef
\D{M_{1a1b2c}}{u^d_2}
-
i^{2f}{}_{2d2c}M_{1a1b2f}
=
{\displaystyle\frac{\partial^4{M}}
{\partial{u^a_1  }    \partial{u^b_1}    \partial{u^c_2}  \partial{u^d_2}}}
-
 i^{1f}{}_{1a1b} \D {M_{1f2c}        }{u^d_2}
-
i^{2f}{}_{2d2c}M_{1a1b2f}.
\ee

It follows that

\be
{\displaystyle\frac{\partial^4{M}}
{\partial{u^a_1  }    \partial{u^b_1}    \partial{u^c_2}  \partial{u^d_2}}}
=
M_{1a1b2c2d}
+
 i^{1f}{}_{1a1b} \D {M_{1f2c}        }{u^d_2}
 +
i^{2f}{}_{2d2c}M_{1a1b2f}.
\ee

{

In the limit
$
u_2\to u_1
$
we have
\ses\\
\be
\D{M}{u^a_1}\to 0 ,
\quad
\D{M}{u^a_2}\to 0,
\ee
\ses
 and
 \be
\Dd{M}{u^a_1}{u^b_1}\to i_{ab},
 \qquad \Dd{M}{u^a_1}{u^b_2} \to - i_{ab},
  \qquad \Dd{M}{u^a_2}{u^b_2} \to i_{ab},
\ee
\ses
together with
\be
 {\displaystyle\frac{\partial^3{M}}
{\partial{u_1^a}    \partial{u_2^b}    \partial{u_2^c}}}
\to
-
  i^e{}_{bc}i_{ea}, \qquad
 {\displaystyle\frac{\partial^3{M}}
{\partial{u_1^a}    \partial{u_1^b}    \partial{u_2^c}}}
\to
 - i^f{}_{ab}i_{fc}
 \ee
(see Section 3.2 in [12]).

Also,
\be
M_{1a1b2c2d}  \to S_{abcd}
\ee
(see (3.2.69) in [12]).
From (E.20) it follows that
\ses\\
\be
{\displaystyle\frac{\partial^4{M}}
{\partial{u^a_1  }    \partial{u^b_1}    \partial{u^c_2}  \partial{u^d_2}}}
\to
S_{abcd}
 -
 i^{f}{}_{ab}
 i^{e}{}_{dc}i_{fe}.
\ee
{

We find
\ses\\
\be
\D{E}{y_2^k}=u^d_{2k}\D{M}{u^d_2},\qquad
\Dd{E}{y_2^k}{y^n_2}=u^d_{2k2n}\D{M}{u^d_2}
+
u^d_{2k}\Dd{M}{u^d_2}{u^t_2}u^t_{2n}
\ee
\ses\\
and
\be
   {\displaystyle\frac{\partial^3{E}}
 {  \partial {y_2^k} \partial {y^n_2} \partial{y^m_1}}}
=
u^d_{2k2n}\Dd{M}{u^d_2}{u_1^v}u^v_{1m}
+
u^d_{2k}
{\displaystyle\frac{\partial^3{M}} {\partial{u^d_2} \partial{u^t_2}  \partial{u^s_1}} }
u^t_{2n}u^s_{1m}.
\ee

\ses

\ses

Moreover,
\ses\\
$$
   {\displaystyle\frac{\partial^4{E}}
 {  \partial {y_2^k} \partial {y^n_2} \partial{y^m_1}  \partial{y^j_1}}}
=
u^d_{2k2n}
{\displaystyle\frac{\partial^3{M}} {\partial{u^d_2} \partial{u^v_1}  \partial{u^s_1}} }
u^v_{1m}
u^s_{1j}
+
u^d_{2k2n}\Dd{M}{u^d_2}{u_1^v}u^v_{1m1j}
$$

\ses

\ses

\ses

\be
+
u^d_{2k}
{\displaystyle\frac{\partial^3{M}} {\partial{u^d_2} \partial{u^t_2}  \partial{u^s_1}} }
u^t_{2n}u^s_{1m1j}
+
u^d_{2k}
{\displaystyle\frac{\partial^4{M}} {\partial{u^d_2} \partial{u^t_2}
 \partial{u^s_1} \partial{u^t_1}} }
u^t_{2n}u^s_{1m}u^t_{1j}.
\ee

{

These observations
entail
 the limits
\ses\\
\be
y_2\to y_1 ~:  \quad \D{E}{y^m_1} \to 0,  \quad \D{E}{y^m_2} \to 0,
\ee

\ses

\ses

\be
\Dd{E}{y^m_1}{y^n_1} \to \fr1{F^2} h_{mn},
\quad \Dd{E}{y^m_2}{y^n_2} \to \fr1{F^2} h_{mn},
\quad \Dd{E}{y^m_1}{y^n_2} \to -\fr1{F^2} h_{mn},
\ee
\ses\\
and
\be
{\displaystyle\frac{\partial^3{E}}
{\partial{y_2^k}    \partial{y_1^m}    \partial{y_2^n}}}
\to
-
\bigl(u^a_{nk}u^b_mi_{ab}+u^a_mu^b_ku^c_ni^e{}_{bc}i_{ea}\bigr),
\ee
\ses\\
plus
\be
{\displaystyle\frac{\partial^3{E}}
{\partial{y_1^k}    \partial{y_1^m}    \partial{y_2^n}}}
\to
-
\bigl(u^a_{mk}u^b_ni_{ab}
+
u^a_nu^b_ku^c_mi^e{}_{bc}i_{ea}
\bigr),
\ee
\ses\\
together with
$$
{\displaystyle\frac{\partial^4{E}}
 {  \partial {y_2^k} \partial {y^n_2} \partial{y^m_1}  \partial{y^j_1}}}
 \to
-
u^a_{kn}
\lf(i_{ab}u^b_{mj}+i^f{}_{cv}i_{af}u^v_mu^c_j\rg)
-
u^a_{mj}
i^f{}_{cv}i_{af}u^v_ku^c_n
$$

\ses

\ses

$$
+
\lf(
S_{abcd}
 -
 i^{1f}{}_{1a1b}
 i^{e}{}_{dc}i_{fe}
\rg)
u^a_ku^b_nu^c_mu^d_j,
$$
\ses
or
\ses\\
\be
{\displaystyle\frac{\partial^4{E}}
 {  \partial {y_2^k} \partial {y^n_2} \partial{y^m_1}  \partial{y^j_1}}}
 \to
-
\lf( u^a_{kn}  +i^a{}_{cv}u^v_ku^c_n\rg)
i_{ab}
\lf(u^b_{mj}+i^b{}_{cv}u^v_mu^c_j\rg)
+
 S_{abcd}
u^a_ku^b_nu^c_mu^d_j.
\ee

{

In this way we arrive at the  reductions
 \ses\\
 \be
{\displaystyle\frac{\partial^3{E}}
{\partial{y_2^k}    \partial{y_1^m}    \partial{y_2^n}}}
\to
-
\bigl(-\fr1{F^3}(h_{nm}l_k+h_{nm}l_k)+\fr1{F^2}  C_{kmn}\bigr)
\ee
\ses
and
\ses\\
\be
{\displaystyle\frac{\partial^3{E}}
{\partial{y_2^k}    \partial{y_1^m}    \partial{y_1^n}}}
\to
-
\bigl(-\fr1{F^3}(h_{nk}l_m+h_{mk}l_n)+\fr1{F^2}  C_{kmn}\bigr).
\ee

\ses

Taking into account (E.13),
we can write
\ses\\
$$
{\displaystyle\frac{\partial^4{E}}
 {  \partial {y_2^k} \partial {y^n_2} \partial{y^m_1}  \partial{y^j_1}}}
 \to
-
\lf(
 \D{\fr1Fh^t_k}{y^n}
 +
\fr1F
 C_{kn}{}^t
\rg)
h_{tq}
\lf(
 \D{\fr1Fh^q_m}{y^j}
 +
\fr1F
 C_{mj}{}^q
\rg)
+
S_{abcd}
u^a_ku^b_nu^c_mu^d_j,
$$
\ses
which  is
\ses\\
$$
F^4
{\displaystyle\frac{\partial^4{E}}
 {  \partial {y_2^k} \partial {y^n_2} \partial{y^m_1}  \partial{y^j_1}}}
 \to
\lf(
l_nh^t_k+l_kh^t_n - F C_{kn}{}^t\rg)
\lf(
-l_jh_{tm}-l_mh_{tj} + F C_{mjt}\rg)
+
F^4
S_{abcd}
u^a_ku^b_nu^c_mu^d_j,
$$
\ses\\
or
$$
F^4
{\displaystyle\frac{\partial^4{E}}
 {  \partial {y_2^k} \partial {y^n_2} \partial{y^m_1}  \partial{y^j_1}}}
 \to
l_n
\lf(
-l_jh_{km}-l_mh_{kj} + F C_{mjk}\rg)
+
l_k
\lf(-l_jh_{nm}-l_mh_{nj} + F C_{mjn}\rg)
$$

\ses

\ses

\be
+F C_{knm}l_j  +F C_{knj}l_m
 -   F^2 C_{kn}{}^t   C_{mjt}
+
F^4
 R _{abcd}
u^a_ku^b_nu^c_mu^d_j.
\ee

{

Now we differentiate the  preservation law  (E.2)
 with respect to $y_1^m$ and
 $y_2^n$ and make $y_2\to y_1$,
 which yields
 \ses\\
$$
\partial_i h_{mn}
-
 h_{mn}
 \fr1{F^2}
\partial_i F^2
+
N^k{}_{im}h_{kn}
+
N^k{}_{in}h_{km}
$$

\ses

\ses

$$
+
\fr1{F}
N^k{}_{i}
\Bigl(-(h_{nk}l_m+h_{mn}l_k)+ FC_{kmn}\Bigr)
+
\fr1{F}
N^k{}_{i}
\Bigl(-(h_{km}l_m+h_{nm}l_k)+F  C_{kmn}\Bigr)
=
-\fr2HH_i h_{mn}.
$$
\ses\\
On so dong we come to the following sought equality
\be
\cD_ih_{mn}
-\fr 2{F}
   h_{mn}
   d_iF
=
-\fr2HH_i h_{mn}
\ee
\ses
with
\be
\cD_ih_{mn}=
\partial_i h_{mn}
+
N^k{}_{i}  \D{h_{mn}}{y^k}
+
N^k{}_{im}h_{kn}
+
N^k{}_{in}h_{km}
\ee
and
\ses\\
\be
d_iF=\partial_i F+N^k{}_il_k.
\ee

{

In the vanishing case
\be
\partial_i F+N^k{}_il_k=0
\ee
we obtain by differentiation the equalities
\be
\partial_i l_n  +\fr1FN^k{}_ih_{kn}
+N^k{}_{in}l_k=0
\ee
\ses\\
and
$$
\partial_i h_{mn}
 +
 N^k{}_{i}\D{h_{kn}}{y^m} +N^k{}_{im}h_{kn}
 +N^k{}_{in}h_{km}
-\fr1F N^k{}_{i}h_{kn}l_m
+ \fr1F N^k{}_{i}l_kh_{mn}
+
FN^k{}_{inm}l_k
=0.
$$
So we can write
\be
\cD_ih_{mn}
+
FN^k{}_{inm}l_k
=0
\ee
together with
\ses\\
\be
FN^k{}_{inm}l_k=
\fr2HH_i h_{mn}.
\ee

{

\ses

\ses

\def\bibit[#1]#2\par{\rm\noindent\parskip1pt
                     \parbox[t]{.05\textwidth}{\mbox{}\hfill[#1]}\hfill
                     \parbox[t]{.925\textwidth}{\baselineskip11pt#2}\par}

\bc {\bf  References}
\ec

\ses

\bibit[1] H. Rund, \it The Differential Geometry of Finsler
 Spaces, \rm Springer, Berlin 1959.

\bibit[2] D. Bao, S. S. Chern, and Z. Shen, {\it  An
Introduction to Riemann-Finsler Geometry,}  Springer, N.Y., Berlin 2000.

\bibit[3] L. Kozma and L. Tam{\'a}ssy,
Finsler geometry without line elements faced to applications,
{   \it Rep. Math. Phys.} {\bf 51} (2003), 233--250.

\bibit[4] L. Tam{\'a}ssy,
 Metrical almost linear connections in $TM$ for Randers spaces,
{\it Bull. Soc. Sci. Lett. Lodz Ser. Rech. Deform } {\bf 51} (2006), 147-152.

\bibit[5] Z. L. Szab{\'o},  All regular Landsberg metrics are Berwald,
{\it Ann Glob Anal Geom } {\bf 34} (2008), 381-386.

\bibit[6] L. Tam{\'a}ssy,  Angle in Minkowski and Finsler spaces,
{\it Bull. Soc. Sci. Lett. Lodz Ser. Rech. Deform } {\bf 49} (2006), 7-14.

\bibit[7] G. S. Asanov,
 Finsleroid  gives rise to the angle-preserving  connection,
 {\it  arXiv:} 0910.0935 [math.DG],  (2009).

\bibit[8] G. S. Asanov,
 Finslerian angle-preserving  connection in two-dimensional
case.  Regular realization,
 {\it  arXiv:} 0909.1641 [math.DG] (2009).

\bibit[9] G. S. Asanov,
 Finsler space connected by angle  in two dimensions. Regular case,
{\it Publ. Math. Debrecen } {\bf 77/1-2} (2010), 245--259.

\bibit[10] G. S. Asanov,
 Finsler connection preserving angle in dimensions $N\ge3$,
  {\it  arXiv:} 1009.1215 [math.DG] (2010).

\bibit[11]
G. S. Asanov,
 Finsler angle-preserving  connection in dimensions $N\ge3$,
 {\it Publ. Math. Debrecen } {\bf 79/1-2} (2011), 181--209.

\bibit[12] J. L. Synge, {\it Relativity: the General Theory,}
 North-Holland, Amsterdam, 1960.

\ses

\end {document}